\documentclass{article}
\usepackage[pagebackref=true]{hyperref}
\hypersetup{colorlinks=true,linkcolor=blue,citecolor=red}
\usepackage[a4paper,margin=1in]{geometry} 
\usepackage{amssymb}
\usepackage{amsmath}
\usepackage{mathabx} 

\newtheorem{theorem}{Theorem}
\newtheorem{definition}{Definition}
\newtheorem{lemma}[theorem]{Lemma}
\newtheorem{rem}[theorem]{Remark}

\newtheorem{cor}[theorem]{Corollary}

\newcommand{\Id }{\mathrm{Id} \,}

\newcommand{\abs}[1]{\lvert#1\rvert} 
\newcommand{\Res}{\mathcal{R}\mathrm{es}}
\newcommand{\HRes}{\mathrm{Res}}

\newcommand{\PhSpace}{\mathcal{P}} 

\newcommand{\SSym}{\mathcal{S}}  

\newcommand{\comment}[1]{\mbox{}}

\def\qed{{\hfill{\vrule height5pt width3pt depth0pt}\medskip}}

\newcommand{\com}[1]{%
  \ifnum\visibleComment=0
  \else \textbf{#1}
  \fi
}
\newcommand{\visibleComment}{0}

\begin{document}
\begin{center}
{\LARGE \textbf{Polynomial normal forms for ODEs near a center-saddle equilibrium point}} \\
 \vskip 0.5cm

{\large Amadeu Delshams}\footnote{Partially supported by the Spanish MINECO/FEDER Grant PID2021-123968NB-I00}\\
Lab of Geometry and Dynamical Systems and IMTech, Universitat Polit{\`e}cnica de Catalunya (UPC) and Centre de Recerca Matem{\`a}tica (CRM),
Barcelona, Spain\\
\texttt{Amadeu.Delshams@upc.edu}

 \vskip\baselineskip
   {\large Piotr Zgliczy\'nski}\footnote{Research has been supported by polish NCN grant
  2016/22/A/ST1/00077 
  }     \\
 Jagiellonian University, Institute of Computer Science and Computational Mathematics, \\
{\L}ojasiewicza 6, 30--348  Krak\'ow, Poland \\
\texttt{umzglicz@cyf-kr.edu.pl}

\vskip 0.5cm
 \today

\end{center}

\begin{abstract}
In this work we consider a saddle-center equilibrium for general vector fields as well as Hamiltonian systems, and we
transform it locally into a polynomial normal form in the saddle variables by a change of coordinates.
This problem was first solved by Bronstein and Kopanskii in 1995, as well as by Banyaga, de la Llave and Wayne in 1996 [BLW] in the saddle case.
The proof used relies on the deformation method used in [BLW], which in particular implies the preservation of the symplectic form for a Hamiltonian system, although our proof is different and, we believe, simpler.
We also show that if the system has sign-symmetry, then the transformation can be chosen so that it also has sign-symmetry.
This issue is important in our study of shadowing non-transverse heteroclinic chains (Delshams and Zgliczynski 2018 and 2024) for the toy model systems (TMS) of the cubic defocusing nonlinear Schr\"odinger equation (NLSE) on $2D$-torus  or similar Hamiltonian PDE, which are used to prove energy transfer in these PDE.
\end{abstract}



\section{Introduction}
The main question addressed in this work is:
Can a vector field near a center-saddle equilibrium be transformed into a polynomial normal form in the saddle variables by a change of coordinates? (finite smoothness is allowed). This problem has been solved by Bronstein and Kopanskii~\cite{BK95,BK96}, as well as by Banyaga, de la Llave and Wayne~\cite{BLW} (they give a complete proof for maps, but the result also applies to ODEs) for a saddle equilibrium point, where it is also shown that the coordinate transformation can be chosen preserving the symplectic form, volume or contact if the original vector field preserves it.


In this work we consider a saddle-center equilibrium for general vector fields as well as Hamiltonian systems, and we apply the deformation method as in~\cite{BLW}. These systems arose in our work on shadowing non-transversal heteroclinic chains~\cite{DSZ,DZ}. Although the technique of the proof is inspired by the approach taken in~\cite{BLW}, our proof is different and, we believe, simpler. We also show that if the system has sign-symmetry, then the transformation can be chosen so that it also has sign-symmetry.
This issue is important in our study of shadowing non-transverse heteroclinic chains~\cite{DSZ,DZ} for the toy model systems (TMS) of the nonlinear cubic defocusing   Schr\"odinger equation on $2$-dimensional torus  (NLSE) or similar Hamiltonian PDE, which are used to prove energy transfer in these PDE.



The plan of the paper is as follows. We first prove the result on general vector fields and Hamiltonian systems, which is stated as Theorem~\ref{thm:polNForm-loccenter} in Section~\ref{sec:PNF}. This proof is concluded
in section~\ref{sec:estm-gen-coheq}. Later on we introduce some additional invariance requirements (the sign-symmetry) and we show in Theorem~\ref{thm:polNForm-centerSubspaces} that the transformation constructed in the general case also has these properties.

\subsection{Notation}


We will denote $\mathbb{N}$, $\mathbb{Z}$, $\mathbb{R}$, $\mathbb{C}$ the natural, integer, real and complex numbers, respectively. We will assume that $0 \in \mathbb{N}$.
For $x \in \mathbb{R}$, we denote $[x]$ its integer part, $\Re x$ its real part and $\Im x$ its imaginary part.
For a matrix $A$ we will denote $A^\top$ its transposed matrix. Finally we will denote $B_n(R)$ the ball of radius $R$ centered at the origin of $\mathbb{R}^n$.

If $u=h(z)$ is a (local) diffeomorphism in $\mathbb{R}^n$, expressing new coordinates $u$ in terms of the old ones denoted by $z$, and $\mathcal{ Z} (z)$ is a vector field, then it is transformed by $u=h(z)$ to the pushforward vector field $h_{\ast}\mathcal{Z}$ defined as
\begin{equation*}
h_{\ast} \mathcal{Z}(u):= Dh(h^{-1}(u))\mathcal{Z}(h^{-1}(u)).
\end{equation*}

Conversely, given $z=g(u)$ (for example $g=h^{-1}$) expressing the old coordinates $z$ in terms of the new coordinates $u$, $\mathcal{ Z} (z)$  is transformed by $z=g(u)$ to the pullback vector field $g^\ast \mathcal{Z} $ defined as
\begin{equation}
  g^\ast \mathcal{Z} (u)=Dg(u)^{-1} \mathcal{Z}(g(u)).  \label{eq:hast}
\end{equation}

The action on functions is as follows
\begin{eqnarray*}
  h_\ast Z (u)&=&Z(h^{-1}(u)), \\ 
  g^\ast Z (u)&=&Z(g(u)).
\end{eqnarray*}

The commutator of two vector fields $X,Y$ is given by
\begin{eqnarray}
  {}[X,Y]_i = \sum_j \left(\frac{\partial Y_i}{\partial z_j} X_j - \frac{\partial X_i}{\partial z_j} Y_j \right), \nonumber\\
  {}[X,Y](z)= DY(z) X(z) - DX(z) Y(z).  \label{eq:commutator}
\end{eqnarray}


For a function $g(\varepsilon,x)$ depending on $\varepsilon$, which is treated as a parameter, we will use  the following notations
\begin{eqnarray*}
  g_\varepsilon(x) &=& g(\varepsilon,x), \\
  Dg_\varepsilon(x) &=& \frac{\partial g}{\partial x}(\varepsilon,x).
\end{eqnarray*}


For $g(z) \in C^r$ we will say that $g(z)=O_k(|z|^q)$ (for $z \to 0$---this might not be written explicitly) for some $k \leq r$ if
\begin{eqnarray*}
  D^j g(z) = O(|z|^{q-j}), \quad j=0,\dots,\min(k,q).
\end{eqnarray*}
If $z=(x,y)$ and $g(z) \in C^r$, we will say $g(z)=O_k(|x|^q)$ for  $z \in W$ (bounded set) and $x \to 0$ for  $k \leq r$ if
\begin{eqnarray*}
  D^j g(z) = O(|x|^{q-j}), \quad j=0,\dots,\min(k,q),
\end{eqnarray*}
where $D$ is the differentiation with respect to $z=(x,y)$, i.e., all partial derivatives are involved.

Given a vector space $V$ with the basis $\{e_i\}_{i=1}^m$ and a multiindex
$\alpha=(\alpha_1,\alpha_2,\allowbreak \dots,\allowbreak \alpha_m) \in \mathbb{N}^m$, for a vector $z = \sum_i z_i e_i$
we define $z^\alpha$ by
\begin{equation*}
  z^\alpha=z_1^{\alpha_1} z_2^{\alpha_2} \cdots z_m^{\alpha_m}.
\end{equation*}
By $|\alpha|$ we will denote the degree of multiindex, defined by $|\alpha|=\sum_i \alpha_i$.

Often, when discussing normal forms we will use complex coordinates $c_k=x_k + i y_k$, and then a part of the phase space will be described by $\mathbb{C}^r$
and the polynomials on $\mathbb{R}^{m} \times \mathbb{C}^{r}$ will be of the following  form
\begin{equation}
  z^\alpha c^\beta {\overline{c}}^{\beta^*}= z_1^{\alpha_1} z_2^{\alpha_2} \cdots z_{m}^{\alpha_{m}} \cdot c_1^{\beta_1} \cdots c_{r}^{\beta_{r}} \cdot {\overline{c}_1}^{\beta^*_1} \cdots {\overline{c}_{r}}^{\beta^*_{r}}. \label{eq:zalpha-compl}
\end{equation}

In such situation the degree of the monomial $z^\alpha c^\beta {\overline{c}}^{\beta^*}$ is given  by $|\alpha|+|\beta|+|\beta^*|$.

\subsubsection{Symplectic forms and Hamiltonian equations}

By $\Omega$ we will denote a symplectic form (a closed non-degenerate differential 2-form) on some phase space $V$.
A map is called symplectic (with respect to $\Omega$) iff it preserves the form $\Omega$.

 The standard (real) symplectic form is given by
\begin{equation*}
\Omega=\sum_{j}  dx_j \wedge dy_j. 
\end{equation*}
Occasionally we will use more general symplectic form, for example $\Omega'=g^\ast \Omega$. 

Given a real function $H(t,z)$ for real $t$ and $z=(x,y)$, and a real symplectic form $\Omega$ the Hamiltonian vector field $X_H(t,z)$ is defined by \cite{A}
\begin{equation}
  \Omega(X_H(t,z),\eta)=D_zH(t,z)\eta, \quad  \forall \eta \in T_z V.  \label{eq:def-ham-vf}
\end{equation}
Therefore the Hamiltonian equations  can be written as
\begin{equation}
  \dot{z} =J(z)  \nabla_{z} H(t,z) \label{eq:def-ham-vf-J}
\end{equation}
where $\nabla_{z} H(t,z)=\left(D_{z} H(t,z)\right)^\top$ and $J(z)$ is the regular antisymmetric matrix defined by (\ref{eq:def-ham-vf}).

In the case of the standard symplectic form  the Hamiltonian equations   of motion are just
\begin{equation}
   \dot{x}_j = \frac{\partial H}{\partial y_j}, \quad \dot{y}_j = - \frac{\partial H}{\partial x_j}. \label{eq:ham-eq}
\end{equation}

Since we are going to consider systems $\dot z=Az +O_2(z)$ in a neighborhood of an equilibrium point, it will be very convenient to consider complex variables
$c_l=x_l + i y_l$ for some $l\in L$ associated to pure imaginary eigenvalues $\nu=\pm i \omega$, $\omega\neq  0$, of the matrix $A$. Then the standard symplectic form $\Omega$ becomes
\begin{equation*}
  \Omega = \sum_{k \in K} dx_k \wedge d y_k + \frac{i}{2}
  \sum_{l \in L} dc_l \wedge d\overline{c}_l. 
\end{equation*}
In such case the Hamiltonian equations of motion (\ref{eq:ham-eq}) become
\begin{eqnarray*}
  \dot{x}_k &=& \frac{\partial H}{\partial y_k}, \quad \dot{y}_k =  -\frac{\partial H}{\partial x_k}, \quad k \in K, \\
  \dot{c}_l &=& -2i \frac{\partial H}{\partial \overline{c}_l}, \quad l \in L.
\end{eqnarray*}

One can get a more compact notation introducing the  complex variables $c_k=x_k + i y_k$ for the indexes $k\in K$ associated to real eigenvalues to get a canonical complex Hamiltonian system
\[
  \dot{c}_j = -2i \frac{\partial H}{\partial \overline{c}_j}
\]
associated to the standard complex symplectic form
$\displaystyle \Omega =\frac{i}{2}\sum_{j} dc_j \wedge d\overline{c}_j$.

In addition, in the presence of saddle variables associated to complex eigenvalues $\nu=\pm a \pm i b$, $ab\neq 0$, of the matrix $A$, complex pairs of variables are handy and are of the form $c_{lm}=x_l+ix_m$ and $g_{lm}=y_l+iy_m$. Then  in  the standard symplectic form the contribution of $dx_l \wedge dy_l + dx_m \wedge dy_m $ becomes $\frac{1}{2} dc_{lm} \wedge d\bar{g}_{lm} + \frac{1}{2} d\bar{c}_{lm} \wedge dg_{lm}$ so that $\Omega$ becomes
\begin{equation*}
   \Omega = \sum_{k \in K} dx_k \wedge d y_k + \frac{1}{2}\sum_{(l,m) \in L' \times L''} \left(dc_{lm} \wedge d\bar{g}_{lm} + d\bar{c}_{lm} \wedge dg_{lm}\right) + \frac{i}{2}
  \sum_{l \in L} dc_l \wedge d\overline{c}_l,  
\end{equation*}
and the equations of motion for the $c_{lm}$ and $g_{lm}$ variables are
\begin{equation*}
  \dot{c}_{lm}=2 \frac{\partial H}{\partial \bar{g}_{lm}}, \quad  \dot{g}_{lm}=-2 \frac{\partial H}{\partial \bar{c}_{lm}}.
\end{equation*}

For a real system associated to a real Hamiltonian, the equations of motion of $\overline{c}_l$, $\overline{c}_{lm}$ and $\overline{g}_{lm}$ are redundant and there is no need to write them, since they are conjugate to those of ${c}_l$, {$c_{lm}$, $g_{lm}$}, respectively. Indeed
\[
  \dot{\overline{c}}_l = 2i \frac{\partial H}{\partial {c}_l}, \quad
  \dot{\overline{c}}_{lm}=-2 \frac{\partial H}{\partial g_{lm}}, \quad  \dot{\overline{g}}_{lm}=-2 \frac{\partial H}{\partial c_{lm}}.
\]

\subsection{Resonant monomials}

Assume that we have a coordinate system $(z_1,\dots,z_n)$ which brings matrix $A \in \mathbb{R}^{n \times n}$ to its complex Jordan normal form. Let $\nu_1,\dots, \nu_n$ be the corresponding eigenvalues. In case some of them are not real, they appear within pairs of conjugate complex numbers, and the Jordan normal form for A takes complex values.

Let us denote by $e_j$, $j=1,\dots,n$ the $n$-dimensional canonical basis, and consider a polynomial vector field
\begin{equation}
  \dot{z}=A z + \sum_{2\leq r\leq N} \sum_{|\alpha|=r} \left(p_\alpha^{1} z^\alpha,\dots,p_\alpha^n z^\alpha \right) \label{eq:A-res-def}
\end{equation}

The following definition is taken from \cite{BK95}.
\begin{definition}\cite[Def. 2.4]{BK95}
Let $j \in \{1,\dots,n\}$ and $\alpha \in \mathbb{Z}_+^n$. The pair $(j,\alpha)$ as well as the monomial $z^\alpha e_j$ of
system~\eqref{eq:A-res-def} is called \emph{resonant} (in the sense of Poincar\'e) if the following equality holds:
\begin{equation*}
  \nu_j=\langle \alpha, \nu \rangle:=\sum_{k=1}^n \alpha_k \nu_k.
\end{equation*}
The integer $r=|\alpha|$ is called the \emph{order of resonance}. \\
We will denote by $\Res (A,k_1,k_2)$ the set of resonant monomials of equation~\eqref{eq:A-res-def} of order greater than or equal to $k_1$ and of order less than or equal to $k_2$.
\end{definition}
We will occasionally say that $z^\alpha$ is resonant for the $j$-th variable if $z^\alpha e_j$ is resonant.

In the case that system~\eqref{eq:A-res-def} is a Hamiltonian of $m=n/2$ degrees of freedom
\begin{equation}
\label{eq:A-res-def-Ham}
H=\sum_{2\leq \abs{\alpha}+ \abs{\beta}+\abs{\beta^*}\leq N+1}a_{\alpha\beta \beta^*} z^{\alpha}c^{\beta} \overline{c}^{\beta^*},
\end{equation}
where $z=(x,y)\in\mathbb{R}^{2m-2r}$, $c\in\mathbb{C}^r$
 and $a_{\alpha\beta \beta^*}=\overline{a_{\alpha\beta^* \beta}}$ for a real Hamiltonian, the eigenvalues
$(\nu_1,\dots,\nu_{m-r},\nu_{m-r+1},\dots,\nu_{m},\nu_{m+1},\dots,\nu_{2m-r},\nu_{2m-r+1},\dots,\nu_{2m})
:=(\nu_{\mathrm{R}}^+,\nu_{\mathrm{C}},\nu_{\mathrm{R}}^-,\nu_{\mathrm{C}}^*)$ of the Hamiltonian matrix $A$ appear in pairs with both signs and can be ordered, for example, so that the $2(m-r)$ complex  ones, $0\leq r\leq m$, satisfy $\nu_i=\bar{\nu}_{m+i}$ for $i=1,\dots,r$. Then the resonant monomials of Hamiltonian~\eqref{eq:A-res-def-Ham} are those that satisfy
$\langle\nu_{\mathrm{R}},\alpha\rangle+\langle\beta,\nu_{\mathrm{C}}\rangle+\langle\beta^*,\nu_{\mathrm{C}}^*\rangle=0$,
where $\nu_{\mathrm{R}}=\left(\nu_{\mathrm{R}}^+,\nu_{\mathrm{R}}^-\right)$ (see, for example,~\cite[Sec.~10.4]{MO})
We will denote by $\HRes (A,k_1,k_2)$ the set of resonant monomials of Hamiltonian~\eqref{eq:A-res-def-Ham} of order greater than or equal to $k_1$ and of order less than or equal to $k_2$.

\com{PZ: I do  not think that the above definition of resonant monomials for Hamiltonian case is correct. For real eigenvalues condition is just $\langle \beta, \nu \rangle=0$,
without any subdivision of $\nu$  into two sets $(\nu,-\nu)$ - all $\nu_j$'s are in this scalar product$\langle \beta, \nu \rangle$.  Moreover, for complex eigenvalues and diagonalizing variables we
have $\nu_j$ corresponds to variable $z_j$, then $\bar{\nu}_j$ corresponds to $\bar{z}_j$.  If $\nu_j=a+i b$, then $\bar{\nu}_j \neq - \nu_j$ if $a \neq 0$. Most likely we should
write like (\ref{eq:zalpha-compl}) in the expression for Hamiltonian (\ref{eq:A-res-def-Ham}) and carefully distinguish between real and complex variables.  }

\com{A: weakly resonant out: It is not used}


\section{The Polynomial Normal Form theorem for ge\-ne\-ral vector fields and Hamiltonian system}
\label{sec:PNF}

\subsection{Normal form around an equilibrium point, the Poincar\'e-Dulac theorem}
Assume that we have a coordinate system $(z_1,\dots,z_n)$ which brings a matrix $A \in \mathbb{R}^{n \times n}$ to its complex Jordan  normal form, and consider the system
\begin{equation}
  \dot{z}=A z + O_2(|z|^2), \label{eq:ODE-gen}
\end{equation}

One would like to bring (\ref{eq:ODE-gen}) to the simplest possible form.  The first well-known result in this direction is
the Poincar\'e-Dulac theorem~\cite{A}.
\begin{theorem}
\label{thm:Dulac-Poinc}
If the vector field (\ref{eq:ODE-gen}) is $C^{Q+1}$ with $Q \geq 2$, then for any $J\subset \mathbb{Z}_+^n$ such that $\Res(A,2,Q) \subset J$ and $2 \leq |\alpha| \leq Q$ for $\alpha\in J$,  system (\ref{eq:ODE-gen})
 after finite number of  analytic changes of variables, can be written as
\begin{equation*}
\dot z= A z + \mathcal{N}_Q(z) +\mathcal{R}_{Q+1}(z),
\end{equation*}
where $\mathcal{N}_Q(z)=\sum_{\alpha \in J}\mathcal{N}_{Q,\alpha} z^\alpha$   and the remainder vector field satisfies
$\mathcal{R}_{Q+1}(z)=O_{Q+1}(\abs{z}^{Q+1})$.
Moreover, if system~\eqref{eq:ODE-gen} is Hamiltonian, the changes of variables can be chosen symplectic, so that the transformed system is also Hamiltonian.
\end{theorem}
Regarding the notation $z^\alpha$, it should be noted that the $n$-dimensional variable $z=(x,c)$ can contain both real variables $x$ associated with real eigenvalues of the matrix $A$, as well as conjugate variables $c , \overline{c}$ associated with complex eigenvalues of the matrix $A$, so any monomial $z^\alpha$ is of the already mentioned form~\eqref{eq:zalpha-compl}, i.e., $z^\alpha=x^{\hat{\alpha}} c^{\beta} {\overline{c}}^{\beta^*}$.

In \cite{BK95} the authors show that also terms which are not weakly resonant can be removed, but then the transformation is of finite smoothness.

\subsection{Normal form for an equilibrium point with  the local center manifold}


We will now consider a system around a saddle-center equilibrium, thus splitting the $n\times n$ matrix $A$ in complex Jordan normal form of~\eqref{eq:ODE-gen} into two matrices, a $n_A\times n_A$ hyperbolic (just called \emph{saddle} from now on) matrix $A$, whose eigenvalues have nonzero real part, and a $n_B\times n_B$ \emph{center} matrix $B$ whose eigenvalues have zero real part, with $ n =n_A+n_B$. Since we are working with real systems, all eigenvalues of $B$ appear in pairs and the Jordan normal form variables for $B$ are complex variables $c\in\mathbb{C}^r$, with $2r =n_B$, whereas the Jordan normal form variables for a matrix $A$ with $k_A$ real eigenvalues and $2l_A$ non-real eigenvalues are of the type $z=(x,\zeta)\in\mathbb{R}^{k_A}\times\mathbb{C}^{l_A}$, with $k_A+2l_A=n_A$,
so that the phase space of our system is described by $(z,c)=(x,\zeta,c)$:

\begin{equation}
 \label{eq:fpcen}
\begin{split}
    \dot{z}&=Az + N_z(z,c) ,\qquad z\in\mathbb{R}^{k_A}\times\mathbb{C}^{l_A}, \\
  \dot{c}&=Bc + N_c(z,c),\qquad c\in\mathbb{C}^r,
  \end{split}
  \end{equation}
with $N(z,c)=(N_z(z,c),N_c(z,c))=O_2(|z|^2 + |c|^2)$.

From now on we will denote our phase space as
\begin{equation}
\PhSpace= \left\{z\in\mathbb{R}^{k_A}\times\mathbb{C}^{l_A},c\in\mathbb{C}^r\right\}. \label{eq:PhSpace}
\end{equation}

We now introduce some notation for the \emph{spectral gap} of $A$.
\begin{definition}
\label{def:Alambmu}
Given a saddle matrix $A$ (i.e., $\Re{\lambda} \neq 0$ for $\lambda \in \mathrm{Sp}(A)$), there exist  four positive real numbers delimiting two minimal intervals $[-\lambda_{\mathrm{max}},-\lambda_{\mathrm{min}}]$ and $[\mu_{\mathrm{min}},\mu_{\mathrm{max}}]$ such that $\Re\mathrm{Sp}(A)\subset[-\lambda_{\mathrm{max}},-\lambda_{\mathrm{min}}]\cup [\mu_{\mathrm{min}},\mu_{\mathrm{max}}]$.
For the sake of brevity, we will just call
$\mathfrak{S}=\mathfrak{S}(A)$ to the vector formed by these numbers characterizing the  \emph{spectral gap} of $A$:
\[
\mathfrak{S}=(\lambda_{\mathrm{min}},\lambda_{\mathrm{max}},\mu_{\mathrm{min}},\mu_{\mathrm{max}})\implies
\Re\mathrm{Sp}(A)\subset[-\lambda_{\mathrm{max}},-\lambda_{\mathrm{min}}]\cup [\mu_{\mathrm{min}},\mu_{\mathrm{max}}].
\]
\end{definition}

\begin{rem}
This definition works well only for a genuine saddle matrix $A$ possessing both eigenvalues with positive real part and eigenvalues with negative real part.
In particular, if $A$ is a Hamiltonian matrix, then $\lambda_{\mathrm{min}}=\mu_{\mathrm{min}}$ and $\lambda_{\mathrm{max}}=\mu_{\mathrm{max}}$.
If $\Re\mathrm{Sp}(A)<0$ then the interval $[\mu_{\mathrm{min}},\mu_{\mathrm{max}}]$ is empty and $\mathfrak{S}(A)=(\lambda_{\mathrm{min}},\lambda_{\mathrm{max}})$.
Analogously, if $\Re\mathrm{Sp}(A)>0$, the interval $[-\lambda_{\mathrm{max}},-\lambda_{\mathrm{min}}]$ is empty and $\mathfrak{S}(A)=(\mu_{\mathrm{min}},\mu_{\mathrm{max}})$.
\end{rem}

We can now state the first main result of this paper.

\com{PZ: do we cite something here?  Mention   Bronstein and Kopanskii~\cite{BK95,BK96} and Banyaga, de la Llave and Wayne~\cite{BLW} , see also added remark after the statement of the theorem}
\begin{theorem}
\label{thm:polNForm-loccenter}
  Consider  system \eqref{eq:fpcen}, such that $(z,c)=0$ is  a non-degenerate  equilibrium point and the subspace $\{z=0\}$ is invariant, and let $\mathfrak{S}$ be the spectral gap of $A$ as introduced in Def.~\ref{def:Alambmu}.

Then for any $k\geq 1$, there exists $Q_0=Q_0(k, \mathfrak{S})$, such that for any $Q\geq Q_0$ there exists
$q_0=q_0(Q,k,\mathfrak{S})$, so that the following assertion of the theorem is valid.


If system \eqref{eq:fpcen} is $C^q$, with $q \geq q_0$ and if $q -1 \geq P \geq Q$ and $q \geq Q+3$, then there exists a $C^k$ change of variables in a neighborhood of the origin,  transforming it to
the system
\begin{equation}
 \label{eq:cen}
\begin{split}
  \dot{z}&=Az + \mathcal{N}(z,c)+ \mathcal{R}(z,c),   \\
  \dot{c}&=Bc + O_2(|z|^2 + |c|^2),
\end{split}
\end{equation}
with
\begin{eqnarray*}
\mathcal{N}(z,c)&=&\sum_{\substack{j = 1,\dots,k_A+l_A,\\
\left(j,\left(\alpha,\alpha^*,\beta,\beta^*\right)\right) \in \Res ((A,B),2,P)}} p^j_{\alpha,\alpha^*,\beta,\beta^*}\, z^{\alpha} \bar{z}^{\alpha^*} c^{\beta} \bar{c}^{\beta^*} e_j , \\ 
  \mathcal{R}(z,c)&=& \sum_{\substack{j=1,\dots,k_A+l_A,\\ |\alpha| + |\alpha^*|=1, |\beta|+|\beta^*|=P+1-Q}} g^j_{\alpha,\alpha^*,\beta,\beta^*}(z,c)
   \,z^{\alpha} \bar{z}^{\alpha^*}c^{\beta} \bar{c}^{\beta^*} e_j \notag \\
  &=&O_{P+1-Q}\left(|c|^{P+1-Q}\right)O_1(z),  
\end{eqnarray*}
where $p^j_{\alpha,\alpha^*,\beta,\beta^*}$ are constants and $g^j_{\alpha,\alpha^*,\beta,\beta^*}(z,c)$ are continuous functions.

 Moreover, if system~\eqref{eq:fpcen} is Hamiltonian, the change of variables can be chosen symplectic, and therefore the transformed system is also Hamiltonian.
\end{theorem}

Our second main result, Theorem~\ref{thm:polNForm-centerSubspaces}, is that if system~\eqref{eq:fpcen} has \emph{sign-symmetry} (see Section~\ref{sec:invWpm}), that is, it is $\SSym_s$-symmetric for any symmetry $ \SSym_s(z_+,z_-,c_1,\dots,c_m)=( s_+z_+,s_-z_-,s_1 c_1,\dots,s_m c_m)$ and any choice of signs $\{s=(s_+,s_-,s_1,\dots,s_m)\}$, then the change of variables can be chosen such that the transformed system also has sign-symmetry.

\begin{rem}
The sign symmetry is a property introduced in~\cite{CK}, which guarantees that the variational equations of system~\eqref{eq:fpcen} along ingoing and outgoing orbits have a block diagonal structure, and is an essential property to prove that heteroclinic connections between different systems~\eqref{eq:fpcen} can be shadowed, even when these heteroclinic connections are not transversal, see~\cite{DZ}. Theorem~\ref{thm:polNForm-centerSubspaces} is adapted in Lemma~\ref{lem:PolyNormForm} to the
toy model system derived in~\cite{CK}.
\end{rem}
\begin{rem}
The first results about polynomial normal forms in the saddle variables for saddle-center equilibria  are due to  Bronstein and Kopanskii~\cite{BK95}. In the pure saddle case the assertions of Theorem~\ref{thm:polNForm-loccenter} are  valid for $P=Q$ and $\mathcal{R}(z)=0$ and were proved using the deformation method by Banyaga, de la Llave and Wayne~\cite{BLW} (in this second paper the authors give a complete proof for maps, but the result also applies to ODEs), showing the preservation of geometric (Hamiltonian, contact, volume preserving) structures. Bronstein and Kopanskii~\cite{BK96} applied later on these methods to show also the preservation of geometric structures.
\end{rem}

\begin{rem}
  Bringing the remainder term to the form stated in the theorem (or  removing it in the pure saddle case) does not require that $\mathcal{N}$ contains only the resonant terms, it fact it can contain any polynomial terms of order up to $P$.
\end{rem}

\begin{rem}
In the transformed system~\eqref{eq:cen}, the term $\mathcal{N}(z,c)$ is the polynomial normal form up to degree $P$ in both variables $z$ and $c$
provided by the Poincar\'e-Dulac Theorem~\eqref{thm:Dulac-Poinc}, whereas the term $\mathcal{R}(z,c)$ contains the linear terms in the variable $z$ with coefficients of order $Q+P-1$ in the variable $c$ of the remainder of Theorem~\eqref{thm:Dulac-Poinc}. It is important to notice that, in general, the linear part of the $z$-equation in~\eqref{eq:fpcen} cannot be reduced to a constant matrix (see, for instance the review~\cite{P} for sufficient conditions), since the frequencies in the matrix $B$ are not assumed to be even incommensurable.
\end{rem}

\begin{rem}
For a Hamiltonian system~\eqref{eq:fpcen}, the Hamiltonian associated to the transformed Hamiltonian system~\eqref{eq:cen}
provided by Theorem~\ref{thm:polNForm-loccenter} takes the form $H=H_1(c)+H_2(z)+N(z,c)+R(z,c)$, where $H_1(c)$ has nondegerate quadratic part, $H_2(z)$ is quadratic (in a Hamiltonian Jordan form) and
\begin{eqnarray*}
N(z,\overline{z},c,\overline{c})&=&\sum_{\left(\alpha,\alpha^*,\beta,\beta^*\right) \in \HRes ((A,B),3,P+1)}
P_{\alpha,\alpha^*,\beta,\beta^*}\, z^{\alpha} \bar{z}^{\alpha^*} c^{\beta} \bar{c}^{\beta^*}  \\ 
R(z,\overline{z},c,\overline{c})&=& \sum_{|\alpha| + |\alpha^*|=2, |\beta|+|\beta^*|=P+1-Q} G_{\alpha,\alpha^*,\beta,\beta^*}(z,\overline{z},c,\overline{c})
   \,z^{\alpha} \bar{z}^{\alpha^*}c^{\beta} \bar{c}^{\beta^*} \notag \\
 &=&O_{P+1-Q}\left(|c|^{P+1-Q}\right)O_2(z^2),  
\end{eqnarray*}
where $P_{\alpha,\alpha^*,\beta,\beta^*}$ are constants and $G_{\alpha,\alpha^*,\beta,\beta^*}(z,c)$ are continuous functions.
\end{rem}

\begin{rem}
\label{rem:Q0q0}
The values of the constants $Q,Q_0,q_0,k$ from Theorem~\ref{thm:polNForm-loccenter}  for the general vector field as follows
\begin{eqnarray}
   Q_0(k) &=&   \max \left(k+ \left[ k \frac{\lambda_{\mathrm{max}}}{\mu_{\mathrm{min}}} +  \frac{\mu_{\mathrm{max}}}{\mu_{\mathrm{min}}}   \right],\left[(k+1)\frac{\mu_{\max}}{\mu_{\mathrm{min}}} \right]  \right) \notag \\
     & &+   \max\left(k +  \left[ k \frac{\mu_{\mathrm{max}}}{\lambda_{\mathrm{min}}} + \frac{\lambda_{\mathrm{max}}}{\lambda_{\mathrm{min}}}   \right], \left[(k+1)\frac{\lambda_{\mathrm{max}}}{\lambda_{\mathrm{min}}} \right]\right), \label{eq:Q0-our}
\end{eqnarray}
and for the Hamiltonian case
\begin{eqnarray}
  Q_0(k) &=&  2k+1 + 2\left[(k+1) \frac{\mu_{\mathrm{max}}}{\mu_{\mathrm{min}}} \right], \label{eq:Q0-ham-our}
\end{eqnarray}
In both cases
\begin{equation*}
  q_0(k,Q)=Q+3,
\end{equation*}
and in the pure saddle case $q_0(k,Q)= Q+2$.
\end{rem}
The comparison of $Q_0$ and $q_0$
with the numbers obtained in previous works \cite{BK95,BK96,BLW} in the case of a saddle is discussed in Appendix~\ref{sec:comparison}.

The proof of the above theorem has two parts. In the first part we will use the Poincar\'e-Dulac theorem to remove some non-resonant terms and in the second part we will
remove the remainder using the deformation method following the idea from \cite{BLW}. The proof of this theorem will start in the next section.

If the matrix $A$ has only real eigenvalues, we have the following corollary.
\begin{cor}
 \label{cor:polNForm-loccenterAreal}
 Under the same assumptions as in Theorem~\ref{thm:polNForm-loccenter}, assume further that all eigenvalues of $A$ are real.

 Then with the same $k,Q,P$ as in Theorem~\ref{thm:polNForm-loccenter}, system~\eqref{eq:fpcen} can be brought locally to the following form
 \begin{eqnarray*}
   \dot{z}&=&Az + \sum_{ j=1,\dots,n_A,
  (j,\alpha) \in \Res (A,2,P)} g^j_{\alpha}(c)  z^\alpha e_j \\
    & & +  \sum_{\substack{j=1,\dots,n_A, |\alpha|=1, |\beta|+|\beta^*|=P+1-Q}}  g^j_{\alpha,\beta,\beta^*}(z,c) z^\alpha c^{\beta} \bar{c}^{\beta^*} e_j, \\
     \dot{c}&=&Bc + O_2(|z|^2 + |c|^2),
\end{eqnarray*}
where $g^j_{\alpha}(c) $ are polynomials and  $g^j_{\alpha,\beta,\beta^*}(z,c)$ are continuous functions.
\end{cor}
\textbf{Proof:}
First note that since $A$ has only real eigenvalues, the variables $\bar{z}$ do not appear.

Let us take a closer look at the sum of monomials in the first equation of~\eqref{eq:cen}.  Since all eigenvalues of $B$ are purely imaginary and those of $A$ are real, then for any
$(j,(\alpha,\beta)) \in \Res ((A,B),2,P)$ we have that $(j,\alpha) \in \Res (A,2,P) $.
Therefore we have
\begin{equation*}
\sum_{\substack{j = 1,\dots,n_A,\\
(j,(\alpha,(\beta,\beta^*))) \in \Res ((A,B),2,P)}} p^j_{\alpha,\beta,\beta^*}\, z^{\alpha}  c^{\beta} \bar{c}^{\beta^*} e_j
=\sum_{(j,\alpha) \in \Res (A,2,P)} g^j_{\alpha}(c)  z^\alpha e_j
\end{equation*}
for the polynomials
$\displaystyle g^j_\alpha(c)=g^j_\alpha(c,\bar{c})
=\sum_{2-\abs{\alpha}\leq \abs{\beta}+\abs{\beta^*}\leq P-\abs{\alpha}}p^j_{\alpha,\beta,\beta^*}\, c^{\beta} \bar{c}^{\beta^*}$.
\qed

\begin{rem}
For a Hamiltonian system~\eqref{eq:fpcen}, when all the eigenvalues of $A$ are real, the Hamiltonian associated to the transformed Hamiltonian system~\eqref{eq:cen}
takes the form $H(z,c)=H_1(c)+H_2(z)+N(z,c)+R(z,c)$,
with
\begin{eqnarray*}
N(z,c,\overline{c})&=&\sum_{\alpha\in \HRes (A,3,P+1)}
P_{\alpha}(c,\bar{c}) z^{\alpha},  \\ 
R(z,c,\overline{c})&=& \sum_{|\alpha|=2, |\beta|+|\beta^*|=P+1-Q}
G_{\alpha,\beta,\beta^*}(z,c)
   \,z^{\alpha} c^{\beta} \bar{c}^{\beta^*}, 
\end{eqnarray*}
where $P_{\alpha}(c,\bar{c})$ are polynomials and $G_{\alpha,\beta,\beta^*}(z,c)$ are continuous functions.
\end{rem}

\section{The first part of the proof of Theorem~\ref{thm:polNForm-loccenter}, removal of non-resonant terms up to some order}
\label{sec:pr-1part}

Let us fix  $P \geq Q \geq 2$.
Let us write a Taylor formula for the $z$-component of our vector field
\begin{eqnarray}
 \dot{z}&=&Az + \sum_{\substack{j = 1,\dots,k_A+l_A,\\
2 \leq |\alpha| + |\alpha^*|+|\beta| + |\beta^*|\leq P}}
p^j_{\alpha,\alpha^*,\beta,\beta^*}\, z^{\alpha} \bar{z}^{\alpha^*} c^{\beta} \bar{c}^{\beta^*} e_j
+ \mathcal{R}(z,c),   \label{eq:cen-z'} \\
 \dot{c}&=&Bc + O_2(|z|^2 + |c|^2), \label{eq:cen-c'}
\end{eqnarray}
where
\begin{equation}
 \mathcal{R}(z,c)=O_{P+1}\left((|z| + |c|)^{P+1}\right).  \label{eq:cen-Rem-estm'}
\end{equation}

Using Theorem~\ref{thm:Dulac-Poinc} we can remove the resonant terms of the $z$-component in system (\ref{eq:cen-z'}--\ref{eq:cen-c'})
up to order $P$ to obtain (with a different remainder term which will denote again by $\mathcal{R}$, and which
satisfies~(\ref{eq:cen-Rem-estm'}))
\begin{eqnarray*}
 \dot{z}&=&Az + \mathcal{N}_1(z,c)+\mathcal{R}(z,c),   \\
 \dot{c}&=&Bc + O_2(|z|^2 + |c|^2), 
\end{eqnarray*}
with
\begin{equation*}
\mathcal{N}_1(z,c)=\sum_{\substack{j = 1,\dots,k_A+l_A,\\
\left(j,\left(\alpha,\alpha^*,\beta,\beta^*\right)\right) \in \Res ((A,B),2,P)}} p^j_{\alpha,\alpha^*,\beta,\beta^*}\, z^{\alpha} \bar{z}^{\alpha^*} c^{\beta} \bar{c}^{\beta^*} e_j .
\end{equation*}

The remainder can be written as sum
\begin{eqnarray*}
\mathcal{R}(c,z) = \sum_{|\alpha|=P+1} \left(\int_0^1 \frac{(1-t)^{P}}{P!} D_{\alpha}\mathcal{R}(t(c,z)) \,\mathrm{d}t\right) (c,z)^\alpha.
\end{eqnarray*}
Now we collect in $\mathcal{R}_2(z,c)$ the terms of $\mathcal{R}(z,c)$ of order greater than or equal
to $Q+1$ \emph{in the variable  $z$}, and name $\mathcal{R}_1(z,c):=\mathcal{R}(z,c)-\mathcal{R}_2(z,c)$.
In other words, split the remainder term $\mathcal{R}(z,c)$ as
\begin{eqnarray}
  \mathcal{R}(z,c)&=& \mathcal{R}_1(z,c) + \mathcal{R}_2(z,c) \nonumber\\
  \label{eq:R1Special}
   \mathcal{R}_1(z,c)&=&O_{P+1-Q}\left(c^{P+1-Q}\right)O_1(z),  \quad \mathcal{R}_2(z,c)=O_{Q+1}(z^{Q+1}).
   \label{eq:R1R2}
\end{eqnarray}

Notice that for a general remainder $\mathcal{R}(z,c)$ it would simply follow that
$\mathcal{R}_1(z,c)=O_{P+1-Q}\left(c^{P+1-Q}\right)$.
The special expression~\eqref{eq:R1Special} comes from the fact that $\mathcal{R}(z,c)$ vanishes for $z=0$ .

\begin{rem}
By assumption~\eqref{eq:cen-Rem-estm'}, $\mathcal{R}(z,c)$ is of order $P+1$ in the variables $(z,c)$, so we could even write
\begin{align*}
\mathcal{R}_1(z,c)&=O_{P+1-Q}\left(c^{P+1-Q}\right)O_1(z)O_{Q-1}\left((|z|+|c|)^{Q-1}\right),\\
\mathcal{R}_2(z,c)&=O_{Q+1}(z^{Q+1})O_{P-Q}\left((|z|+|c|)^{P-Q}\right),
\end{align*}
but as we are going to bound the factors containing $O_{Q-1}$ and $O_{P-Q}$ by constants later on, the estimates~\eqref{eq:R1R2} will be enough for us. On the other hand, since $\mathcal{R}(0,c)$ vanishes for any $c$, the same happens to the $O(|c|^{P+1})$ term in $\mathcal{R}(z,c)$
\begin{equation*}
 \left(\int_0^1 \frac{(1-t)^{P}}{P!} D_{c}^{(P+1)}\mathcal{R}(tc,tz) \,\mathrm{d}t\right) c^{P+1}.
\end{equation*}
\end{rem}

As a result of these transformations we obtain the following system
\begin{eqnarray}
  \dot{z}&=&Az + \mathcal{N}(z,c) + \mathcal{R}(z,c), \label{eq:z'-1part} \\
  \dot{c}&=&Bc +  O_2(|z|^2 + |c|^2), \label{eq:c'-1part}
\end{eqnarray}
where
\begin{eqnarray}
  \mathcal{N}(z,c)&=&\mathcal{N}_1(z,c)+\mathcal{R}_1(z,c), \nonumber \\
  \mathcal{R}(z,c)&=&\mathcal{R}_2(z,c)=O_{Q+1}(|z|^{Q+1}). \label{eq:ROq+1}
\end{eqnarray}

Moreover, due to the saddle character of the matrix $A$, it splits as
\[
A=\begin{pmatrix}A_u&0\cr0&A_s\end{pmatrix},
\]
and the matrices $A_u$, $A_s$, and $B$ are written in an appropriate Jordan form with
\begin{eqnarray}
  m_l(A_u) >0, \qquad   \mu_{log}(A_s) < 0, \qquad   \mu_{log}(B)<\delta, \qquad \mu_{log}(-B)<\delta  \label{eq:nAuAsB}
\end{eqnarray}
for any $\delta>0$ small enough fixed in advance, where $m_l$ and $\mu_{log}$ are the logarithmic norms  recalled in Appendix~\ref{app:sec-lognorm}.

The logarithmic norm $\mu_{log}(B)$ might be not equal to zero in the presence of non-trivial Jordan blocks, however in such situation it can be made arbitrarily close to $0$, by choosing a linear coordinate
system in the center direction so that the off-diagonal terms are very small.

 In the Hamiltonian case, all these linear transformations to Jordan normal form can be chosen to be symplectic, see~\cite{W,A2}. In order to have the vector field to be in the form (\ref{eq:z'-1part}--\ref{eq:ROq+1}) for Hamiltonian case we proceed as follows. For Hamiltonian $H(z,c)-H(c,0)$ we consider the Taylor formula with remainder of order $P+2$ (notice that we consider order greater than in the case of the general vector field, because the vector fields is obtained from the derivatives of Hamiltonian) and remove all
 or some non-resonant terms up to order $P+1$ and then we split the reminder gathering all terms with order $Q+2$ or higher in $z$.  As a result our
Hamiltonian takes the following form
\begin{equation}
H(z,c)=H_{c,2}(c) + N_{c}(c) +H_2(z)+N_z(z,c)+R_1(z,c) + R_2(z,c),  \label{eq:Ham-prep}
\end{equation}
 with $H_2(z)$ and $H_{c,2}(c)$ being  quadratic Hamiltonians in suitable Jordan forms and
\begin{eqnarray}
H(c,0) &=& H_{c,2}(c) + N_{c}(c) \notag\\
N_c(c) &=& O_{q+1}(|c|^3),\notag \\
N_z(z,\overline{z},c,\overline{c})&=&\sum_{\left(\alpha,\alpha^*,\beta,\beta^*\right) \in \HRes ((A,B),3,P+1)}
P_{\alpha,\alpha^*,\beta,\beta^*}\, z^{\alpha} \bar{z}^{\alpha^*} c^{\beta} \bar{c}^{\beta^*}\notag \\
R_1(z,\overline{z},c,\overline{c})&=& \sum_{|\alpha| + |\alpha^*|=2, |\beta|+|\beta^*|=P+1-Q} G_{\alpha,\alpha^*,\beta,\beta^*}(z,\overline{z},c,\overline{c})
   \,z^{\alpha} \bar{z}^{\alpha^*}c^{\beta} \bar{c}^{\beta^*} \notag \\
 &=&O_{P+1-Q}\left(|c|^{P+1-Q}\right)O_2(z^2),\notag \\
R_2(z,c)&=&O_{Q+2}(|z|^{Q+2}), \label{eq:Ham-prepR2}
\end{eqnarray}
where $P_{\alpha,\alpha^*,\beta,\beta^*}$ are constants and $G_{\alpha,\alpha^*,\beta,\beta^*}(z,c)$ are continuous functions.

Now our goal is to apply the deformation method to remove the remainder term $\mathcal{R}_1$ (or $R_1$ in Hamiltonian setting). This starts in the next section.

\section{ Derivation of the cohomological equation}

The process of removing the remainder  term $\mathcal{R}$ in (\ref{eq:z'-1part}--\ref{eq:c'-1part})
through the deformation method, involves several stages and is concluded at the end of section~\ref{sec:estm-gen-coheq}.
The first step is to derive a \emph{cohomological equation}.

Consider a vector field dependent on a parameter $\varepsilon$
\begin{equation}
 \mathcal{Z}_\varepsilon= \mathcal{Z}_0 + \varepsilon \mathcal{R}(z).  \label{eq:defZ-eps}
\end{equation}

Following the idea of the deformation method taken from \cite{BLW}, to remove the term $\varepsilon \mathcal{R}(z)$ from equation~\eqref{eq:defZ-eps} we would like to find a family of diffeormorphisms  $g_\varepsilon$ such that 
\begin{equation}
  g_{\varepsilon}^\ast \mathcal{Z}_\varepsilon = \mathcal{Z}_0. \label{eq:gupastZe}
\end{equation}

Since $g_\varepsilon$ is a smooth family of diffeomorphisms, we can find  a family of vector
fields $\mathcal{G}_\varepsilon$  such that
\begin{eqnarray}
  \frac{d}{d \varepsilon} g_\varepsilon = \mathcal{G}_\varepsilon  \circ g_\varepsilon.  \label{eq:geps-ode}
\end{eqnarray}
This more carefully written means that
\begin{eqnarray*}
  \frac{\partial g}{\partial \varepsilon} (\varepsilon,x) = \mathcal{G}(\varepsilon,g(\varepsilon,x)),  
\end{eqnarray*}
so that
\begin{equation*}
 \mathcal{G}(\varepsilon,x)=\frac{\partial g}{\partial \varepsilon} (\varepsilon,g_\varepsilon^{-1}(x)).
\end{equation*}
Note that (\ref{eq:geps-ode}) is a  non-autonomous o.d.e.


In \cite{BLW} the main focus was on diffeomorphims and  for vector fields the authors  claim (without a proof) that
for the pushforward $\left(g_\varepsilon\right)_\ast$ it holds that
\begin{equation}
    \frac{d}{d \varepsilon} \left(g_\varepsilon\right)_\ast \mathcal{Z}_\varepsilon= \left(g_\varepsilon\right)_{\ast}
    \left([\mathcal{Z}_\varepsilon,\mathcal{G}_\varepsilon] + \frac{d}{d\varepsilon}\mathcal{Z}_\varepsilon\right), \label{eq:der-g-eps}
  \end{equation}
which is not  true, since the right formula is
\begin{equation*}
    \frac{d}{d\varepsilon} (g_{\varepsilon \ast}) \mathcal{Z}_\varepsilon = (g_{\varepsilon \ast})\left([\mathcal{Z}_\varepsilon,\mathcal{G}_\varepsilon]  +  \frac{d \mathcal{Z}_\varepsilon}{d \varepsilon}\right)
 + [(g_{\varepsilon \ast})\mathcal{Z}_\varepsilon,\mathcal{G}_\varepsilon - (g_{\varepsilon \ast})\mathcal{G}_\varepsilon  ].  
\end{equation*}

However it turns out that equation (\ref{eq:der-g-eps}) works for the pullback $g^\ast$. Namely, we have
\begin{lemma}
\label{lem:corr-vect-field}
  Assume that  $g \in C^2$ and $\mathcal{G} \in C^1$    are such that
    $g_0=\Id$
   and
   \begin{eqnarray}
  \frac{d}{d \varepsilon} g_\varepsilon = \mathcal{G}_\varepsilon  \circ g_\varepsilon.  \label{eq:heps-ode}
\end{eqnarray}
Assume that $\mathcal{Z}(\varepsilon,z)$ is $C^1$.   Then
  \begin{equation}
      \frac{d (g^\ast_\varepsilon \mathcal{Z}_\varepsilon)}{d \varepsilon}= g_\varepsilon^\ast \left([\mathcal{Z}_\varepsilon,\mathcal{G}_\varepsilon] + \frac{d \mathcal{Z}_\varepsilon}{d \varepsilon}\right). \label{eq:g*corr-vect-field}
  \end{equation}
\end{lemma}
\noindent
\textbf{Proof:}
From (\ref{eq:heps-ode}) we obtain
\begin{equation}
\frac{\partial }{\partial \varepsilon} \frac{\partial }{\partial x} g(\varepsilon,x) = \frac{\partial \mathcal{G}}{\partial x}(\varepsilon,g(\varepsilon,x)) \frac{\partial g}{\partial x}(\varepsilon,x)=D \mathcal{G}_\varepsilon(g(\varepsilon,x))Dg_\varepsilon(x).
\label{eq:hode-var}
\end{equation}

Let us denote
\begin{equation*}
 Y(\varepsilon,u)=\left(g^\ast_\varepsilon \mathcal{Z}_\varepsilon\right)(u). 
\end{equation*}
Then by (\ref{eq:hast})
\begin{eqnarray*}
  Dg_\varepsilon(u) Y(\varepsilon,u)=\mathcal{Z}_\varepsilon(g_\varepsilon(u)).
\end{eqnarray*}
Observe that $Y(\varepsilon,u)$ is $C^1$.

We differentiate the above equality with respect to $\varepsilon$ . For the left hand side we obtain (we use (\ref{eq:hode-var}))
\begin{eqnarray*}
  \frac{d}{d \varepsilon} \left( Dg_\varepsilon(u) Y(\varepsilon,u) \right)= \left(\frac{\partial }{\partial \varepsilon} Dg_\varepsilon(u)\right) Y(\varepsilon,u) + Dg_\varepsilon(u) \frac{\partial Y}{\partial \varepsilon}(\varepsilon,u) \\
  =  \ D\mathcal{G}_\varepsilon(g_\varepsilon(u)) Dg_\varepsilon(u) Y(\varepsilon,u) + Dg_\varepsilon(u) \frac{\partial Y}{\partial \varepsilon}(\varepsilon,u) \\
    =D\mathcal{G}_\varepsilon(g_\varepsilon(u)) \mathcal{Z}_\varepsilon(g_\varepsilon (u)) + Dg_\varepsilon(u) \frac{\partial Y}{\partial \varepsilon}(\varepsilon,u),\\
\end{eqnarray*}
while for the right hand side we have
\begin{eqnarray*}
    \frac{d}{d \varepsilon}\left( \mathcal{Z}_\varepsilon(g_\varepsilon(u)) \right)&=&\frac{\partial \mathcal{Z}}{\partial \varepsilon}(\varepsilon,g_\varepsilon(u))
       + D\mathcal{Z}_\varepsilon(g_\varepsilon(u)) \frac{\partial g}{\partial \varepsilon}(\varepsilon,u)\\
       &=&\frac{\partial \mathcal{Z}}{\partial \varepsilon}(\varepsilon,g_\varepsilon(u))
       + D\mathcal{Z}_\varepsilon(g_\varepsilon(u)) \mathcal{G}_\varepsilon (g_\varepsilon(u)).
\end{eqnarray*}
Therefore we obtain (we use also (\ref{eq:commutator}))
\begin{eqnarray*}
  Dg_\varepsilon(u) \frac{\partial Y}{\partial \varepsilon}(\varepsilon,u) = -D\mathcal{G}_\varepsilon(g_\varepsilon(u)) \mathcal{Z}_\varepsilon(g_\varepsilon(u))
       + D\mathcal{Z}_\varepsilon(g_\varepsilon(u)) \mathcal{G}_\varepsilon (g_\varepsilon(u) \\
       +  \frac{d \mathcal{Z}_\varepsilon}{d \varepsilon}(g_\varepsilon(u))  \\
       = \left([\mathcal{G}_\varepsilon,\mathcal{Z}_\varepsilon] +  \frac{d \mathcal{Z}_\varepsilon}{d \varepsilon}\right)(g_\varepsilon(u)) .
\end{eqnarray*}
This establishes (\ref{eq:g*corr-vect-field}).
\qed

Therefore we obtain that the desired conjugation~\eqref{eq:gupastZe} is equivalent to the \emph{cohomological equation}
\begin{equation*}
[\mathcal{G}_\varepsilon,\mathcal{Z}_\varepsilon]= -\mathcal{R}
\quad (\text{or }
 [\mathcal{Z}_\varepsilon,\mathcal{G}_\varepsilon]= \mathcal{R}),
\end{equation*}
which is now a \emph{linear} equation for $\mathcal{G}_\varepsilon$.


\subsection{Cohomological equation for Hamiltonian  systems}

Assume that we have a symplectic form $\Omega$ and let $J(z)$ be the associated matrix defining the vector Hamiltonian vector field (see (\ref{eq:def-ham-vf},\ref{eq:def-ham-vf-J})).
We assume that we have an $\varepsilon$-dependent  family of Hamiltonians $Z_\varepsilon(z)=Z(\varepsilon,z)$ inducing a $\varepsilon$-dependent family
of Hamiltonian vector fields $\mathcal{Z}_\varepsilon$.

We would like to find a family of symplectic transformations (with respect to the form $\Omega$) $z=h_\varepsilon(u)$  such that
\begin{equation*}
   h_\varepsilon^* Z_\varepsilon = Z_\varepsilon \circ h_\varepsilon =Z_0.
\end{equation*}

We will seek $h_\varepsilon$ as the time shift $\varepsilon$ along the trajectory of some $\varepsilon$-dependent Hamiltonian $H_\varepsilon$.

\begin{lemma}
Assume that $Z(\varepsilon,z)$ is $C^1$ and $\mathcal{Z}_\varepsilon(z)=J(z) (DZ_\varepsilon(z))^\top$.
  Assume that  $h(\varepsilon,z) \in C^2$ and $H(\varepsilon,z) \in C^2$   are such that
    $h_0=\Id$
   and
   \begin{eqnarray*}
  \frac{d}{d \varepsilon} h_\varepsilon = J \cdot (D H_\varepsilon)^\top  \circ h_\varepsilon.  
\end{eqnarray*}
Then
\begin{equation*}
\frac{d}{d \varepsilon} \left( h_{\varepsilon}^\ast Z_\varepsilon \right) = h_{\varepsilon}^\ast \left(\frac{\partial Z_\varepsilon}{\partial \varepsilon} - D H_\varepsilon \cdot \mathcal{Z}_\varepsilon \right).
\end{equation*}

\end{lemma}
\noindent
\textbf{Proof:}
We have
\begin{eqnarray*}
  \frac{d}{d \varepsilon} Z_\varepsilon(h_\varepsilon(u)) &=& \frac{\partial Z}{\partial \varepsilon}(\varepsilon,h_\varepsilon(u))
   + DZ_\varepsilon(h_\varepsilon(u)) \cdot \frac{d}{d\varepsilon} h_\varepsilon(u) \\
   &=& \frac{\partial Z_\varepsilon}{\partial \varepsilon}(h_\varepsilon(u)) + DZ_\varepsilon(h_\varepsilon(u)) \cdot \left(J(h_\varepsilon(u))  (DH_\varepsilon(h_\varepsilon(u)))^\top \right)\\
   &=& \frac{\partial Z_\varepsilon}{\partial \varepsilon}(h_\varepsilon(u)) - \left(J(h_\varepsilon(u)) ( DZ_\varepsilon(h_\varepsilon(u)))^\top\right)^\top \cdot  (DH_\varepsilon(h_\varepsilon(u)))^\top  \\
   &=& \frac{\partial Z_\varepsilon}{\partial \varepsilon}(h_\varepsilon(u)) - \left(\mathcal{Z}_\varepsilon(h_\varepsilon(u))\right)^\top \cdot  (DH_\varepsilon(h_\varepsilon(u)))^\top  \\
   &=& \frac{\partial Z_\varepsilon}{\partial \varepsilon}(h_\varepsilon(u))  - D H_\varepsilon( h_\varepsilon(u)) \cdot \mathcal{Z}_\varepsilon( h_\varepsilon(u)).
\end{eqnarray*}

\qed

Therefore the cohomological equation to be solved  in the Hamiltonian context is
\begin{equation}
  DH_\varepsilon \cdot \mathcal{Z}_\varepsilon -  \frac{\partial  Z_\varepsilon}{\partial \varepsilon}=0.  \label{eq:coh-ham-derZ}
\end{equation}
Observe that equation (\ref{eq:coh-ham-derZ}) does not depend on the particular formula for $\Omega$. The only way $\Omega$ factors in (\ref{eq:coh-ham-derZ})
is the relation between $Z$ and $\mathcal{Z}$.

\begin{rem}
In fact in the proof of the above lemma we obtained
 \begin{equation*}
   \frac{d}{d \varepsilon} Z_\varepsilon(h_\varepsilon(u)) = \left(\frac{d Z_\varepsilon}{d \varepsilon} + \left\{Z_\varepsilon,H_\varepsilon\right\}\right)(h_\varepsilon(u)), 
 \end{equation*}
 where the \emph{Poison bracket}  is given by
\begin{equation*}
   \left\{Z_\varepsilon,H_\varepsilon\right\}=\Omega(I (dZ_\varepsilon), I(dH_\varepsilon)),
\end{equation*}
where $\mathcal{Z}_\varepsilon=I (dZ_\varepsilon)=J\cdot (dZ_\varepsilon)^\top$ and $\mathcal{H}_\varepsilon=I (dH_\varepsilon)=J\cdot (dH_\varepsilon)^\top$ are defined by (\ref{eq:def-ham-vf}).
\end{rem}

\section{Strategy for solving the cohomological equation}

In this section we represent the points of our phase space by $z=(x,y,c)$, where $c$ are the center variables and $(x,y)$ the saddle variables, $x$ unstable and $ y$ stable.

We set (compare to (\ref{eq:z'-1part}--\ref{eq:c'-1part}))
\begin{equation}
 \mathcal{Z}_\varepsilon(z)= \left(A(x,y) + \mathcal{N}_{x,y}(z) + \varepsilon \mathcal{R}(z), B(c) + \mathcal{N}_c(z) \right). \label{eq:Z-with-rem}
\end{equation}

 From  equation (\ref{eq:g*corr-vect-field}) in Lemma~\ref{lem:corr-vect-field} (and (\ref{eq:commutator})) it follows that to bring (\ref{eq:Z-with-rem}) to the normal form (i.e. to remove $\mathcal{R}$) we need first to solve for $\mathcal{G}_\varepsilon$ the following \emph{cohomological equation}
\begin{equation}
  D \mathcal{G}_\varepsilon \cdot \mathcal{Z}_\varepsilon = D \mathcal{Z}_\varepsilon  \cdot \mathcal{G}_\varepsilon  - \mathcal{R}, \label{eq:coh-flows}
\end{equation}
where with some abuse of notation we redefined $\mathcal{R}$ by setting $\mathcal{R}(z)=(\mathcal{R}(z),0)$.

In the symplectic setting we assume that the Hamiltonian function $Z_\varepsilon$ is of the following form
 \begin{equation}
  Z_\varepsilon(z)=A(z) + N(z)+\varepsilon R(z), \label{eq:ham-Z-with-rem}
\end{equation}
where (compare (\ref{eq:Ham-prep})) $A(z)=H_{2,c}(c) + H_2(x,y)$, $N(z)= N_c(c)+ N_z(x,y,c) + R_1(x,y,c)$ and $R(z)=R_2(x,y,c)$.

The cohomological equation for a Hamiltonian $G_\varepsilon$ will be (see~(\ref{eq:coh-ham-derZ}))
\begin{equation}
  DG_\varepsilon  \cdot \mathcal{Z}_\varepsilon =  R. \label{eq:coh-ham}
\end{equation}

In both equations (\ref{eq:coh-flows}) and (\ref{eq:coh-ham}) the variable $\varepsilon$ plays a role of a parameter. In fact we have  a family of equations parameterized by $\varepsilon$, but we want to construct the solution which will depend
smoothly on $\varepsilon$.  Observe that the solution of (\ref{eq:coh-flows}) (or (\ref{eq:coh-ham})) cannot be unique, since we can always add $f(\varepsilon) \mathcal{Z}_\varepsilon$ (or $f(\varepsilon) Z_\varepsilon$) for any function $f(\varepsilon)$ to get a new solution.

In the case of general vector field the strategy to construct the solution of the cohomological equation is as follows (for Hamiltonian system the procedure is analogous)

\begin{itemize}
\item Preparation of ``compact data";
after this step we will have a modified vector field $\tilde{\mathcal{Z}}_\varepsilon$ and $\tilde{\mathcal{R}}$, coinciding with the original one in some neighborhood of the origin, but with a linear vector field $\tilde{\mathcal{Z}}_\varepsilon$ far from the origin and a remainder $\tilde{\mathcal{R}}$ with compact support, see section~\ref{subsec:prep-compt-data}.
After this step, we will drop the tildes in the modified equations.

\item Straightening of the local center-stable and center-unstable manifolds via a transformation $T_\varepsilon$ (see Section~\ref{subsec:straight-man}). We obtain a new vector
field $\mathcal{Z}_\varepsilon=T_{\varepsilon \ast} \mathcal{Z}_\varepsilon$ and a new remainder $\mathcal{R}_\varepsilon=T_{\varepsilon \ast} \mathcal{R}$.

\item Splitting of the remainder term, $T_\ast \mathcal{R}=\mathcal{R}_1+\mathcal{R}_2$,
where $\mathcal{R}_1=O_{\ell_1+1}(|x|^{\ell_1+1})$ \newline and  $\mathcal{R}_2=O_{\ell_2+1}(|y|^{\ell_2+1})$ (see Sec.~\ref{subsec:prep-R1-R2}).

\item Resolution of (\ref{eq:coh-flows}) (or  (\ref{eq:coh-ham})) with remainders $\mathcal{R}_1$ and $\mathcal{R}_2$ (see Section~\ref{subsec:solve-fsp}), obtaining $\mathcal{G}_1(z)$ and $\mathcal{G}_2(z)$, respectively.
\item Finally,  the solution is given by $\mathcal{G}=\left(T^{-1}\right)_{\ast}(\mathcal{G}_1 + \mathcal{G}_2)$.
\end{itemize}
The whole process of construction of solution is smooth in $\varepsilon$ and other parameters if they are present in $\mathcal{R}$ and $\mathcal{Z}_\varepsilon$. This is important, because the solution
of the cohomological equation defines  a  vector field (depending on $\varepsilon$), for which the time shift from $\varepsilon=0$ to $\epsilon=1$ along a trajectory defines the desired transformation (with $\varepsilon$ being the time variable).


\subsection{Preparation of compact data}
\label{subsec:prep-compt-data}
Consider system (\ref{eq:z'-1part}--\ref{eq:c'-1part}) obtained at the end of the first part of the proof of Theorem~\ref{thm:polNForm-loccenter}
in Section \ref{sec:pr-1part}.
 We intend  to modify the vector field, away from the origin so that all nonlinearities will have a compact support.

Recall (see (\ref{eq:PhSpace})) that our phase-space is $\PhSpace= \left\{(z\in\mathbb{R}^{k_A}\times\mathbb{C}^{l_A},c\in\mathbb{C}^m)\right\}$.  We have the splitting $z=(x,y)$, where $x \in \mathbb{R}^{n_u}$
and $y \in \mathbb{R}^{n_s}$ with $k_A + 2l_A=n_u + n_s$.

Let $\eta: \PhSpace  \to [0,1]$ be a $C^\infty$ function such that for some $0<r_0 < r_1$ it holds that
\begin{eqnarray*}
  \eta(z)=1, \quad \|z\| \leq r_0, \\
  \eta(z)=0, \quad \|z\| \geq r_1.
\end{eqnarray*}

Let the constants $K_\eta$, $K_N$ be such that
\begin{eqnarray*}
 \|D\eta(z)\| &\leq& K_\eta,  \quad \|D^2\eta(z)\| \leq K_\eta  \quad \forall z \in \PhSpace \\
  \|N(z)\| &\leq& K_N \|z\|^2, \quad \mbox{if} \quad  \|z\| \leq r_1, \\
  \|DN(z)\| &\leq& K_{N} \|z\|, \quad \mbox{if} \quad  \|z\| \leq r_1,  \\
   \|D^2N(z)\| &\leq& K_{N}, \quad \mbox{if} \quad  \|z\| \leq r_1.
\end{eqnarray*}

  Let $\sigma>0$ be some small number.
Let us set
\begin{equation*}
   \widetilde{N}(x,y,c)= N(x,y,c)\eta((x,y,c)/\sigma),\quad  \widetilde{\mathcal{R}}(x,y,c)=\mathcal{R}(x,y,c)\eta((x,y,c)/\sigma).
\end{equation*}

The modified system is
\begin{eqnarray}
 \dot{x}&=&A_u x + \widetilde{N}_x(x,y,c) + \varepsilon \widetilde{\mathcal{R}}_x(x,y,c)=A_u x + \widetilde{M}_x(\varepsilon,x,y,c),  \label{eq:dotx-lcd-app2} \\
  \dot{y}&=&A_s y + \widetilde{N}_y(x,y,c)+ \varepsilon \widetilde{\mathcal{R}}_y(x,y,c)=A_s y + \widetilde{M}_y(\varepsilon,x,y,c),\nonumber \\
  \dot{c}&=&Bc + \widetilde{N}_c(x,y,c). \label{eq:dotc-lcd-app2}
\end{eqnarray}
Let $\widetilde{M}(\varepsilon,x,y,c)=(\widetilde{M}_x(\varepsilon,x,y,c),\widetilde{M}_y(\varepsilon,x,y,c),0)$.
Obviously, we have (see \eqref{eq:ROq+1}) (where all the $O(\cdot)$ terms are for their arguments converging to $0$)
\begin{eqnarray}
  \widetilde{M}(\varepsilon, 0,0,c)&=&0, \quad \widetilde{N}_{x,y}(0,0,c)=0, \label{eq:tNzc00=0} \\
  \widetilde{M}(\varepsilon, x,y,c)&=&O_2(|(x,y,c)|^2), \quad \widetilde{N}(x,y,c)=O_2(|(x,y,c)|^2),\nonumber \\
  \widetilde{\mathcal{R}}(x,y,c)&=&O_{Q+1}(|(x,y)|^{Q+1}).\nonumber 
\end{eqnarray}

Note that we have defined $\widetilde{M}$ as a notation for the nonlinearities appearing in the saddle directions. Indeed, $\widetilde{M}$ also depends on $\varepsilon$ , but since we have bounds that are uniform for $\varepsilon \in [0,1]$ , we will usually omit this dependence on $\varepsilon$ in the sequel.

We have for $\sigma \leq 1$ and any $z$ (for suitable constant $K'_N$)
\begin{eqnarray}
 \|\widetilde{N}(z)\| &\leq& \sup_{\|z\| \leq \sigma r_1} \|N(z)\| \leq K_N \sigma^2 r_1^2  \leq K'_{N} \sigma^2, \label{eq:tildeN-estm} \\
 \|D \widetilde{N}(z)\| &\leq& \sup_{\|z\| \leq \sigma r_1} \left(\frac{1}{\sigma} \|D\eta(z/\sigma)\| \cdot \|N(z)\|  + |\eta(z/\sigma)| \cdot \|DN(z)\|  \right) \nonumber   \\
   &\leq& \sigma^{-1} K_\eta K_N \sigma^2 r_1^2 + K_{N} \sigma r_1= \sigma (K_\eta K_N r_1^2 + K_{N} r_1) \leq K'_N \sigma, \label{eq:DtildeN-estm} \\
  \|D^2 \widetilde{N}(z)\| &\leq&  \sup_{\|z\| \leq \sigma r_1} \left(\frac{1}{\sigma^2} \|D^2\eta(z/\sigma)\| \cdot \|N(z)\|  + \frac{2}{\sigma} \|D\eta(z/\sigma)\| \cdot \|DN(z)\| \right.  \nonumber \\
  &&\left. + \eta(z/\sigma) \|D^2N(z)\| \right)
  \leq K_\eta K_N r_1^2 + 2 K_\eta K_N r_1 + K_N=K'_N. \nonumber 
\end{eqnarray}

Analogously, for some constant $K_\mathcal{R}$, we obtain that for any $w=(x,y,c)$
\begin{eqnarray}
  \|\widetilde{\mathcal{R}}(w)\| &\leq& K_\mathcal{R} \sigma^{Q+1}, \nonumber \\
 \|D\widetilde{\mathcal{R}}(w)\| &\leq& K_{\mathcal{R}} \sigma^{Q}, \label{eq:DtilR-estm}\\
 \|D^2\widetilde{\mathcal{R}}(w)\| &\leq& K_{\mathcal{R}} \sigma^{Q-1}. \nonumber
\end{eqnarray}


For the Hamiltonian system we do the following modification of the Hamilton function $Z$  (compare (\ref{eq:ham-Z-with-rem})):
\begin{equation}
  \widetilde{Z}(z)= A(z) + N(z)\eta(z/\sigma) + \varepsilon R(z)\eta(z/\sigma)= A(z) + \tilde{N}(z) + \varepsilon \tilde{R}(z).  \label{eq:mod-ham}
\end{equation}
It is easy to see that the induced differential equation can be written as (\ref{eq:dotx-lcd-app2}--\ref{eq:dotc-lcd-app2}) satisfying conditions (\ref{eq:tNzc00=0}--\ref{eq:DtilR-estm}) for the nonlinear terms.

Let us denote
$$W^c=\{z=0\}.$$
With some abuse of notation we will also  treat  $W^c$ as the $\mathbb{C}^r$ component of $\PhSpace$, for example in a lemma below we write $[0,1] \times \mathbb{R}^{n_u} \times W^c $
as a domain of some function or in $D(\delta)=\overline{B}_{n_u}(\delta) \times \overline{B}_{n_s}(\delta)\times W^c$ to represent some neighborhood of $W^c$.

\begin{lemma}
\label{lem:iso-block-mod-vekf}
Let us consider the system (\ref{eq:dotx-lcd-app2}--\ref{eq:dotc-lcd-app2}) and denote by $\varphi(t,x,y,c)$ the induced flow. Let $L < 1$.

Then there exists $\sigma_0 >0$ so that, for any $0 < \sigma \leq \sigma_0$, it holds that

 \begin{itemize}
  \item  There exist functions $y^\textrm{u}: [0,1] \times \mathbb{R}^{n_u} \times W^c  \to \mathbb{R}^{n_s}$
and $x^\textrm{s}: [0,1] \times \mathbb{R}^{n_s} \times W^c  \to  \mathbb{R}^{n_u}$ in $C^{q-1}$,
 such that for any $\varepsilon \in [0,1]$
\begin{eqnarray*}
  W^{cu}=\{(x,y^\textrm{u}(\varepsilon,x,c),c), \quad (x,c) \in \mathbb{R}^{n_u} \times W^c\}, \\
 W^{cs}=\{(x^\textrm{s}(\varepsilon,y,c),y,c), \quad (y,c) \in \mathbb{R}^{n_s} \times W^c \},
\end{eqnarray*}
and for any $\varepsilon_1$,$\varepsilon_2$, $c_1$,$c_2$, $x_1,x_2$ and $y_1,y_2$, it holds that
\begin{eqnarray}
  \|y^\textrm{u}(\varepsilon_1,x_1,c_1) - y^\textrm{u}(\varepsilon_2,x_2,c_2)\| &\leq&
      L \left(|\varepsilon_1 - \varepsilon_2|+ \|c_1-c_2\|+\|x_1 - x_2\|\right),  \label{eq:Lip-yu}\\
    \|x^\textrm{s}(\varepsilon_1,y_1,c_1) - x^\textrm{s}(\varepsilon_2,y_2,c_2) \| &\leq&
       L  \left(|\varepsilon_1 - \varepsilon_2|+ \|c_1-c_2\|+\|y_1 - y_2\|\right). \label{eq:Lip-xs}
\end{eqnarray}
\item There exists a constant $K_W=K_W(L,\sigma)$, such that for all $j=1,\dots,q-1$,
\begin{equation}
   \|D^j x^{\textrm{s}}\| \leq K_W, \quad  \|D^j y^{\textrm{u}}\| \leq K_W. \label{eq:Dwcscu-bnd}
\end{equation}
This function $K(L,\sigma)$ is non-decreasing with respect to $L$ and $\sigma$.
\item For any $\delta >0$, introduce the set
$D(\delta)=\overline{B}_{n_u}(\delta) \times \overline{B}_{n_s}(\delta)\times W^c$. Then if $(x,y,c) \notin W^{cs}$, then there exists $t_0\geq 0$ such that $\varphi(t,x,y,c) \notin D(\delta)$ for $t \geq t_0$ and $\varepsilon \in [0,1]$;
and if $(x,y,c) \notin W^{cu}$, then there exists $t_0 \leq 0$ such that $\varphi(t,x,y,c) \notin D(\delta)$ for $t \leq t_0$ and $\varepsilon \in [0,1]$.
\end{itemize}
\end{lemma}
\textbf{Proof:}
We will begin by proving that for any $\delta > 0$ the set $D(\delta)=\overline{B}_{n_u}(\delta) \times \overline{B}_{n_s}(\delta)\times W^c$ is an \emph{isolating block} (see Def.~\ref{def:isolating-segment} in  Appendix~\ref{subsec:NHIM-ode}) for this system  for any $\varepsilon \in [0,1]$.
To prove this, we have to check  exit and entry conditions, i.e. (\ref{eq:isoblock-exit}) and (\ref{eq:isoblock-entry}), from the Appendix.

We have for $(x,y,c) \in D(\delta)$ (we use (\ref{eq:eucl-log-norm},\ref{eq:eucl-m-l-norm}) and (\ref{eq:DtildeN-estm}, \ref{eq:DtilR-estm})) that
\begin{eqnarray*}
(A_u x,x) &\geq& m_l(A_u) \|x\|^2 ,\\
 (A_s y,y) &\leq& \mu_{log}(A_s) \|y\|^2, \\
  \|\widetilde{N}(x,y,c)\| &=& \|\widetilde{N}(x,y,c) - \widetilde{N}(c,0,0)\| \leq \|D\widetilde{N}\| \sqrt{2} \delta \leq \sqrt{2}K'_{N} \sigma \delta , \\
  \|\widetilde{\mathcal{R}}(x,y,c)\| &=& \|\widetilde{\mathcal{R}}(x,y,c) - \widetilde{\mathcal{R}}(c,0,0)\| \leq \|D\widetilde{\mathcal{R}}\| \sqrt{2} \delta \leq \sqrt{2}K_{\mathcal{R}} \sigma^Q \delta.
\end{eqnarray*}
Hence, if  $(x,y,c)\in D(\delta)^- $ (i.e. $\|x\|=\delta$)
\begin{eqnarray*}
  (\dot{x},x)&=&(A_u x,x) + (\widetilde{N}(x,y,c),x) +  \varepsilon (\widetilde{\mathcal{R}}(x,y,c),x),  \\
  &\geq& \delta^2 \left( m_l(A_u) -\sqrt{2}K'_{N} \sigma  -  \sqrt{2}K_{\mathcal{R}} \sigma^Q \right),
\end{eqnarray*}
and for  $(x,y,c) \in D(\delta)^+$ (i.e. $\|y\|=\delta$)
\begin{eqnarray*}
  (\dot{y},y) &=&(A_s y,y) + (\widetilde{N}(x,y,c),y) +  \varepsilon (\widetilde{\mathcal{R}}(x,y,c),y) \\
    &\leq& -\delta^2 \left( -\mu_{log}(A_s) -\sqrt{2}K'_{N} \sigma  - \sqrt{2}K_{\mathcal{R}} \sigma^Q  \right).
\end{eqnarray*}
Hence if $\sigma $ is such that
\begin{eqnarray*}
 m_l(A_u) > \sqrt{2}K'_{N} \sigma  +  \sqrt{2}K_{\mathcal{R}} \sigma^Q,  \\
  -\mu_{log}(A_s) > \sqrt{2}K'_{N} \sigma  + \sqrt{2}K_{\mathcal{R}} \sigma^Q,
\end{eqnarray*}
then conditions (\ref{eq:isoblock-exit}) and (\ref{eq:isoblock-entry}) are satisfied and therefore $D(\delta)$ is an isolating block. In view of (\ref{eq:nAuAsB}), this holds if $\sigma$ is small enough.

To deal with $W^{cu}$ and $W^{cs}$ and their dependence on $\varepsilon$
we add to system (\ref{eq:dotx-lcd-app2}--\ref{eq:dotc-lcd-app2}) the  equation
\begin{equation}
 \dot{\varepsilon}= 0. \label{eq:doteps-lcd-app}
\end{equation}
This allows us to obtain the dependence with respect to the parameter $\varepsilon$ from the general arguments of \cite{CZ15,CZ17} on the existence of a normally hyperbolic invariant manifold (NHIM) and its center-unstable and center-stable manifolds, which are recalled in the Appendix~\ref{sec:NHIM}.
Note that the loss of one unit in the degree of differentiability of $ W^{cu}$ and $ W^{cs}$ comes in general situation from the loss of differentiability of the center manifold $ W^{c}$.

Now we check the rate conditions (see Def.~\ref{def:rate-cond-ode} and Theorem~\ref{thm:nhim-ode} in  Appendix~\ref{subsec:NHIM-ode}) for system (\ref{eq:dotx-lcd-app2}--\ref{eq:dotc-lcd-app2},\ref{eq:doteps-lcd-app}) up to order $q-1$.
To be in agreement with the notation used in Theorem~\ref{thm:nhim-ode} we set $\Lambda=[0,\varepsilon] \times W^c$. Hence the center direction denoted there by $\lambda$
is now $(\varepsilon,c)$.

On $\tilde{D}(\delta)=[0,\varepsilon]  \times D(\delta)$ we have
\begin{align*}
\overrightarrow{\mu _{s,1}}&=& \mu_{log}(A_s) + O(\sigma)+\frac{1}{L}O(\sigma)  ,&&
\overrightarrow{\mu _{s,2}}&=& \mu_{log}(A_s) + O(\sigma) +L O(\sigma), \\
\overrightarrow{\xi _{u,1}}&=& m_{l}(A_u) - O(\sigma) -\frac{1}{L} O(\sigma)  , &&
\overrightarrow{\xi _{u,1,P}}&=&m_{l}(A_u) - O(\sigma) -\frac{1}{L}O(\sigma)  \\
\overrightarrow{\mu _{cs,1}}&=& \mu_{log}(B) + O(\sigma)+L  O(\sigma)  ,  &&
\overrightarrow{\mu _{cs,2}}&=& \mu_{log}(B) + O(\sigma) +\frac{1}{L} O(\sigma) \\
\overrightarrow{\xi _{cu,1}}&=& -O(\sigma) -L O(\sigma)  , &&
\overrightarrow{\xi _{cu,1,P}}&=&-O(\sigma) -L O(\sigma),  \\
\overrightarrow{\xi _{u,2}} &=& m_{l}(A_u) -O(\sigma) - L O(\sigma), &&
\overrightarrow{\xi _{cu,2}}&=&-O(\sigma) -\frac{1}{L} O(\sigma)
\end{align*}
It is easy to see that for any $k$ and $L$, there exists $\sigma_0=\sigma_0(k,L)$, such that for $\sigma \leq \sigma_0$ the rate conditions of order $k$ are satisfied (we need
$k=q-1$).

Therefore $y^u$ and $x^s$, which in the notation of Theorem~\ref{thm:nhim-ode} are $w^{cu}$ and $w^{cs}$, respectively, satisfy the Lipschitz condition with constant $L$. This establishes inequalities (\ref{eq:Lip-yu},\ref{eq:Lip-xs}).

Since all partial derivatives up to the order of the regularity class of our modified vector field are globally bounded, we obtain condition (\ref{eq:Dwcscu-bnd}) of Theorem~\ref{thm:nhim-ode}.

The above arguments apply to local center-stable and center-unstable manifolds  inside $D(\delta)$. But since $\delta$ is arbitrary and the bounds on derivatives of our vector field are global,
then all the above estimates also apply to the global manifolds $W^{cs}$ and $W^{cu}$.

\qed

\begin{lemma}
\label{lem:hyp-behaviour}
There exists $\sigma_0 >0$ so that, for any $0 < \sigma \leq \sigma_0$, the following holds:

 There exist constants  $E$  and $C$, such that  for any $R > \delta >0$ and for any $\varepsilon \in [0,1]$ $z=(x,y,c) \in D(R)$, holds
  \begin{itemize}
    \item if $\|x\| \geq \|y\|$, $(x,y) \neq 0$ and $\varphi(t,x,y,c) \in D(R)$ for $t \in [0,T]$, then $\|\pi_x\varphi(t,x,y,c)\| > \|\pi_y\varphi(t,x,y,c)\|$ for $t \in (0,T]$
      and
      \begin{equation}
      \frac{d}{dt}\|x(t)\|  \geq E \|x(t)\|, \quad t \in (0,T].  \label{eq:xCone-growth}
      \end{equation}

      In particular, $T \leq \frac{1}{E}\ln \left( \frac{R}{\|x(0)\|} \right)$.
    \item if   $\|x(t)\| \leq \|y(t)\|$, and $\varphi(t,x,y,c) \in D(R)$ for $t \in [0,T]$, then  $\|y(t)\|\leq \|y(0)\|e^{-tC}$.

    In particular, if $t> \frac{1}{E}\ln \left( \frac{\|y(0)\|}{\delta} \right)$, then $\varphi(t,x,y,c) \in D(\delta)$.
  \end{itemize}
\end{lemma}
\textbf{Proof:}
To prove the first assertion we need to establish two facts. First, the forward invariance of the cone $\|x\| \geq \|y\|$ and
second, the expansion in this cone.
We have  for $\|x\| \geq \|y\|$ (we use (\ref{eq:eucl-log-norm},\ref{eq:eucl-m-l-norm})) (compare to the proof of Lemma~\ref{lem:iso-block-mod-vekf}).

We have
\begin{eqnarray*}
 \frac{1}{2}\frac{d \|x\|^2}{dt} =(\dot{x},x)=(A_u x,x) + \left(\widetilde{N}(x,y,c) \varepsilon \widetilde{\mathcal{R}}(x,y,c),x\right)  \\
 \geq m_l(A_u) \|x\|^2 -\left(\|D\widetilde{N} +  \|D\widetilde{\mathcal{R}}\|\right) \cdot \|(x,y)\|\cdot \|x\|
\end{eqnarray*}
and
\begin{eqnarray*}
  \frac{1}{2}\frac{d \|y\|^2}{dt}=(\dot{y},y) =(A_s y,y) + \left(\widetilde{N}(x,y,c) + \varepsilon \widetilde{\mathcal{R}}(x,y,c),y\right) \\
    \leq \mu_{log}(A_s) \|y\|^2 +\left(\|D\widetilde{N} +  \|D\widetilde{\mathcal{R}}\|\right) \cdot \|(x,y)\|\cdot \|x\|.
\end{eqnarray*}
Therefore for $(x,y,c)$ such that $\|x\| \geq \|y\|$, $x\neq 0$   (we use  (\ref{eq:DtildeN-estm}, \ref{eq:DtilR-estm})) holds
\begin{eqnarray*}
  \frac{1}{2} \frac{d}{dt}  \left(\|x\|^2 - \|y\|^2 \right) \geq  m_l(A_u) \|x\|^2 - \mu_{log}(A_s) \|y\|^2  - 2\left(\|D\widetilde{N} +  \|D\widetilde{\mathcal{R}}\|\right) \cdot \|(x,y)\|\cdot \|x\| \\
  \geq \left(m_l(A_u) -  \mu_{log}(A_s) + 2 \sqrt{2}\left(K_N'\sigma + K_\mathcal{R}\sigma^Q \right) \right) \|x\|^2 >0
\end{eqnarray*}
if $\sigma \leq \sigma_0$.  Therefore the cone $\|x\| \geq \|y\|$ is forward invariant.

Now we deal with the expansion in the $x$-direction. From previous computations we obtain in the cone $\|x\|\geq \|y\|$
\begin{eqnarray*}
 \frac{1}{2}\frac{d \|x\|^2}{dt}  \geq  m_l(A_u) \|x\|^2 -\sqrt{2}(K_N'\sigma + K_\mathcal{R}\sigma^Q) \|x\|^2
 \geq E \|x\|^2
\end{eqnarray*}
where $ E = m_l(A_u)  -\sqrt{2}(K_N'\sigma + K_\mathcal{R}\sigma_0^Q) >0$. We chose $\sigma_0$, small enough for $E$ to be positive.  Observe
that this implies  (\ref{eq:xCone-growth}).

Now we turn to second assertion about the decay in the cone $\|x\| \leq \|y\|$.  In this cone from previous computations we obtain
\begin{eqnarray*}
  \frac{1}{2}\frac{d \|y\|^2}{dt}\leq \left(\mu_{log}(A_s) +\sqrt{2}K'_{N} \sigma  + \sqrt{2}K_{\mathcal{R}} \sigma^Q  \right) \|y\|^2.
\end{eqnarray*}
Since $\mu_{log}(A_s)<0$ we can find $\sigma_0$, such that $C=-\mu_{log}(A_s) -\sqrt{2}\left(K'_{N} \sigma  + \sqrt{2}K_{\mathcal{R}}\right) > 0$. This implies the second assertion.

\qed

We fix $0<\sigma \leq \sigma_0$ so that the assertions of  Lemmas~\ref{lem:iso-block-mod-vekf} and~\ref{lem:hyp-behaviour} hold true and from now on we work with the modified vector field, but dropping tildes.


We will need some extra  bounds for the functions $y^\textrm{u}$ and $x^\textrm{s}$ as well as for their derivatives up to order two. It is easy to see that if $B$ does not contain any nontrivial Jordan blocks, then $y^\textrm{u}(x,c)=0$ and $x^\textrm{s}(y,c)=0 $ if $|c| > O(\sigma r_1)$. This is an immediate consequence of the following fact: as $|\pi_c \varphi(t,x,y,c)=e^{Bt}c| > O(\sigma r_1)$ for $t \in \mathbb{R}$, the trajectory remains in the domain where the dynamics is just linear; $x'=A_u x$, $y'= A_s y$, $c'=Bc$.

However, we  want to include the situation with nontrivial Jordan blocks in the center direction and we also  need information about the derivatives of $y^\textrm{u}$ and $x^\textrm{s}$
with respect to $x$ and $y$, respectively, when $|x| \to \infty$ and $|y| \to \infty$.

\begin{lemma}
\label{lem:bnds-yu-xs}
In the context of Lemma~\ref{lem:iso-block-mod-vekf},
\begin{eqnarray}
  |y^\textrm{u}(x,c)| \leq    L \sigma r_1 \min\left(1,\left(\frac{|x|}{\sigma r_1}\right)^{\frac{\mu_{log}(A_s)}{\mu_{log}(A_u)}}\right)  , \quad \forall (x,c), \label{eq:yu-bnd} \\
  |x^\textrm{s}(y,c)| \leq    L \sigma r_1 \min\left(1,\left(\frac{|y|}{\sigma r_1}\right)^{\frac{\mu_{log}(-A_u)}{\mu_{log}(-A_s)}}\right)  , \quad \forall (y,c). \nonumber 
\end{eqnarray}

Moreover, for some $w_c>1$ and $w_z >0$ and $w=\min(w_c,w_z)$,

\begin{eqnarray}
\left|\frac{\partial y^\textrm{u}}{\partial x}(x,c)\right| \leq \min(L,L O(\sigma^{w}) (|c|^{w_c} + |x|^{w_z})^{-1}),  \label{eq:dyudx-bnd}  \\
  \left|\frac{\partial x^\textrm{s}}{\partial y}(y,c)\right| \leq \min(L,L O(\sigma^{w}) (|c|^{w_c} + |y|^{w_z})^{-1}),  \nonumber
\end{eqnarray}
and
\begin{eqnarray*}
  \left|\frac{\partial^2 y^\textrm{u}}{\partial^2 x}(x,c)\right| \leq \min(M,O(\sigma^w) (|c|^{w_c} + |x|^{w_z})^{-1}) )  \\
  \left|\frac{\partial^2 x^\textrm{s}}{\partial^2 y}(y,c)\right| \leq \min(M,O(\sigma^w) (|c|^{w_c} + |y|^{w_z})^{-1}) ).
\end{eqnarray*}
\end{lemma}
\textbf{Proof:}
It is enough to consider just $y^\textrm{u}$, since the argument for $x^\textrm{s}$ is analogous,
just reversing the time direction and changing $x$ to $y$.

Let us denote by $S$ the support of the nonlinearities in system  (\ref{eq:dotx-lcd-app2}--\ref{eq:dotc-lcd-app2}). We know that
\begin{equation*}
  S \subset \overline{B}_{n_u}(0,\sigma r_1) \times  \overline{B}_{n_s}(0,\sigma r_1)\times \overline{B}_{2m}(0,\sigma r_1).
\end{equation*}
Let $(x,y^{\textrm{u}}(x,c),c) \in W^{cu}$ be such that in the $x$-direction we have $|x|=\sigma r_1$. Then $|y^{\textrm{u}}(x,c)| \leq L \sigma r_1$ and
its image forward in time does not enter $S$, hence is equal to $(e^{t A_u}x,e^{t A_s}y^{\textrm{u}}(x,c),e^{tB}c) \in W^{cu}$, which means that
\begin{eqnarray*}
  e^{t A_s}y^{\textrm{u}}(x,c)=y^{\textrm{u}}(e^{t A_u}x,e^{tB}c).
\end{eqnarray*}
Since the graph transform (used in proving the existence of $W^{cu}$) of $W^{cu}$ is equal to  $W^{cu}$, therefore for $(\bar{x},\bar{c})$ with $|\bar{x}|> \sigma r_1$, there exists $(x,c)$ such that
\begin{equation*}
  (e^{t A_u}x,e^{t A_s}y^{\textrm{u},e^{tB}c,}(x,c))=(\bar{x},y^{\textrm{u}}(\bar{x},\bar{c})),
\end{equation*}
and therefore we get (because $\|e^{A_s t}\| \leq e^{\mu_{log}(A_s)t} < 1$ for $t>0$)
\begin{equation*}
  |y^\textrm{u}(x,c)| \leq L \sigma r_1, \quad \forall (x,c).  
\end{equation*}
On the other hand, for $|x| > \sigma r_1$ the time to reach $|x|$ from $\sigma r_1$ is estimated from below by \newline
 $\displaystyle \frac{1}{\mu_{log}(A_u)} \ln \left(\frac{|x|}{\sigma r_1} \right)$, so we get
\begin{equation*}
  |y^\textrm{u}(x,c)| \leq L \sigma r_1 \left(\frac{|x|}{\sigma r_1}\right)^{\frac{\mu_{log}(A_s)}{\mu_{log}(A_u)}}.
\end{equation*}

Now we proceed to find the estimates for the derivatives of $y^\textrm{u}$.
It follows immediately from  (\ref{eq:Lip-yu}) in Lemma~\ref{lem:iso-block-mod-vekf} that
\begin{equation*}
  \left|\frac{\partial y^{\textrm{u}}}{\partial c}(x,c) \right|,  \left|\frac{\partial y^{\textrm{u}}}{\partial x}(x,c) \right| \leq L, \quad \forall (x,c). 
\end{equation*}

Consider a point  $p=(x,y^\textrm{u}(x,c),c)  \in W^{cu}$, $p \notin S$. Let $t(p)>0$ be such that for $t \in [0,t(p)]$ $\varphi(-t,p) \notin S$. Then
$\varphi(-t,p)=(e^{-tA_u}x,e^{-tA_s}y^\textrm{u}(x,c),e^{-tB}c) \in W^{cu}$, therefore
\begin{eqnarray*}
  e^{-t(p)A_s}y^\textrm{u}(c,x)=y^\textrm{u}(e^{-t(p)A_u}x,e^{-t(p)B}c),
\end{eqnarray*}
and finally
\begin{equation}
  y^\textrm{u}(x,c)= e^{t(p)A_s}y^\textrm{u}(e^{-t(p)A_u}x,e^{-t(p)B}c). \label{eq:yu-rel}
\end{equation}
Observe that $t(p)$ is also good for all points in the neighborhood of $p$, so when differentiating the above equation we can treat $t(p)$ as a constant, if needed.

Depending on whether $|c|$ or $|x|$ are larger than $\sigma r_1$ we can take the following expressions $t_{x}(p)$ or $t_c(p)$ for $t(p)$
\begin{equation*}
  t_{x}(p)= \frac{1}{\mu_{log}(A_u)} \ln \left(\frac{|x|}{ \sigma r_1} \right), \quad t_c(p)= \frac{1}{\mu_{log}(B)} \ln \left(\frac{|c|}{\sigma r_1} \right).
\end{equation*}

By differentiation of  (\ref{eq:yu-rel})
\begin{eqnarray*}
   \frac{\partial y^\textrm{u}}{\partial x}(x,c)&=& e^{t(p)A_s} \left(\frac{\partial y^\textrm{u}}{\partial c}\left(e^{-t(p)A_u}x,e^{-t(p)B}c\right)\right) \cdot e^{-t(p)A_u},
\end{eqnarray*}
hence
\begin{eqnarray*}
  \left|\frac{\partial y^\textrm{u}}{\partial x}(x,c) \right| &\leq& L e^{t(p) (\mu_{log}(A_s) + \mu_{log}(-A_u) )}.
\end{eqnarray*}
For the second derivatives we obtain the following bounds (we use (\ref{eq:Dwcscu-bnd}))
\begin{eqnarray*}
   \left|\frac{\partial^2 y^\textrm{u}}{\partial^2 x}(x,c) \right| &\leq& K_W e^{t(p) (\mu_{log}(A_s) + 2\mu_{log}(-A_u))}.
\end{eqnarray*}
Observe that all exponents $w$ appearing in $e^{t(p) \cdot w}$ are negative (if $\mu_{log}(-B)$ is small enough).

Let us evaluate exponentials of the type $e^{t(p) \cdot w}$ for $t(p)=t_{x}(p), t_{c}(p)$.
We have
\begin{eqnarray*}
  e^{t_{x}(p) w}&=&(\sigma r_1)^{-\frac{w}{\mu_{log}(A_u)}} |x|^\frac{w}{\mu_{log}(A_u)}, \\
  e^{t_c(p) w}&=&(\sigma r_1)^{-\frac{w}{\mu_{log}(B)}} |c|^\frac{w}{\mu_{log}(B)},
\end{eqnarray*}
since for each exponent $w \leq \mu_{log}(A_s) + 2 \mu_{log}(-B)<0$, hence $w_c=-\frac{w}{\mu_{log}(B)} >1$ provided $\mu_{log}(B)$ is small enough.

Regarding the decay with increasing $|c|$ we see that the first order derivatives  will be less than or equal to $LO(\sigma^{w_c})|c|^{-w_c}$, and for the second order we have a bound  $O(\sigma^{w_c})|c|^{-w_c}$.

Let us see now what will be the decay in $x$ direction.
\begin{itemize}
\item  $\displaystyle \frac{\partial y^\textrm{u}}{\partial x}(x,c) $. We have $w=\mu_{log}(A_s) + \mu_{log}(-A_u)$, so
$\displaystyle w_z=\frac{-\mu_{log}(A_s)}{\mu_{log}(A_u)} - \frac{\mu_{log}(-A_u)}{\mu_{log}(A_u)} $,
\item $\displaystyle \frac{\partial^2 y^\textrm{u}}{\partial x_i x_j}(x,c) $. We have $w=\mu_{log}(A_s) + 2\mu_{log}(-A_u)$, so
$\displaystyle w_z=\frac{-\mu_{log}(A_s)}{\mu_{log}(A_u)} - 2\frac{\mu_{log}(-A_u)}{\mu_{log}(A_u)} $.
\end{itemize}
In both cases it can happen that $w_z <1$.   For example if $\mu_{\mathrm{min}} < \dots < \mu_{\mathrm{max}}$ are eigenvalues of $A_u$, then $-\mu_{log}{-A_u}=\mu_{\mathrm{min}}$, and $w_z=\lambda_{\mathrm{min}}/\mu_{\mathrm{max}} + \mu_{\mathrm{min}}/\mu_{\mathrm{max}}$.

\qed

\subsection{Straightening the  invariant manifolds}
\label{subsec:straight-man}
We continue with the modified vector field obtained in Section~\ref{subsec:prep-compt-data}.
The center-unstable and center-stable  invariant manifolds    can be straightened in suitable coordinates.
The following transformation does this
  \begin{equation}
    T_\varepsilon(x,y,c)=(x-x^{\textrm{s}}(\varepsilon,y,c),y-y^{\textrm{u}}(\varepsilon,x,c),c).  \label{eq:S-man-tran}
   \end{equation}

\begin{lemma}
\label{lem:T-str-diffeo}
Under the same assumptions as in Lemma~\ref{lem:iso-block-mod-vekf}, $T_\varepsilon(x,y,c)$ is $C^{q-1}$ and for every $\varepsilon \in [0,1]$ the transformation $T_\varepsilon: \PhSpace  \to \PhSpace$ is a diffeomorphism.
Moreover,
\begin{equation}
 T_\varepsilon^{-1}(x,y,c)=(x+O(L\sigma), y+O(L\sigma),c). \label{eq:TT-inv-close-Id}
\end{equation}
\end{lemma}
\textbf{Proof:}
The regularity follows from Lemma~\ref{lem:iso-block-mod-vekf}.  We need to show that $T_\varepsilon$ is a bijection.  To show that it is onto we fix $c$ and observe that
for $\|x\|=R$ and $\|y\| \leq R$ it holds that (we use (\ref{eq:Lip-xs}))
\begin{eqnarray*}
  \|\pi_x T_\varepsilon (x,y,c)\| \geq \|x\| - \|x^{\textrm{s}}(\varepsilon,y,c)\| \geq R - L R = (1-L)R,
\end{eqnarray*}
and analogously for $\|x\|\leq R$ and $\|y\| = R$, from (\ref{eq:Lip-yu}) we get
\begin{eqnarray*}
  \|\pi_y T_\varepsilon (x,y,c)\| \geq (1-L)R.
\end{eqnarray*}
From this inequality, using the local Brouwer degree argument, we obtain that
\begin{equation*}
  \overline{B}_{n_u}((1-L)R) \times   \overline{B}_{n_s}((1-L)R) \subset \pi_{x,y} T_\varepsilon(\overline{B}_{n_u}(R), \overline{B}_{n_s}(R),c).
\end{equation*}
Therefore $T_\varepsilon$ is onto. 

Now we prove the injectivity. Take any $(x_1,y_1) \neq (x_2,y_2)$.  We have (we use (\ref{eq:Lip-xs},\ref{eq:Lip-yu}))
\begin{eqnarray*}
  \| \pi_x \left(T_\varepsilon(x_1,y_1,c) - T_\varepsilon(x_2,y_2,c) \right) \| \geq \|x_1 - x_2\| - L\|y_1-y_2\|, \\
   \| \pi_y \left(T_\varepsilon(x_1,y_1,c) - T_\varepsilon(x_2,y_2,c) \right) \| \geq \|y_1 - y_2\| - L\|x_1-x_2\|.
\end{eqnarray*}
Since $L<1$, if $\|x_1-x_2\| \geq \|y_1-y_2\|$, then $ \|x_1 - x_2\| - L\|y_1-y_2\| >0$ and analogously in the other case. Hence $T_\varepsilon$ is an injection.


Condition (\ref{eq:TT-inv-close-Id}) follows immediately from (\ref{eq:yu-bnd}) in Lemma~\ref{lem:bnds-yu-xs}.

\qed


\subsubsection{Cohomological equation after straightening}

Observe that the transformation $T_\varepsilon$ given by (\ref{eq:S-man-tran}) gives new variables in terms of the old ones, therefore to transform vector fields and Hamiltonians we use the pushforward $T_{\varepsilon \ast}$. We will often drop $\varepsilon$ in the sequel.

Cohomological equation (\ref{eq:coh-flows})   becomes
\begin{equation*}
 [ T_{\varepsilon \ast} \mathcal{Z},  T_{\varepsilon \ast}\mathcal{G}]=-T_{\varepsilon \ast} \mathcal{R}
\end{equation*}
or in the Hamiltonian case,  (\ref{eq:coh-ham}) becomes
\begin{equation*}
  (D (T_{\varepsilon \ast} G)) \cdot  T_{\varepsilon \ast} \mathcal{Z} = T_{\varepsilon \ast} R.
\end{equation*}

Hence the cohomological equation retains its form, but the vector field $\mathcal{Z}$ and the remainder $\mathcal{R}$ (or $R$)  change.

In the Hamiltonian case $T$ might not be not symplectic in the presence of center variables. Therefore the transformed symplectic
form  $T_\ast \Omega$ no longer has the standard form. This influences the relationship between the Hamiltonian function and the induced vector field.

\subsubsection{The vector field after straightening}

\begin{lemma}
\label{lem:vf-afterstr}
For a general vector field, $T_{\varepsilon \ast} \mathcal{Z}_\varepsilon$ and $T_{\varepsilon \ast} \mathcal{R}$ are $C^{q-2}$ as functions of $(\varepsilon,x,y,c)$,
whereas in the Hamiltonian case  $T_{\varepsilon \ast} Z_\varepsilon$ and $T_{\varepsilon \ast} R$ are $C^{q-1}$.

The vector field $\mathcal{Z}_\varepsilon$ after straightening, i.e. $T_{\varepsilon \ast} \mathcal{Z}_\varepsilon$, becomes
\begin{eqnarray}
  \dot{x}&=&g_x(\varepsilon,x,y,c)x,  \label{eq:dotx-after-str} \\
  \dot{y}&=&g_y(\varepsilon,x,y,c)y, \label{eq:doty-after-str} \\
  \dot{c}&=&Bc + M_c(\varepsilon ,x,y,c),  \label{eq:dotc-after-str}
\end{eqnarray}
where   $g_x,g_y \in C^{q-3}$,  $M_c$ has compact support and
\begin{eqnarray}
 M_c(\varepsilon,x,y,c)&=&O_2(|(x,y,c)|^2),  \label{eq:Mcstr} \\
  (g_{x},g_{y})(\varepsilon,x,y,c)&=&(A_u,A_s) +  O(L) + O(\sigma) + (O(L)y,O(L)x). \label{eq:gxystr}
\end{eqnarray}
All  $O()$ terms are uniform with respect to $\varepsilon \in [0,1]$ in (\ref{eq:Mcstr}) and with respect to $(\varepsilon, x,y,c) \in [0,1] \times \PhSpace$.

Moreover,
\begin{equation}
T_{\varepsilon \ast} \mathcal{R}(x,y,c)=O_{Q+1}(|(x,y)|^{Q+1}),  \label{eq:T*RO}
\end{equation}
or $T_{\varepsilon \ast} R(x,y,c)=O_{Q+1}(|(x,y)|^{Q+1})$ in the Hamiltonian case.
\end{lemma}
\textbf{Proof:}
Recall that in Theorem~\ref{thm:polNForm-loccenter} we assumed that
\begin{equation}
  q \geq Q+3.  \label{eq:qgeqQ3}
\end{equation}

The regularity of $T_{\varepsilon \ast} \mathcal{Z}_\varepsilon$ and $T_{\varepsilon \ast} \mathcal{R}$ is an immediate consequence of Lemma~\ref{lem:T-str-diffeo}.

The factorization of $x$ in (\ref{eq:dotx-after-str}) and $y$ in (\ref{eq:doty-after-str}) is due to the fact that in the new coordinates $W^{cs}=\{x=0\}$ and $W^{cu}=\{y=0\}$, respectively. For this we need that the transformed vector field is at least of class $C^1$, so we need $q-2 \geq 1$, which is granted by (\ref{eq:qgeqQ3}).

The straightening implies that in the new variables we obtain the equations
\begin{eqnarray*}
  \dot{x}&=&f_x(\varepsilon,x,y,c)= \left(\int_0^1 \frac{\partial f_x}{\partial x}(\varepsilon,c,t x,y)dt\right) \cdot x,   \\
  \dot{y}&=&f_y(\varepsilon,x,y,c)= \left(\int_0^1 \frac{\partial f_y}{\partial y}(\varepsilon,c,x,ty)dt\right) \cdot y,  \\
  \dot{c}&=&f_c(\varepsilon,x,y,c),
\end{eqnarray*}
where
\begin{equation*}
  (f_x,f_y,f_c)(\varepsilon,x,y,c)=DT_\varepsilon(T_\varepsilon^{-1}(x,y,c)) \mathcal{Z}_\varepsilon(T_\varepsilon^{-1}(x,y,c)).
\end{equation*}
In the sequel we drop $\varepsilon$ because we have bounds which are valid for all $\varepsilon$.

 Our first goal is the computation
of   $\displaystyle \frac{\partial }{\partial x}\left(\frac{\partial T_x}{\partial c}(T^{-1}(x,y,c))\right)$,
 $\displaystyle \frac{\partial }{\partial x}\left(\frac{\partial T_x}{\partial x}(T^{-1}(x,y,c))\right)$
 as well as of $\displaystyle \frac{\partial }{\partial x}\left(\frac{\partial T_x}{\partial y}(T^{-1}(x,y,c))\right)$.

From (\ref{eq:Lip-yu},\ref{eq:Lip-xs}) in Lemma~\ref{lem:iso-block-mod-vekf} it follows that
\begin{equation}
  \left\| \frac{\partial x^\textrm{s}}{\partial y} \right\|,  \left\| \frac{\partial x^\textrm{s}}{\partial c} \right\|,  \left\| \frac{\partial y^\textrm{u}}{\partial x} \right\|,  \left\| \frac{\partial y^\textrm{u}}{\partial c} \right\| \leq L. \label{eq:derxuysbnd}
\end{equation}
Recall that  $L$ is a parameter, which can be chosen arbitrarily small.

From (\ref{eq:derxuysbnd}) and (\ref{eq:S-man-tran}) it follows that  for any $(x,y,c)$ it holds that
\begin{eqnarray}
  DT = \left[ \begin{array}{ccc}
          I_{x}, & O(L), & O(L) \\
          O(L)  & I_y, & O(L) \\
           0, & 0, & I_c
       \end{array} \right], \qquad
     DT^{-1} = \left[ \begin{array}{ccc}
         I_{x}+O(L), & O(L), & O(L) \\
          O(L)  & I_y+O(L), & O(L) \\
           0, & 0, & I_c
       \end{array} \right],   \label{eq:DT-DT-1}
\end{eqnarray}
where by $I_c,I_{x},I_y$ we denote the identity matrix acting on $c$,$x$,$y$ variables, respectively.

From Lemma~\ref{lem:bnds-yu-xs} we can obtain  better bounds for $\displaystyle \frac{\partial T^{-1}_y}{\partial x}(x,y,c)$ (and $\displaystyle \frac{\partial T^{-1}_x}{\partial y}(x,y,c)$).
Let us denote, for $L,\sigma,c_1,z_1 \in \mathbb{R}$,
\begin{equation*}
  E(L,\sigma,c_1,z_1)= \min(L,L \sigma^{w_z} (c_1^{w_c} + z_1^{w_z})^{-1}),
\end{equation*}
where $w_z$ and $w_c$ are as Lemma~\ref{lem:bnds-yu-xs}.
Then from  (\ref{eq:dyudx-bnd}) in Lemma~\ref{lem:bnds-yu-xs}  we have, using the formula $(I-N)^{-1}=I+N + N^2 + \cdots$ with $\displaystyle \frac{\partial T_{x,y}}{\partial (x,y)}=I_{x,y} +N $ and $\displaystyle N=\left[\begin{array}{cc} 0 & \frac{\partial x^s}{\partial y}  \\ \frac{\partial y^u}{\partial x} & 0  \end{array}\right]$,
\begin{eqnarray}
  \left|\frac{\partial T^{-1}_y}{\partial x}(x,y,c)\right|= O(E(L,\sigma,|c|,|x|)).  \label{eq:dT-1xdy}
\end{eqnarray}

We  compute $\frac{\partial }{\partial x}\left(\frac{\partial T_x}{\partial c}(T^{-1}(x,y,c))\right)$. From (\ref{eq:S-man-tran}) we see that $\frac{\partial T_c^{-1}}{\partial x}=0$ and $\frac{\partial^2 T_x}{\partial c \partial x}=0$, hence from (\ref{eq:Dwcscu-bnd}), (\ref{eq:DT-DT-1}) and (\ref{eq:dT-1xdy}) we obtain
\begin{eqnarray*}
  \frac{\partial }{\partial x}\left(\frac{\partial T_x}{\partial c}(T^{-1}(x,y,c))\right) &=& \frac{\partial^2 T_x}{\partial c \partial c}(T^{-1}(x,y,c)) \cdot \frac{\partial T_c^{-1}}{\partial x}(x,y,c)
  + \frac{\partial^2 T_x}{\partial c \partial x}(T^{-1}(x,y,c)) \cdot \frac{\partial T_x^{-1}}{\partial x}(x,y,c)\\
   &&+ \frac{\partial^2 T_x}{\partial c \partial y}(T^{-1}(x,y,c)) \cdot \frac{\partial T_y^{-1}}{\partial x}(x,y,c) \nonumber \\
  &=&-\frac{\partial^2 x^\textrm{s}}{\partial c \partial y}(T_{y,c}^{-1}(x,y,c)) \cdot \frac{\partial T^{-1}_y}{\partial x}(x,y,c)\\
  &=&K_W O(E(L,\sigma,|c|,|x|)) = O(E(L,\sigma,|c|,|x|)).
\end{eqnarray*}

Using (\ref{eq:Dwcscu-bnd}) and (\ref{eq:dT-1xdy}) we obtain
\begin{eqnarray*}
  \frac{\partial }{\partial x}\left(\frac{\partial T_x}{\partial y}(T^{-1}(x,y,c))\right)=\frac{\partial^2 T_x}{\partial y \partial c}(T^{-1}(x,y,c)) \cdot \frac{\partial T_c^{-1}}{\partial x}(x,y,c) (=0)\\
  + \frac{\partial^2 T_x}{\partial y \partial x}(T^{-1}(x,y,c)) (=0) \cdot \frac{\partial T_x^{-1}}{\partial x}(x,y,c)
   + \frac{\partial^2 T_x}{\partial y \partial y}(T^{-1}(x,y,c)) \cdot \frac{\partial T_y^{-1}}{\partial x}(x,y,c) \\
  =-\frac{\partial^2 x^\textrm{s}}{\partial y \partial y}(T_{y,c}^{-1}(x,y,c)) \cdot \frac{\partial T^{-1}_y}{\partial x}(x,y,c)
  =K_W O(E(L,\sigma,|c|,|x|))= O(E(L,\sigma,|c|,|x|)).
\end{eqnarray*}

Summarizing, we get
\begin{eqnarray}
  \frac{\partial }{\partial x}\left(\frac{\partial T_x}{\partial c}(T^{-1}(x,y,c))\right)&=& O(E(L,\sigma,|c|,|x|)), \label{eq:ddxTxc} \\
  \frac{\partial }{\partial x}\left(\frac{\partial T_x}{\partial x}(T^{-1}(x,y,c))\right) &=& 0,  \label{eq:ddxTxx} \\
   \frac{\partial }{\partial x}\left(\frac{\partial T_x}{\partial y}(T^{-1}(x,y,c))\right) &=& O(E(L,\sigma,|c|,|x|)).\label{eq:ddxTxy}
\end{eqnarray}

In order to compute $T_\ast \mathcal{Z}(x,y,c))$ we will treat separately the linear part of $\mathcal{Z}$ and the nonlinear one given by $\widetilde{M}$.
The linear part is given by $\mathcal{L}(x,y,c)=(A_u x, A_s y,Bc)$, and it is immediate that
\begin{equation*}
  (T_\ast \mathcal{L})_c(x,y,c)=Bc.
\end{equation*}

For the $x$-component (for $y$-component the computations are analogous) we get
\begin{eqnarray*}
  (T_\ast \mathcal{L})_x(x,y,c)&=& \frac{\partial T_x}{\partial c}(T^{-1}(x,y,c)) Bc
 +  \frac{\partial T_x}{\partial x}(T^{-1}(x,y,c))(A_u T^{-1}_x(x,y,c))\\
   &&+ \frac{\partial T_x}{\partial y}(T^{-1}(x,y,c))(A_s T^{-1}_y(x,y,c)).
\end{eqnarray*}
 We are interested in the derivative of $(T_\ast \mathcal{L})_x(x,y,c)$ with respect to $x$ . Let us investigate each term of the above sum separately.
 For the first term we have following  estimate (we use (\ref{eq:ddxTxc}))
\begin{eqnarray*}
  \frac{\partial }{\partial x} \left( \frac{\partial T_x}{\partial c}(T^{-1}(x,y,c)) Bc \right) = O(|c| E(L,\sigma,|c|,|x|))=O(L).
\end{eqnarray*}

For the second term we have (we use (\ref{eq:DT-DT-1},\ref{eq:ddxTxx}))
\begin{eqnarray*}
  \frac{\partial }{\partial x} \left(\frac{\partial T_x}{\partial x}(T^{-1}(x,y,c))(A_u (T^{-1}_x(x,y,c)))\right)
  =A_u \frac{\partial T_x^{-1}}{\partial x}(x,y,c)=A_u + O(L).
\end{eqnarray*}

For the third term we get
\begin{eqnarray*}
  \frac{\partial }{\partial x} \left( \frac{\partial T_x}{\partial y}(T^{-1}(x,y,c))(A_s T^{-1}_y(x,y,c)) \right)&=&
  \frac{\partial }{\partial x} \left( \frac{\partial T_x}{\partial y}(T^{-1}(x,y,c))\right)\cdot A_s T^{-1}_y(x,y,c) \\
   &&+  \frac{\partial T_x}{\partial y}(T^{-1}(x,y,c)) \cdot A_s \frac{\partial T^{-1}_y}{\partial x}(x,y,c),
\end{eqnarray*}
and from  (\ref{eq:ddxTxy}) and (\ref{eq:TT-inv-close-Id}) we have

\begin{eqnarray*}
   \frac{\partial }{\partial x} \left( \frac{\partial T_x}{\partial y}(T^{-1}(x,y,c))\right)\cdot A_s T^{-1}_y(x,y,c)
    &=& O(E(L,\sigma,|c|,|x|)) A_s (y+O(\sigma L)) =O(L)\cdot y + O(\sigma^2 L^2),\\
  \frac{\partial T_x}{\partial y}(T^{-1}(x,y,c)) \cdot A_s \frac{\partial T^{-1}_y}{\partial x}(x,y,c) &=&  - \frac{\partial x^\textrm{s}}{\partial y}(T^{-1}_{c,y}(x,y,c))A_s O(L) =O(L^2)
\end{eqnarray*}

Therefore we obtain
\begin{eqnarray*}
 \frac{\partial }{\partial x} \left( (T_\ast \mathcal{L})_x(x,y,c)\right)=O(L) + (A_u + O(L)) + (O(L) y   + O(L^2))
 = A_u + O(L)y + O(L).
\end{eqnarray*}

From bounds (\ref{eq:tildeN-estm},\ref{eq:DtildeN-estm},\ref{eq:DtilR-estm}) it easily follows that
\begin{eqnarray*}
  T_\ast \widetilde{M}(x,y,c) = O(\sigma^2), \\
  D T_\ast \widetilde{M}(x,y,c) = O(\sigma),
\end{eqnarray*}
and $T_\ast \widetilde{M}(x,y,c))$ has compact support.

This finishes the proof of (\ref{eq:dotx-after-str},\ref{eq:doty-after-str},\ref{eq:dotc-after-str},\ref{eq:Mcstr},\ref{eq:gxystr}).

It remains to prove (\ref{eq:T*RO}).
We will start with the function $R$, because the result is valid also for components of the vector field $\mathcal{R}$.

From   (\ref{eq:S-man-tran}) and inequalities (\ref{eq:Dwcscu-bnd}) in Lemma~\ref{lem:iso-block-mod-vekf} one can easily infer that there exists a constant $K_1$, such that
\begin{eqnarray}
  \|D^j T\|\leq K_1, \quad   \|D^jT^{-1}\| \leq K_1, \quad j=1,\dots, q-1. \label{eq:DjT-1-estm}
\end{eqnarray}

Since $T^{-1}(c,0,0)=(c,0,0)$ we obtain from (\ref{eq:DT-DT-1}) that
\begin{equation*}
  T^{-1}(x,y,c)=(x+O(L)\cdot(x,y),y+O(L)\cdot (x,y),c),  
\end{equation*}
hence there exists a constant $K_2$, such that
\begin{eqnarray}
  \|\pi_{x,y}T^{-1}(x,y,c)\| \leq K_2 \|(x,y)\|. \label{eq:T-inv-bnd}
\end{eqnarray}
From this it is immediate that
\begin{eqnarray*}
  |T_\ast R(x,y,c)| \leq K  \|\pi_{x,y}T^{-1}(x,y,c)\|^{Q+1} \leq K K_2^{Q+1} \|(x,y)\|^{Q+1}.
\end{eqnarray*}
Observe that the derivatives  of $T_\ast R$ are of the following form, for $j=1,\dots,Q+1$,
\begin{eqnarray*}
  D^j T_\ast R &=& D^j R() (DT^{-1})^j + D^{j-1} R() (DT^{-1})^j (D^2T^{-1}) + \dots + DR() (D^jT^{-1}).
\end{eqnarray*}
Using (\ref{eq:DjT-1-estm}) and (\ref{eq:T-inv-bnd}) we see that (\ref{eq:T*RO})  holds for the function $R$.

From this we also have that in the case of the  vector field $\mathcal{R}$
\begin{equation*}
  \mathcal{R}(T^{-1}(x,y,c)) = O_{Q+1}(|(x,y)|^{Q+1}).
\end{equation*}
The $x$-component of $T_\ast \mathcal{R}$ is given by (recall that $\mathcal{R}_c=0$)
\begin{eqnarray*}
  (T_\ast \mathcal{R})_x (x,y,c) =  \mathcal{R}(T^{-1}(x,y,c))_x +  \frac{\partial T_x}{\partial y}(T^{-1}(x,y,c)) \mathcal{R}(T^{-1}(x,y,c))_y.
\end{eqnarray*}
Observe that due to (\ref{eq:DjT-1-estm}) the  function $\frac{\partial T_x}{\partial y}(T^{-1}(x,y,c))$  has bounded derivatives up to order $q-2$ (condition (\ref{eq:qgeqQ3}) gives us the needed regularity). Hence we can infer in the same way as in the case of functions that $(T_\ast \mathcal{R})_x (x,y,c)=O_{Q+1}(|(x,y)|^{Q+1})$. The same holds for the $y$-component. This finishes the proof of (\ref{eq:T*RO}) for the vector field $\mathcal{R}$.
\qed

\begin{definition}
\label{def:chi-tilr1}
Let $\tilde{r}_1$ be such that the support of $T_\ast \mathcal{R}$  (or $T_\ast R$ in Hamiltonian case)  is contained in $E(\tilde{r}_1)=\overline{B}_{n_u}(0,\tilde{r}_1) \times \overline{B}_{n_s}(0,\tilde{r}_1) \times \overline{B}_{2m}(0,\tilde{r}_1)$ for all $\varepsilon \in [0,1]$. Then we denote by $\chi$ the \emph{indicator function} for
$E(\tilde{r}_1)$: $\chi(z)=1$ if $z \in D(\tilde{r}_1)$ and $\chi(z)=0$ otherwise.
\end{definition}

\subsubsection{Consequences of straightening of invariant manifolds.}

Let $\mathcal{Z}$ be the vector field obtained after straightening and let  $\varphi(t,z)$ be the local flow induced by it.
Let  us define
\begin{equation*}
  A_u(\varepsilon,c)=g_x(\varepsilon,0,0,c), \quad A_s(\varepsilon,c)=g_s(\varepsilon,0,0,c), 
\end{equation*}
where the functions $g_x$ and $g_y$ defined in Lemma~\ref{lem:vf-afterstr}.

From Lemma~\ref{lem:vf-afterstr} it follows that there exists functions $h_1(L)=O(L)$ and $h_2(\sigma)=O(\sigma)$, such that for all $\varepsilon \in [0,1]$ it holds that
\begin{eqnarray*}
  m_l(A_u) - h_1(L) - h_2(\sigma) \leq \inf_{c \in W^c} m_{l}(A_u(\varepsilon,c))  \leq \sup_{c \in W^c} \mu_{log}(A_u(\varepsilon,c)) \leq \mu_{log}(A_u) + h_1(L) + h_2(\sigma) \\
  m_l(A_s)  - h_1(L) - h_2(\sigma) \leq \inf_{c \in W^c} m_{l}(-A_s(\varepsilon,c))  \leq  \sup_{c \in W^c} \mu_{log}(-A_s(\varepsilon,c)) \leq \mu_{log}(-A_s) + h_1(L) + h_2(\sigma)
\end{eqnarray*}
where $m_l(A)$ and $\mu_{log}(A)$  are the logarithmic norms defined in Appendix~\ref{app:sec-lognorm} in Def.~\ref{def:lognorm}.

We assume that the  numbers $\mu_{\min}, \mu_{\max}, \lambda_{\min}, \lambda_{\max}$ are such that (we decrease $L$ and $\sigma$ if needed)
\begin{eqnarray}
 0 < \mu_{\min} \leq   m_l(A_u) - h_1(L) - h_2(\sigma), \qquad
\mu_{log}(A_u) + h_1(L) + h_2(\sigma)  \leq \mu_{\max} ,   \label{eq:mu-min}   \\
 0 < \lambda_{\min}\leq  m_l(A_s)  - h_1(L) - h_2(\sigma), \qquad  \mu_{log}(-A_s) + h_1(L) + h_2(\sigma) \leq  \lambda_{\max}. \label{eq:lambda-max}
\end{eqnarray}

Observe that the numbers $\mu_{\min}$, $\mu_{\max}$, $\lambda_{\min}$, $\lambda_{\max}$ can be made arbitrarily close to the numbers $\mu_{\min}(A)$, $\mu_{\max}(A)$, $\lambda_{\min}(A)$, $\lambda_{\max}(A)$  given in Def.~\ref{def:Alambmu} if we take $L$ and $\sigma$ small enough.

\begin{definition}
Let $N \subset \PhSpace$. For   $z \in N$ we define $T^\pm_N(z)$ by
\begin{eqnarray*}
T^+_N(z)=\sup\{T \in \mathbb{R}_+ \cup \{\infty\}: \quad \varphi(t,z)\in N, \quad t \in [0,T)\}, \\
T^-_N(z)=\inf\{T \in \mathbb{R}_- \cup \{-\infty\} : \quad \varphi(t,z)\in N, \quad t \in (T,0]\}.
\end{eqnarray*}
\end{definition}

\begin{lemma}
\label{lem:x-estm}
For any $\delta >0$ let $N=D(\delta)=\overline{B}_{n_u}(\delta) \times \overline{B}_{n_s}(\delta)\times W^c$ and
let  $\mu_1'$ be any constant such that $0<\mu_1' <\mu_{\mathrm{min}}$ (see  (\ref{eq:mu-min})).

 Then there exists a constant $C_x(\mu_1',N)$ such that for all $z \in N$ it holds that
\begin{equation}
 |\pi_x\varphi(t,z)| \leq C_x(\mu_1',N) |x| e^{\mu_1' t}, \quad t \in (T^-_N(z),0].\label{eq:x-estm}
\end{equation}
\end{lemma}
\textbf{Proof:}
From Lemma~\ref{lem:vf-afterstr} we have
\begin{equation}
  \dot{x}=g_x(\varepsilon,x,y,c)x=(A_u + O(L) + O(\sigma) + O(L)y)x. \label{eq:x'-straight}
\end{equation}
Let us fix $\delta$. If we take $0<\delta'\leq \delta$, such that $m_l(A_u + O(L) + O(\sigma) + O(L)y) > \mu_1'$, then we obtain for all $z \in D(\delta')$
  \begin{equation}
 |\pi_x\varphi(t,z)| \leq  |x| e^{\mu_1' t}, \quad t \in (T^-_{D(\delta')}(z),0].\label{eq:x-delta'-estm}
\end{equation}
Observe that the set $D(\delta')$ isolates $W^c$ (i.e., it is a maximal invariant set in $N$), hence all backward orbits starting in $N$ either converge to $W^c$
or leave it.  From Lemma~\ref{lem:hyp-behaviour}
 it follows that
 each backward orbit starting in $z \in N$ can be split into at most three parts: 1) for $t \in [0,T_1]$ $\varphi(-t,z) \notin D(\delta')$ and $\varphi(-T_1,z) \in D(\delta')$, 2) for $t \in [T_1,T_2)$ with $T_2=\infty$ for $z \in W^{cu}$ $\varphi(-t,z) \in D(\delta')$, 3) if $T_2 < \infty$, then there exist $T_3$ such that for $t \in (T_2,T_3]$
$\varphi(-t,z) \in N \setminus D(\delta')$ and $\varphi(-t-\eta,z) \notin N$ for $\eta>0$ arbitrary small.

Moreover, the time spend in the first and third parts are bounded from above by some constant $T_{13}$ and due to (\ref{eq:x'-straight}) we would have the following estimate
on first part  (if present)
\begin{equation*}
  |\pi_x \varphi(-t,z)| \leq  |x| e^{\mu_{log}(-g_x(N))t} \quad  t \in [0,T_1],
\end{equation*}
and analogously for the third part (if needed). This combined  with (\ref{eq:x-delta'-estm}) plus bounds on $T_1$ and $T_3-T_2$ gives us (\ref{eq:x-estm}) for a suitable constant $C_x(\mu_1',N)$.
\qed



\begin{lemma}
\label{lem:xy-estm}
Let us take any $\delta>0$.
 Then for any $\mu_1'$, $\lambda_1'$ such that $0<\mu_1' <\mu_{\mathrm{min}}$ and $0<\lambda_1' <\lambda_{\mathrm{min}}$, there exist constants
 $C_x(\mu_1')$ and $C_y(\lambda_1')$, such that for all $z \in D(\delta)$,
 \begin{eqnarray}
    |\pi_x\varphi(t,z)  \chi(\varphi(t,z))| \leq C_x(\mu_1') |x| e^{\mu_1' t}, \quad t \in (-\infty,0], \nonumber\\ 
     |\pi_y\varphi(t,z) \chi(\varphi(t,z))| \leq C_y(\lambda_1') |y| e^{-\lambda_1' t}, \quad t \in [0,\infty). \label{eq:y-estm-ext}
 \end{eqnarray}
\end{lemma}
\textbf{Proof:}
We use Lemma~\ref{lem:x-estm} with $N=D(\tilde{r}_1)$. The proof of (\ref{eq:y-estm-ext}) is analogous.
\qed

\subsection{Preparation of $R_1$ and $R_2$}
\label{subsec:prep-R1-R2}

We will not make any difference between the vector field case, where the remainder $\mathcal{R}(z)$ is a vector, and  the Hamiltonian case, where $R(z)$ is a real number. In both cases
we will use the same letter $R$. Below our $R$ is $T_{\varepsilon \ast}R$.
\begin{lemma}
\label{lem:prepR1R2-subtle}
Assume that  $R(\varepsilon,z) \in C^{q}$, and $Q$, $q$, $\ell_1,\ell_2$  are such that
\begin{eqnarray}
  q &\geq& Q+1,  \label{eq:q>=Q+1}\\
  R(\varepsilon,z)&=&O_{Q+1}(|(x,y)|^{Q+1}),\nonumber \\
  \ell_1 + \ell_2 &\leq& Q,  \label{eq:l1l2<=Q} \\
  2 \ell_2 + 1&\leq& Q. \label{eq:2l2<=Q}
\end{eqnarray}

Then there exist functions $R_1(\varepsilon, z)$ and $R_2(\varepsilon,z)$ such that
\begin{eqnarray}
  R_1,R_2 &\in& C^{q-\ell_2} ,\nonumber \\
  R_1(\varepsilon,z)&=&O_{\ell_1+1}(|x|^{\ell_1+1}), \quad R_2(\varepsilon,z)=O_{\ell_2+1}(|y|^{\ell_2+1}), \label{eq:R1R2Ol-bis} \\
  R(\varepsilon,z)&=&R_1(\varepsilon,z) + R_2(\varepsilon,z).\nonumber
\end{eqnarray}
Moreover, the supports of $R_1$ and $R_2$ are contained in the support of $R$.
\end{lemma}
\textbf{Proof:}
Below we will not write $\varepsilon$, $c$ in the arguments of $R$ and its derivatives.

For a given $x$ consider a Taylor formula with an integral remainder for the function $ y \mapsto R(x,\cdot)$
\begin{eqnarray}
  R(x,y) &=& \sum_{j=0}^{\ell_2}\frac{1}{j!}D_y^{j}R(x,0)(y)^{j}
    + \mathrm{Rem}(x,y)(y)^{\ell_2+1}, \label{eq:RTexp-y}
\end{eqnarray}
where  $\mathrm{Rem}(x,y)(y)^{\ell_2+1}$ is  given by
\begin{equation*}
   \mathrm{Rem}(z)= \left(\int_0^1 \frac{(1-t)^{\ell_2}}{\ell_2!} D_y^{\ell_2+1}R(x,ty) dt\right)(y)^{\ell_2+1} ,
\end{equation*}
where $D_y^{\ell_2+1}R(x,p)(y)^{\ell_2+1}$ is the $(\ell_2+1)$-linear map $D_y^{\ell_2+1}R(x,p)$ at $y=p$ applied to  the argument $(y,\dots,y)$ with $y$ appearing $\ell_2+1$ times.



We have for $s=0,\dots,\ell_2$
\begin{eqnarray*}
   D^s_y R(x,0)(y)^s \in C^{q-s}, \quad D^s_y R(x,0)(y)^s = O_{Q+1-s}\left(|x|^{Q+1-s}\right).
\end{eqnarray*}
We set
\begin{equation*}
  R_1(x,y)=R(x,0) + D_y R(x,0)y + \dots + \frac{1}{\ell_2!}D_y^{\ell_2}R(x,0)(y)^{\ell_2}.
\end{equation*}
Then
\begin{equation*}
 R_1 \in C^{q-\ell_2}, \qquad R_1(x,y)=O_{Q+1-\ell_2}(|x|^{Q+1-\ell_2}).
\end{equation*}


Since from (\ref{eq:l1l2<=Q}) we have $Q+1-\ell_2 \geq \ell_1+1$, then $R_1(x,y)=O_{\ell_1+1}(|x|^{\ell_1+1})$.

Since $R_2 = R - R_1$ we get
\begin{equation*}
  R_2(x,y)=\mathrm{Rem}(x,y) \in C^{q-\ell_2}.
\end{equation*}

It remains to show that  $R_2(x,y)=O_{\ell_2+1}(|y|^{\ell_2+1})$.
Observe that
\begin{equation*}
  R_2(x,y)= I(x,y)(y)^{\ell_2 +1},\qquad I(x,y) \in C^{q-\ell_2 -1}.
\end{equation*}
From (\ref{eq:q>=Q+1},\ref{eq:2l2<=Q}) we have that $q-\ell_2 - 1\geq \ell_2 + 1$, so that $I(x,y) \in C^{\ell_2 + 1}$. From this we obtain the second condition in (\ref{eq:R1R2Ol-bis}).
\qed



\section{Solving the cohomological equation}
\label{sec:estm-gen-coheq}

Before we begin let us comment that to solve the cohomological equation we use the method of characteristics.   The only difficulty is around $W^c$, but we ask for $\mathcal{R}$ (or $R$) to vanish on $W^c$, which makes it possible to solve it.   Once a solution in the neighborhood of $W^c$ is obtained  we can in a unique way
extend it to some neighborhoods of $W^{cs}$ and $W^{cu}$.

To deal with the cases of a general vector field and Hamiltonians in a uniform way, from now on we will not use different fonts for the vector fields. The nature of each object will be  clear from the context.

We  consider the equation (see (\ref{eq:coh-flows},\ref{eq:coh-ham}))
\begin{equation}
  D G (z) Z(z) =  M(z) G(z) + R(z), \label{eq:coh-gen}
\end{equation}
where $z \in \mathbb{R}^n$,  $G: \mathbb{R}^n \to \mathbb{R}^m$, $Z: \mathbb{R}^n \to \mathbb{R}^n$, $M(z) \in \mathrm{Lin}(\mathbb{R}^m,\mathbb{R}^m)$ and $R: \mathbb{R}^n \to \mathbb{R}^m$.
$Z$, $M$, $R$ are given functions, and we are looking for $G$.

In our case we have either  $n=m$ and $M(z)=DZ(z)$ for a general vector field or $m=1$ and $M(z)=0$ for the Hamiltonian case.

We assume that
\begin{equation*}
  Z \in C^{q(Z)}, \quad R \in C^{q(R)}. 
\end{equation*}

We apply the \emph{method of characteristics} with  characteristic lines given as
solutions of
\begin{equation}
  \dot{z}=Z(z).  \label{eq:dcoh-ode-n}
\end{equation}
Let $\varphi$ be a local flow induced by (\ref{eq:dcoh-ode-n}). It is well known that $\varphi(t,z)$ is in $C^{q(Z)}$.


Let $S(t,z) \in \mathrm{Lin}(\mathbb{R}^m,\mathbb{R}^m)$ be a solution of
\begin{equation*}
  \frac{d}{dt}S(t,z)= -S(t,z) M(\varphi(t,z)). 
\end{equation*}
with initial condition $S(0,z)=I$. In the Hamiltonian case $S(t,z)=1$.

Observe that
\begin{equation}
S(t,z)=D\varphi(t,z)^{-1},  \label{eq:S=V-1}
\end{equation}
where $D\varphi(t,z)$ satisfies the following equation
\begin{equation*}
  \frac{d}{dt}D\varphi(t,z)=  M(\varphi(t,z))D\varphi(t,z).  
\end{equation*}

Assume  that $G$  satisfies (\ref{eq:coh-gen}), then
\begin{eqnarray*}
  \frac{d}{dt} \left( S(t,z) G(\varphi(t,z)) \right)&=&
    - S(t,z) M(\varphi(t,z))G(\varphi(t,z))
    + S(t,z) DG(\varphi(t,z)) \frac{d}{dt} \varphi(t,z)\\
    &=& S(t,z)\left(DG(\varphi(t,z) Z(\varphi(t,z)) - M(\varphi(t,z)) G(\varphi(t,z)) \right)\\
    &=& S(t,z) \left(DG\cdot Z - M \cdot G \right)(\varphi(t,z))= S(t,z) R(\varphi(t,z)),
\end{eqnarray*}
hence we obtain
\begin{equation}
  \frac{d}{dt} \left( S(t,z) G(\varphi(t,z)) \right)=S(t,z) R(\varphi(t,z))  \label{eq:SG-diff-eq}
\end{equation}
and therefore
\begin{equation}
  S(t_2,z) G(\varphi(t_2,z))-  S(t_1,z) G(\varphi(t_1,z))=\int_{t_1}^{t_2}S(w,z)R(\varphi(w,z)) dw. \label{eq:SG-gen-integral}
\end{equation}
Observe that if $G$ is smooth and $G$ satisfies (\ref{eq:SG-diff-eq}) or (\ref{eq:SG-gen-integral}), then $G$ solves (\ref{eq:coh-gen}).

\subsection{Estimates on $\displaystyle \frac{\partial^r  S}{\partial z^r}$ and $\displaystyle \frac{\partial^r \varphi}{\partial z^r}$ }

The goal of this section is to derive bounds on  $\frac{\partial^r  S}{\partial z^r}$ and $\frac{\partial^r \varphi}{\partial z^r}$.

Note that from (\ref{eq:S=V-1}) it follows that (for the case of a  general vector field)
\begin{equation*}
  S(t,z)=\frac{\partial \varphi}{\partial z}(-t,z),
\end{equation*}
therefore  estimates for $S(t,z)$ and its derivatives can be obtained from those on $\dfrac{\partial \varphi}{\partial z}$.

\begin{lemma}
\label{lem:varZ-fsp-estm}


Let $N\subset \PhSpace $ be any convex set.
Assume that
\begin{equation}
  \mu_{log}(DZ(z)) \leq \alpha \neq 0, \quad z \in N.  \label{eq:log-norm-Z-full-sp}
\end{equation}

Then for any $z \in N$
and for any $k \in \mathbb{Z}_+$, such that  $1 \leq k \leq q(Z)$, it  holds that
\begin{itemize}
\item  if $\alpha >0$
\begin{equation*}
  \left|\frac{\partial^k \varphi}{\partial z^k}(t,z) \right| \leq C_k e^{k \alpha t}, \quad 0 \leq t \leq T^+_N(z), 
\end{equation*}
\item if $\alpha <0$
\begin{equation*}
  \left|\frac{\partial^k \varphi}{\partial z^k}(t,z)  \right| \leq C_k e^{ \alpha t}, \quad 0 \leq t \leq T^+_N(z) , 
\end{equation*}
\end{itemize}
 for some constants $C_k$, with $C_1=1$.
\end{lemma}
\textbf{Proof:}
Observe that  $D^j Z$ for $1 \leq j \leq q(Z)$ are bounded on $N$.

The variational equations of increasing orders are
\begin{eqnarray*}
  \frac{d}{dt} \frac{\partial \varphi}{\partial z}(t,z)&=& DZ(\varphi(t,z)) \frac{\partial \varphi}{\partial z}(t,z) \\
   \frac{d}{dt} \frac{\partial^2 \varphi}{\partial z^2}(t,z)&=& DZ(\varphi(t,z)) \frac{\partial^2 \varphi}{\partial z^2}(t,z) +
     D^2Z(\varphi(t,z)) \left(\frac{\partial \varphi}{\partial z}(t,z),\frac{\partial \varphi}{\partial z}(t,z) \right).
\end{eqnarray*}
For higher orders we have
\begin{equation*}
  \frac{d}{dt} \frac{\partial^k \varphi}{\partial z^k}(t,z)= DZ(\varphi(t,z)) \frac{\partial^k \varphi}{\partial z^k}(t,z) + R_k(t,z),
\end{equation*}
where $R_k$ is the sum of terms of the following form
\begin{equation*}
  R_{k,j}(t,z)=C D^{m}Z(\varphi(t,z))\left(\frac{\partial^{k_1} \varphi}{\partial z^{k_1}}(t,z),\frac{\partial^{k_2} \varphi}{\partial z^{k_2}}(t,z), \dots, \frac{\partial^{k_m} \varphi}{\partial z^{k_m}}(t,z)\right),  
\end{equation*}
where $C$ is a constant, $k_i >0$ and $\sum_{i=1}^m k_i=k$.

We will use formula (\ref{eq:ln-duhamel}) from Theorem~\ref{thm:lognorm-Duhamel} in the Appendix~\ref{app:sec-lognorm} applied to the variational equations of $k$-th order to estimate $ \left\| \frac{\partial^k \varphi}{\partial z^k}(t,z)\right\|$.

For $\dfrac{\partial \varphi}{\partial z}(t,y)$ we  obtain from (\ref{eq:log-norm-Z-full-sp}) and Theorem~\ref{thm:log-norm-derflow}
\begin{equation*}
  \left\| \frac{\partial \varphi}{\partial z}(t,z) \right\| \leq  e^{\alpha t}
\end{equation*}
and then inductively  (observe that $D^j Z$ for $1 \leq j \leq q(Z)$ are bounded on $N$ and the initial condition for $\dfrac{\partial^k \varphi}{\partial z^k}$ vanishes for $t=0$)
\begin{eqnarray*}
  \left\|\frac{\partial^k \varphi}{\partial z^k}(t,z) \right\| \leq \sum_j \int_0^t  e^{\alpha(t-s)} \left\| R_{k,j}(s,z) \right\|ds .
\end{eqnarray*}

For  $\alpha <0$,
from our induction assumption  it follows that each  term $R_{kj}$ will give a contribution bounded by
\begin{equation*}
\tilde{C}\int_0^t e^{\alpha(t-s)}\left(e^{\alpha s}\right)^k ds \leq \tilde{C}e^{\alpha t} \int_0^t e^{\alpha(k-1)s}ds \leq \frac{\tilde{C} e^{\alpha t}}{-\alpha(k-1)}.
\end{equation*}
For $\alpha >0$ this contribution will be bounded by
\begin{equation*}
\tilde{C}\int_0^t e^{\alpha(t-s)}\left(e^{k \alpha s}\right)\leq \tilde{C}e^{\alpha t} \int_0^t e^{\alpha(k-1)s}ds < \frac{\tilde{C} e^{k \alpha t}}{\alpha (k-1)}.
\end{equation*}
\qed
\begin{rem}
\label{rem:varZ-fsp-estm-r1}
If the vector field $\mathcal{Z}$ depends smoothly  on parameters, the above bounds also applies to partial derivatives involving parameters at the price of taking arbitrary $\alpha'> \alpha$, $\alpha'>0$ and increasing constants $C_k$.
\end{rem}
\textbf{Proof:}
 We add equation $\dot{\varepsilon}=0$  to the system. In a suitable norm on the extended phase space  $\PhSpace \times \Lambda$, where $\Lambda$ is some compact set in parameter space we will have $\alpha=\sup_\varepsilon \sup_{z \in N} \mu_{log}(DZ_\varepsilon(z))$ to be approximately equal to $\displaystyle \alpha'=\sup_\varepsilon \sup_{z \in N} \mu_{log}\left(\frac{\partial Z_\varepsilon(z)}{\partial (z,\varepsilon)}(z,\varepsilon)\right)$.
\qed

\begin{rem}
\label{rem:varZ-fsp-estm-r2}
If $N=\overline{B}_{n_u}(\delta) \times \overline{B}_{n_s}(\delta)\times W^c$  is a neighborhood of $W^c$, then  in the above lemma  instead of $\alpha=\sup_{z \in N}\mu_{log}(DZ(z))$ we can take
any $\alpha' > \sup_{z \in W^c}\mu_{log}(DZ(z))$ at the price  of increasing the constants $C_k$.
\end{rem}
\textbf{Proof:}
The argument is the same as in the proof of Lemma~\ref{lem:xy-estm}.
\qed

\subsection{Solving  for particular remainder---backward in time }
\label{subsec:solve-fsp}
We solve (\ref{eq:coh-gen}) with $R_1(z)=O_{\ell+1}(|x|^{\ell+1})$.
We set (in the case of a general vector field)
\begin{equation*}
  M(z)=DZ(z), \qquad M(z) \in C^{q(Z)-1}, 
\end{equation*}
and in the Hamiltonian case
\begin{equation*}
 M(z)=0. 
\end{equation*}
We are in the setting discussed at the beginning of this section and  we use (\ref{eq:SG-gen-integral}) with $t_2=0$ and $t_1 \to -\infty$. Let us assume that
\begin{equation}
  S(t_1,z) G(\varphi(t_1,z)) \to 0, \quad t_1 \to -\infty,  \label{eq:SG-good-limit}
\end{equation}
which has to be verified later, and then we obtain
\begin{equation}
  G(z)=\int_{-\infty}^0 S(t,z)R(\varphi(t,z))dt. \label{eq:Gsolx}
\end{equation}

\begin{lemma}
\label{lem:sol-full-space}
Assume that $Z=Z_\varepsilon$ and $R=R_\varepsilon=R_1$ are as obtained after straightening of invariant manifolds and the decomposition of $T_\ast R$  described in Lemma~\ref{lem:prepR1R2-subtle} .
Let $\mu_{\min}, \mu_{\max}, \lambda_{\min}, \lambda_{\max}$ satisfy (\ref{eq:mu-min}--\ref{eq:lambda-max}).

Assume that for some positive integer $\ell$ it holds that
\begin{eqnarray}
q(R) &\geq& \ell+1,  \label{eq:qRl1} \\
 R(z)&=& O_{\ell +1}(|x|^{\ell+1}),   \label{eq:Rls-estm} \\
  Z(z) &\in& C^{q(Z)}.\nonumber
\end{eqnarray}

Assume that for the general vector field
\begin{equation}
  \ell + 1 \leq q(Z)-1, \label{eq:elqZ}
\end{equation}
and in the Hamiltonian case
\begin{equation}
 \ell + 1 \leq q(Z).  \label{eq:elqZ-2}
\end{equation}

Assume that  $k \geq 1$  satisfies
\begin{itemize}
\item for the general vector field
\begin{equation}
  k < \min \left(\frac{\mu_{\mathrm{min}} (\ell+1) - \mu_{\mathrm{max}} }{\mu_{\mathrm{min}} + \lambda_{\mathrm{max}}}, \frac{\mu_{\mathrm{min}}(\ell +1) - \mu_{\mathrm{max}}}{\mu_{\mathrm{max}}}\right), \label{eq:cond2}
\end{equation}
\item for the Hamiltonian vector field
\begin{equation}
  k < \frac{\mu_{\mathrm{min}} (\ell+1)}{\mu_{\mathrm{min}} + \lambda_{\mathrm{max}}}. \label{eq:cond2-ham}
\end{equation}
\end{itemize}

Let us define $G_\varepsilon$ by (\ref{eq:Gsolx}) for all $(\varepsilon,z) \in [0,1] \times \PhSpace$. Then $G(\varepsilon,z) \in C^k$ and $G_\varepsilon$ is a solution of (\ref{eq:coh-gen}).

Moreover,
\begin{eqnarray}
  G(\varepsilon,z)&=&O_{k}(|x|^{\ell+1}),\quad \mbox{for} \ |x| \to 0,  \label{eq:GOk} \\
  \|G(\varepsilon,z)\| &\leq& K |x|^{\ell+1}, \quad \forall z, \label{eq:G-large-z} \\
  \|DG(\varepsilon,z)\| &\leq& K |x|^{\ell}, \quad \forall z. \label{eq:DG-large-z}
\end{eqnarray}
\end{lemma}
\noindent
\textbf{Proof:}
In this proof we consider $\varepsilon$ as a part of the phase space, hence it enters as one of the components in the variable $z$ and with the vector field $Z$, $G$ ,$R$ extended so that it vanishes in the $\varepsilon$ direction. On this extended space the cohomological equation is still a cohomological equation with respect to variables $(\varepsilon,z)$.

Let  $\xi=1$ for a general vector field and $\xi=0$ in the Hamiltonian case.

Let $\alpha \approx \mu_{\mathrm{min}}$, $\gamma \approx \lambda_{\mathrm{max}}$ be such that
\begin{equation*}
0<\alpha < \mu_{\mathrm{min}}, \quad \gamma>\lambda_{\mathrm{max}}, \quad \beta > \mu_{\mathrm{max}}.
\end{equation*}
Let   $\gamma$ be a bound on the  logarithmic norm for the backward in time evolution close to $W^c$, which should be close
 to $\mu_{log}(-A_s(c))=\lambda_{\max} \approx \mu_{log}(-A_s)=\lambda_{\mathrm{max}}(A)$, $\gamma > \lambda_{\max}$ and
  let $\beta$ be a logarithmic norm for forward in time evolution
 close to $W^c$, which should be $> \mu_{\max}$ but close.

For any $z $, from Lemmas~\ref{lem:xy-estm} and~\ref{lem:varZ-fsp-estm} (see also Remarks~\ref{rem:varZ-fsp-estm-r1} and \ref{rem:varZ-fsp-estm-r2} )  we have the following bounds for $t \in (-\infty,0]$:
\begin{eqnarray}
    \|\pi_x\varphi(t,z)\| \chi(\varphi(t,z)) &\leq& C_\alpha |x| e^{\alpha t}, \label{eq:x-phi-estm} \\
    \left\| \frac{\partial^j \varphi}{\partial z^j}(t,z) \right\|  &\leq& C_{\gamma,j} e^{-j \gamma t}, \quad j=1,\dots,q(Z),  \label{eq:D^jvarphi-bnd} \\
    \left\| \frac{\partial^j S}{\partial z^j}(t,z)  \right\|  &\leq& C'_{\beta,j} e^{-(j+1) \beta t}, \quad j=0,\dots,q(Z)-1.  \label{eq:D^jS-bnd}
\end{eqnarray}
where $ \chi$ is the indicator function for a ball containing support of $R$  (see Def.~\ref{def:chi-tilr1} and Lemma~\ref{lem:xy-estm}).

Observe that from  (\ref{eq:cond2},\ref{eq:cond2-ham}) and assumption $k \geq 1$  it follows that
\begin{equation}
 \xi \mu_{\mathrm{max}} + \lambda_{\mathrm{max}} <  \ell \mu_{\mathrm{min}}. \label{eq:cond1}
\end{equation}

From (\ref{eq:cond1}) for $\beta \to \mu_{\mathrm{max}}$, $\gamma \to \lambda_{\mathrm{max}}$ and $\alpha \to \mu_{\mathrm{min}}$ we have
\begin{equation}
 \xi \beta + \gamma <  \ell \alpha, \label{eq:bet-gam-alp-l}
\end{equation}
From (\ref{eq:cond2},\ref{eq:cond2-ham}) we should also have
\begin{equation}
  k < \frac{\alpha (\ell+1) - \xi \beta }{\alpha + \gamma}, \label{eq:chsolfp-reg}
\end{equation}
and for the case of a general vector field additionally the following holds
\begin{equation}
   k < \frac{\alpha (\ell+1) -  \beta }{\beta}. \label{eq:k-sec-ineq}
\end{equation}

Since $R$ has support in $\{|y| \leq \tilde{r}_1\}$, then the integral in (\ref{eq:Gsolx}) is as smooth as $R$ and $S$ for $z=(x,y,c)$ with $y \neq 0$, as  each backward trajectory spends only a finite interval of time in the set $\{|y| \leq \tilde{r}_1\}$.

Since
\begin{equation*}
   G(z)=\lim_{T \to -\infty} \int_{T}^0 S(t,z)R(\varphi(t,z))dt, 
\end{equation*}
we will show that $\int_{T}^0 S(t,z)R(\varphi(t,z))dt$ converges uniformly in $C^k$-norm as $T \to -\infty$.

From (\ref{eq:Rls-estm}) and since $R$ has support contained in $D(\tilde{r}_1)$ it follows that
\begin{equation}
 \| D^s R(\varepsilon,z) \| \leq K_s \|x\|^{\ell + 1- s}, \quad s=0,\dots,\ell+1, \quad \forall (\varepsilon,z). \label{eq:DsR-estm}
\end{equation}

From (\ref{eq:x-phi-estm}) and (\ref{eq:DsR-estm}) it follows that  for constants $K_{\alpha,s}=C_\alpha^{\ell+1-s} K_s$
and for any $z$ it holds that
\begin{equation}
\left\lvert (D^sR) \left(\varphi(t,z)\right)\right\rvert\leq K_{\alpha,s}\abs{x}^{\ell+1-s}\mathrm{e}^{(\ell+1-s)\alpha t} \text{ for } t\leq 0, \quad s=0,\dots,\ell+1.  \label{eq:Rbnd-ls}
\end{equation}

From the above expression with $s=0$ and from (\ref{eq:D^jS-bnd}) with $j=0$ (recall that for the Hamiltonian case we have $S=I$) we have for $t \leq 0$ and some constant $C$
\begin{eqnarray*}
  \lvert S(t,z)R(\varphi(t,z)) \rvert \leq  C \mathrm{e}^{((\ell+1)\alpha-\xi \beta) t} \abs{x}^{\ell+1}.  
\end{eqnarray*}
Since from (\ref{eq:bet-gam-alp-l}) it follows that
\begin{equation}
   \alpha (\ell+1) - \xi \beta >0, \label{eq:ls-beta-alpha}
\end{equation}
we see that the improper integral defining $G(z)$ in eq.~(\ref{eq:Gsolx}) is convergent.
Moreover, we have the following bound for some constant $K$, valid for all $z$
\begin{equation}
 |G(z)|\leq K |x|^{\ell +1}. \label{eq:G-bnd-x}
\end{equation}

To complete the proof that $G$ is a solution of  cohomological equation (\ref{eq:coh-gen}), we need to show that $G \in C^1$ and that condition (\ref{eq:SG-good-limit}) holds.

First we deal with (\ref{eq:SG-good-limit}). 
We have  from (\ref{eq:D^jS-bnd},\ref{eq:x-phi-estm},\ref{eq:G-bnd-x}),  for $t < 0$
\begin{eqnarray*}
  |S(t,z)G(\varphi(t,z))| &\leq& K |S(t,z)| \cdot |\pi_x \varphi(t,z)|^{\ell +1}\\
   &\leq&  K \mathrm{e}^{-\xi \beta t} \abs{x}^{\ell+1}\mathrm{e}^{(\ell+1)\alpha t}
  \leq K  \abs{x}^{\ell+1} \mathrm{e}^{-(\xi \beta -(\ell+1)\alpha) t},
\end{eqnarray*}
hence it converges to $0$ for $t \to -\infty$ due to (\ref{eq:ls-beta-alpha}).

The derivatives  of $G$  (of order $k \leq \ell+1$)  will be given by expressions of the following form
\begin{eqnarray}
  D^k G(z)=\int_{-\infty}^0 \sum_{j=0}^k C_j D_z^j S(t,z) \cdot D_z^{k-j}\left(R(\varphi(t,z))\right) dt,  \label{eq:DkG-fsp}
\end{eqnarray}
where for the Hamiltonian case only the term with $j=0$ is present and $S=I$.

Bounds for $\left\| D_z^j S(t,z) \right\|$ are given by (\ref{eq:D^jS-bnd}). We need to find estimates for $D_z^{k-j}\left(R(\varphi(t,z))\right)$.

It is easy to see that  $D_z^{j}\left(R(\varphi(t,z))\right)$ is the sum of terms of the following form, with $s=1,\dots,j$
\begin{equation*}
  (D^{s}R)(\varphi(t,z))\left(\frac{\partial^{k_1} \varphi}{\partial z^{k_1}}(t,z),\frac{\partial^{k_2} \varphi}{\partial z^{k_2}}(t,z), \dots, \frac{\partial^{k_m} \varphi}{\partial z^{k_m}}(t,z)\right),
\end{equation*}
where  $k_i >0$ and $\sum_{i=1}^s k_i=j$.

From the above using (\ref{eq:Rbnd-ls},\ref{eq:D^jvarphi-bnd})  it is easy to  see that for  $t \leq 0$ we have, for some constant $C$,
  \begin{eqnarray}
  \left\| D_z^j (R(\varphi(t,z)))  \right\| \leq C |x|^{\ell+1 - j} e^{\alpha(\ell+1 -j)t} e^{-j \gamma t}, \quad j=0,1,\dots,\ell+1.  \label{eq:DjzRphi}
\end{eqnarray}
By combining (\ref{eq:DkG-fsp},\ref{eq:D^jS-bnd},\ref{eq:DjzRphi}) we obtain the following bound for the case of a general vector field for $k \leq \ell+1$
\begin{eqnarray}
  \|D^k G(z)\| \leq C \int_{-\infty}^0  \sum_{j=0}^k e^{-\beta(j+1)t} |x|^{\ell+1 - (k-j)} e^{\alpha(\ell+1 -(k-j))t} e^{-(k-j) \gamma t} dt. \label{eq:DkG-integral}
\end{eqnarray}

 In the Hamiltonian case   only the term with $j=0$ is present in (\ref{eq:DkG-integral}) with $\beta=0$, hence   for the convergence of the above improper integral we need that
 \begin{equation*}
    \alpha(\ell+1)  -k (\alpha + \gamma) >0,
 \end{equation*}
 which is implied by (\ref{eq:chsolfp-reg}).

In the case of a general vector field   we need that for $j=0,\dots,k$ it holds that
 \begin{equation*}
  \alpha(\ell+1 -(k-j))-\beta(j+1)- \gamma (k-j) =
  \alpha(\ell+1) -k(\alpha + \gamma) -\beta - j(\beta - (\alpha + \gamma)) >0.
\end{equation*}
It is easy to see that this is enough to satisfy the above inequality for $j=0$ and $j=k$.  We obtain
\begin{eqnarray}
  \alpha(\ell+1) -k(\alpha + \gamma) -\beta >0,
  \label{eq:back-1cond} \\
   \alpha(\ell+1)  -(k+1)\beta  >0.  \label{eq:back-2cond}
\end{eqnarray}
Observe that (\ref{eq:back-1cond}) is implied by (\ref{eq:chsolfp-reg}) and condition (\ref{eq:back-2cond}) by (\ref{eq:k-sec-ineq}).
\qed

\subsection{Solving  for a particular remainder---forward in time }
In this section we discuss solving equation (\ref{eq:coh-gen})
but this time we assume that the remainder term
 $R$ decays with $y$.

In (\ref{eq:SG-gen-integral}) we pass to the limit $t_2 \to \infty$ and set $t_1=0$ to get
\begin{equation}
  G(z)=-\int_{0}^{\infty}S(w,z)R(\varphi(w,z)) dw, \label{eq:Gsol}
\end{equation}
provided the above integral is convergent  and
\begin{equation*}
   S(t_2,z) G(\varphi(t_2,z)) \to 0 , \quad t_2 \to \infty.  
\end{equation*}

\begin{lemma}
\label{lem:sol-full-space-frw}
Assume that $Z$ and $R=R_2$ are as obtained after straightening of invariant manifolds.

Let $\mu_{\min}, \mu_{\max}, \lambda_{\min}, \lambda_{\max}$ satisfy (\ref{eq:mu-min}--\ref{eq:lambda-max}).
Assume that
\begin{eqnarray}
q(R) &\geq& \ell+1,  \label{eq:qRl1-frw} \\
    R(z)&=&O_{\ell+1}(|y|^{\ell+1}), \nonumber \\ 
  Z(z) &\in& C^{q(Z)}.\nonumber
\end{eqnarray}

Assume that
\begin{equation}
  \ell + 1 \leq q(Z)-1 , \label{eq:el1qZ-frw}
\end{equation}
and in the Hamiltonian case
\begin{equation}
 \ell + 1 \leq q(Z). \label{eq:el1qZfrw-2}
\end{equation}

Assume that  $k \geq 1$  satisfies
\begin{itemize}
\item for a general vector field
\begin{equation*}
  k < \min \left(\frac{\lambda_{\mathrm{min}} (\ell+1) - \lambda_{\mathrm{max}} }{\lambda_{\mathrm{min}} + \mu_{\mathrm{max}}}, \frac{\lambda_{\mathrm{min}}(\ell +1) - \lambda_{\mathrm{max}}}{\lambda_{\mathrm{max}}}\right), 
\end{equation*}
\item for a Hamiltonian vector field
\begin{equation*}
  k < \frac{\lambda_{\mathrm{min}} (\ell+1)}{\lambda_{\mathrm{min}} + \mu_{\mathrm{max}}}. 
\end{equation*}
\end{itemize}

Let us define $G_\varepsilon(z)$ by (\ref{eq:Gsol}) for all $(\varepsilon,z) \in [0,1] \times \PhSpace$. Then $G(\varepsilon, z) \in C^k$ and $G$ is a solution of (\ref{eq:coh-gen}).

Moreover,
\begin{eqnarray}
 G(\varepsilon,z)&=&O_{k}(|y|^{\ell+1}),\quad \mbox{for} \ |x| \to 0,  \label{eq:GfOk} \\
  \|G(\varepsilon,z)\| &\leq& K |y|^{\ell+1}, \quad \forall z, \label{eq:Gf-large-z} \\
  \|DG(\varepsilon,z)\| &\leq& K |y|^{\ell}, \quad \forall z. \label{eq:DGf-large-z}.
\end{eqnarray}

\end{lemma}
The proof is the same as for Lemma~\ref{lem:sol-full-space} simply reversing the direction of time, and exchanging $ x \leftrightarrow y$, $\lambda \leftrightarrow \mu$.

\subsection{Conclusion of the proof of Theorem~\ref{thm:polNForm-loccenter} }

\subsubsection{Inequalities related to the regularity}

We need   to show that the  inequalities relating $k$ (the regularity class of $G$, which is the solution of the cohomological equation) to $q$ (the regularity of the vector field, $\mathcal{Z}, \mathcal{R} \in C^q$) and to $Q$, $\ell_1$, $\ell_2$, can be satisfied. These inequalities appear in Lemmas~\ref{lem:vf-afterstr}, \ref{lem:prepR1R2-subtle}, \ref{lem:sol-full-space} and~\ref{lem:sol-full-space-frw}.

Let us notice that in the Hamiltonian case the solution of cohomological equation gives us a Hamiltonian defining the vector field for the deformation method, hence $G \in C^{k+1}$. For the same reason we have $R, Z \in C^{q+1}$ and  $R(z,c)=O_{Q+2}(|z|^{Q+2})$, (see (\ref{eq:Ham-prepR2})), hence when applying Lemmas~\ref{lem:sol-full-space} and~\ref{lem:sol-full-space-frw} in the Hamiltonian setting we should replace there $k$ by $k+1$ and $Q$ by $Q+1$.

From Lemma~\ref{lem:vf-afterstr} we obtain that $T_{\ast} \mathcal{Z},T_{\ast} \mathcal{R} \in C^{q-2}$, and in the Hamiltonian case, $T_{\ast}R \in C^{q-1}$.  In the pure saddle case, because the stable and unstable manifolds are as smooth as the vector field, we have $T_{\ast} \mathcal{Z},T_{\ast} \mathcal{R} \in C^{q-1}$, and $T_\ast R \in C^q$ in the Hamiltonian case.

For Lemma~\ref{lem:prepR1R2-subtle} applied to the vector field after straightening of invariant manifolds
to obtain $R_1,R_2$ in $C^{q - 2 -\ell_2}$ and in $C^{q - 1 -\ell_2}$ in Hamiltonian case, and $C^{q-1-\ell_2}$ in pure saddle case ($C^{q-\ell_2}$ in Hamiltonian case),
 we need the following inequalities for the general vector field case (in the pure saddle case  the first inequality will be $q-1 \geq Q+1$ )
\begin{eqnarray}
  q-2 &\geq& Q+1,  \label{eq:q>=Q+1-bis} \\
  \ell_1 + \ell_2 &\leq& Q,  \label{eq:l1l2<=Q-bis} \\
  2 \ell_2 + 1 &\leq& Q, \label{eq:2l2<=Q-bis}
\end{eqnarray}
while  in the Hamiltonian case we have the following inequalities  (in the pure saddle case  the first inequality will be $q \geq Q+2$)
\begin{eqnarray}
 q-1 &\geq& Q+2,  \label{eq:Hq>=Q+1-bis} \\
  \ell_1 + \ell_2 &\leq& Q+1,  \label{eq:Hl1l2<=Q-bis} \\
  2 \ell_2 + 1 &\leq& Q+1. \label{eq:H2l2<=Q-bis}
\end{eqnarray}

 Observe that from (\ref{eq:q>=Q+1-bis},\ref{eq:l1l2<=Q-bis}) it follows that conditions (\ref{eq:qRl1},\ref{eq:elqZ},\ref{eq:elqZ-2})
 and (\ref{eq:qRl1-frw},\ref{eq:el1qZ-frw},\ref{eq:el1qZfrw-2}) from Lemmas ~\ref{lem:sol-full-space} and ~\ref{lem:sol-full-space-frw}, respectively, are satisfied.

 From Lemmas~\ref{lem:sol-full-space} and~\ref{lem:sol-full-space-frw} we obtain the following  inequalities for the general vector field
\begin{eqnarray}
  k &<& \min \left(\frac{\mu_{\mathrm{min}} (\ell_1+1) - \mu_{\mathrm{max}} }{\mu_{\mathrm{min}} + \lambda_{\mathrm{max}}}, \frac{\mu_{\mathrm{min}}(\ell_1 +1) - \mu_{\mathrm{max}}}{\mu_{\mathrm{max}}}\right), \label{eq:cond2-bis} \\
  k &<&  \min \left(\frac{\lambda_{\mathrm{min}} (\ell_2+1) - \lambda_{\mathrm{max}} }{\lambda_{\mathrm{min}} + \mu_{\mathrm{max}}}, \frac{\lambda_{\mathrm{min}}(\ell_2 +1) - \lambda_{\mathrm{max}}}{\lambda_{\mathrm{max}}}\right), \label{eq:frw-cond2-bis}
\end{eqnarray}
and for the Hamiltonian case
\begin{eqnarray}
  k+1 &<& \frac{\mu_{\mathrm{min}} (\ell_1+1)}{\mu_{\mathrm{min}} + \lambda_{\mathrm{max}}}, \label{eq:cond2-ham-bis} \\
   k+1 &<& \frac{\lambda_{\mathrm{min}} (\ell_2+1)}{\lambda_{\mathrm{min}} + \mu_{\mathrm{max}}}. \label{eq:cond2-ham-frw-bis}
\end{eqnarray}
 The order of setting parameters is as follows. We pick any $k \geq 1$. Then we find the minimal $\ell_1$ satisfying (\ref{eq:cond2-bis}),
 \begin{equation*}
   \ell_1(k) =  \max \left(k+ \left[ k \frac{\lambda_{\mathrm{max}}}{\mu_{\mathrm{min}}} +  \frac{\mu_{\mathrm{max}}}{\mu_{\mathrm{min}}}   \right],\left[(k+1)\frac{\mu_{\max}}{\mu_{\mathrm{min}}} \right]  \right) , 
 \end{equation*}
 and in the Hamiltonian case we use  (\ref{eq:cond2-ham-bis})
 \begin{equation*}
   \ell_1^H(k)= k+1 + \left[(k+1) \frac{\lambda_{\mathrm{max}}}{\mu_{\mathrm{min}}} \right]. 
 \end{equation*}

  Analogously we find the minimal $\ell_2$ satisfying (\ref{eq:frw-cond2-bis}).
 \begin{equation*}
   \ell_2(k) = \max\left(k +  \left[ k \frac{\mu_{\mathrm{max}}}{\lambda_{\mathrm{min}}} + \frac{\lambda_{\mathrm{max}}}{\lambda_{\mathrm{min}}}   \right], \left[(k+1)\frac{\lambda_{\mathrm{max}}}{\lambda_{\mathrm{min}}} \right]\right) 
 \end{equation*}
 and for the Hamiltonian case from (\ref{eq:cond2-ham-frw-bis})
 \begin{equation*}
   \ell_2^H(k)= k+1 + \left[(k+1) \frac{\mu_{\mathrm{max}}}{\lambda_{\mathrm{min}}} \right].
 \end{equation*}

From (\ref{eq:l1l2<=Q-bis}) for general vector field we obtain
 \begin{eqnarray*}
   Q\geq Q_0 &=& \ell_1(k) + \ell_2(k) \\
   &=&  \max \left(k+ \left[ k \frac{\lambda_{\mathrm{max}}}{\mu_{\mathrm{min}}} +  \frac{\mu_{\mathrm{max}}}{\mu_{\mathrm{min}}}   \right],\left[(k+1)\frac{\mu_{\max}}{\mu_{\mathrm{min}}} \right]  \right) + \\
   & & \max\left(k +  \left[ k \frac{\mu_{\mathrm{max}}}{\lambda_{\mathrm{min}}} + \frac{\lambda_{\mathrm{max}}}{\lambda_{\mathrm{min}}}   \right], \left[(k+1)\frac{\lambda_{\mathrm{max}}}{\lambda_{\mathrm{min}}} \right]\right)
 \end{eqnarray*}
 For the Hamiltonian case  from (\ref{eq:Hl1l2<=Q-bis}) we obtain
 \begin{eqnarray*}
   Q \geq Q_0 = \ell^H_1(k) + \ell^H_2(k) - 1= 2k + 1 + \left[(k+1) \frac{\lambda_{\mathrm{max}}}{\mu_{\mathrm{min}}} \right] + \left[(k+1) \frac{\mu_{\mathrm{max}}}{\lambda_{\mathrm{min}}} \right] \\
   =  2k + 1 + 2\left[(k+1) \frac{\mu_{\mathrm{max}}}{\mu_{\mathrm{min}}} \right].
 \end{eqnarray*}

 Observe that we ignored  condition (\ref{eq:2l2<=Q-bis}) (and (\ref{eq:H2l2<=Q-bis})), but  we can always rearrange the decomposition $R=R_1+R_2$, so that $\ell_2 \leq \ell_1$
 and then from (\ref{eq:q>=Q+1-bis}) (and (\ref{eq:Hq>=Q+1-bis})) it follows that $q_0=Q_0+3$ and
 in the pure saddle case $q_0=Q_0+2$.

 This finishes the proof of Remark~\ref{rem:Q0q0}.

\subsubsection{Transformation removing the remainder}
Consider the general vector field case.
From the solution of cohomological equation for $\varepsilon \in [0,1]$ we obtain the vector field $\mathcal{G}(\varepsilon,z)=T_{\varepsilon \ast}^{-1}(G_1(\varepsilon,z)+G_2(\varepsilon,z))$, which is at least $C^k$ in $(\varepsilon,z)$.

$\mathcal{G}(\varepsilon,z)$ is defined for $\varepsilon \in [0,1]$ and arbitrary $z$ and satisfies for all $\varepsilon \in [0,1]$ the following estimates,
see  (\ref{eq:GOk},\ref{eq:G-large-z},\ref{eq:DG-large-z}) in Lemma~\ref{lem:sol-full-space}  and (\ref{eq:GfOk},\ref{eq:Gf-large-z},\ref{eq:DGf-large-z}) in Lemma~\ref{lem:sol-full-space-frw}
\begin{eqnarray}
  \mathcal{G}(\varepsilon,z)&=&O_k(|(x,y)|^{\ell +1}), \quad (x,y) \to 0,  \label{eq:G-bnd-close-0} \\
  \|\mathcal{G}(\varepsilon,z)\| &\leq& K \|(x,y)\|^{\ell +1}, \quad \forall z, \label{eq:G-bnd-big-z} \\
  \|D\mathcal{G}(\varepsilon,z)\| &\leq& K \|(x,y)\|^{\ell}, \quad \forall z, \label{eq:DG-bnd-big-z}
\end{eqnarray}
where $k \leq \ell+1$, and  $\ell=\min(\ell_1(k),\ell_2(k))$ with $\ell_i(k)$ defined in the previous section.

Note that $T^{-1}_{\varepsilon \ast}$ preserves conditions (\ref{eq:G-bnd-close-0},\ref{eq:G-bnd-big-z},\ref{eq:DG-bnd-big-z}), which were proved
for $G_1+G_2$.

The differential equation defining $g(\varepsilon,z)$ is
\begin{equation}
  \frac{d}{d\varepsilon} g(\varepsilon,z)=\mathcal{G}(\varepsilon,g(\varepsilon,z)). \label{eq:ode-def-g}
\end{equation}

From (\ref{eq:G-bnd-close-0}) it follows that the points from the center manifold are fixed points of $g_1$ and $Dg_1=Id$ for such points.

In view of (\ref{eq:G-bnd-big-z}) the solution may be defined only for a finite time (i.e., $\varepsilon$) interval.  To estimate the length of this
interval we proceed as follows. As
\begin{eqnarray*}
 \frac{d \|\pi_{xy}g(\varepsilon,z)\|}{d \varepsilon} \leq K \|\pi_{xy}g(\varepsilon,z)\|^{\ell+1},
\end{eqnarray*}
we compare $\|\pi_{xy}g(\varepsilon,z)\|$  to the solution of the equation
\begin{equation*}
  v'=K v^{\ell+1}, \quad v(0)=\|\pi_{xy} z\|,
\end{equation*}
which is given by (we assume that $v(0)>0$)
\begin{equation*}
  v(t)=\frac{v(0)}{(1 - K\ell t v(0)^{\ell})^{1/\ell}},
\end{equation*}
and we obtain
\begin{eqnarray*}
  \|\pi_{x,y}g(\varepsilon,z) \|  \leq \frac{\|z\|}{(1 - K\ell t \|\pi_{xy} z\|^{\ell})^{1/\ell}}.
\end{eqnarray*}
Therefore if
\begin{equation*}
  \|\pi_{xy} z\| \leq \frac{1}{(K \ell)^{1/\ell}}, 
\end{equation*}
then $g(1,z)$ is defined.

The same reasoning applies to the inverse map of $g_1$, which is a solution of (\ref{eq:ode-def-g}) moving backward in time.



This finishes the proof of Theorem~\ref{thm:polNForm-loccenter}.

\section{Sign symmetry}
\label{sec:invWpm}

Assume that (compare (\ref{eq:PhSpace})) $\PhSpace=\{(z,c), z=(z_+,z_-) \in \mathbb{R}^{n_1} \times \mathbb{R}^{n_2}, c \in \mathbb{C}^{2m}\}$, i.e.,
we identify $\mathbb{R}^{k_A} \times \mathbb{C}^{l_A}$ with   $\mathbb{R}^{n_1} \times \mathbb{R}^{n_2}$.

For $l=1,\dots,m$ we define the hyperplane $W_l$ as
\begin{equation*}
W_l=\{ (z_+,z_-,c) \in \PhSpace  :  c_l=0\},
\end{equation*}
and $W_\pm$ as
\begin{eqnarray*}
  W_+=\{ (z_+,z_-,c) \in \PhSpace :  z_+=0\}, \\
  W_-=\{ (z_+,z_-,c) \in \PhSpace :  z_-=0\}.
\end{eqnarray*}

We also define the following subspaces
\begin{eqnarray*}
V_j &=& \{ (z_+,z_-,c) \in \PhSpace: \, z_\pm=0, c_l =0,\ l\neq j \}=\bigcap_{l\neq j} W_l  \cap W_- \cap W_+, \\
V_+ &=& \{ (z_+,z_-,c) \in \PhSpace: \, z_-=0, c_l =0,\ 1\leq l \leq m \}=\bigcap_{l=1}^m W_l  \cap W_-, \\
V_-&=& \{ (z_+,z_-,c) \in \PhSpace: \, z_+=0, c_l =0,\ 1\leq l \leq m \}=\bigcap_{l=1}^m W_l  \cap W_+.
\end{eqnarray*}

Observe that our phase space has the following decomposition
\begin{equation*}
 \PhSpace  = V_+ \oplus V_- \oplus V_1 \oplus \cdots \oplus V_m. 
\end{equation*}

Let $X$ and $Y$ be the unstable and stable manifolds of the linearization at $0$. We assume that we have decompositions $X=X_+ \oplus X_-$ and $Y=Y_+ \oplus Y_-$.
We also assume that
\begin{equation*}
  V_+=X_+ \oplus Y_+, \quad V_-=X_- \oplus Y_-.
\end{equation*}
We introduce the variables   $x_{\pm}$, $y_\pm$ by  $z_+=(x_+,y_+) \in X_+ \oplus Y_+$ and $z_-=(x_-,y_-) \in X_- \oplus Y_-$.  We will use coordinates $x_{i\pm}$ and $y_{i\pm}$
in $X_\pm$ and $Y_\pm$, so that $z_-$ is given by $x_{i-}$ and $y_{i-}$, etc.  We will use multiindices $\alpha \in \mathbb{N}^{k_1}$ and $\beta \in \mathbb{N}^{k_2}$, to define $z_+^\alpha z_-^\beta$ and
   $\dfrac{\partial^{|\alpha| R}}{\partial z_+^\alpha}$, etc.

\begin{definition}

For any choice of signs $s_+,s_-,s_j \in \{-1,1\}$ consider the map $\SSym_s$ be given by
\begin{equation*}
  \SSym_s(z_+,z_-,c_1,\dots,c_m)=( s_+z_+,s_-z_-,s_1 c_1,\dots,s_m c_m).
\end{equation*}
We say that a differential equation or map or symplectic form has a \emph{sign-symmetry}  if it is $\SSym_s$-symmetric for any choice of signs $\{s=(s_+,s_-,s_1,\dots,s_m)\}$.

\end{definition}

 In the Hamiltonian case we assume that
\begin{equation}
  \Omega = \sum_j dx_{j-} \wedge dy_{j-} + \sum_j dx_{j+} \wedge dy_{j+} + \frac{i}{2} \sum_k dc_k \wedge d\bar{c}_k.  \label{eq:OmegaSym}
\end{equation}
Observe that if $\Omega$ has sign-symmetry, $\Omega(\SSym_s(z,c))=\Omega(z,c)$ for any $s$ and $(z,c) \in \PhSpace$.

Observe that the sign-symmetry implies the invariance of the spaces $W_l$ for a vector field or map and also the subspaces $V_l$ and $V_\pm$ are invariant. Moreover, the origin $0$ is a fixed point of an o.d.e. or map
under consideration.

It is easy to see that
\begin{itemize}
\item the inverse map of $G$  with  sign-symmetry has also this property
\item the composition of two maps with sign-symmetry has also this property,
\item a vector field with  sign-symmetry  induces a flow for which any time shift has also these properties.
\end{itemize}

Our goal is to establish the following extension of Theorem~\ref{thm:polNForm-loccenter}.

\begin{theorem}
\label{thm:polNForm-centerSubspaces}
Consider  system~\eqref{eq:fpcen} with the same assumptions as in Theorem~\ref{thm:polNForm-loccenter}.
Then the change of coordinates bringing it to polynomial normal form~\eqref{eq:cen}
has  sign-symmetry if system~\eqref{eq:fpcen} has sign-symmetry and is symplectic if the transformed system was Hamiltonian.
\end{theorem}
For the  proof of Theorem~\ref{thm:polNForm-centerSubspaces} we show that each step in the construction of coordinate change in Theorem~\ref{thm:polNForm-loccenter} is performed so that the sign-symmetry are preserved.
Obviously the regularity of the coordinate change is the same as obtained in Theorem~\ref{thm:polNForm-loccenter}.
Corollary~\ref{cor:polNForm-loccenterAreal} is also valid in the case of sign-symmetry.

 \subsection{The toy model system for transfer of energy in NLSE}
\label{subsec:Q0-spec-system}

The motivation of the present work was the toy model system derived in \cite{CK} in the context of the transfer of energy in cubic de\-fo\-cusing non-linear
Schr\"o\-din\-ger equation (NLSE), see also \cite{DZ} and references given there.

The phase space is described by $z=(x_-,y_-,x_+,y_+,c) \in \mathbb{R}^4 \times \mathbb{C}^N$. The toy model system is
\begin{eqnarray}
\dot x_-&=&\hphantom{-}\lambda x_-+O_2(|z|^2),\notag\\
\dot y_-&=&-\lambda y_-+O_2(|z|^2),\notag\\
\dot x_+&=&\hphantom{-}\lambda x_+ +O_2(|z|^2),\label{eq:4Dsys}\\
\dot y_+&=&-\lambda y_+ +O_2(|z|^2),\notag \\
\dot c_l &=&\hphantom{-} i \nu_k c_l + O(|z)|^2), \quad l=1,\dots,N  \notag
\end{eqnarray}
where $\lambda>0$, $\nu_l >0$.
We assume that $\{(x_-,y_-,x_+,y_+)=0\}$ is invariant and the system has the sign-symmetry.

The following result follows from Corollary~\ref{cor:polNForm-loccenterAreal} adapted to include the sign-symmetry.
\begin{lemma}
\label{lem:PolyNormForm}
For any $k\geq 1$ if system~\eqref{eq:4Dsys} is $C^q$ with $q$ sufficiently
 large, then there exists a $C^k$ change of variables in a neighborhood of the origin,
 transforming system~\eqref{eq:4Dsys} to the system
 \begin{equation*}
\begin{split}
  \dot{x}_-&= \hphantom{-}\lambda x_- +  N_{x_-}(z), \\
  \dot{y}_-&= -\lambda y_- +  N_{y_-}(z),  \\
  \dot{x}_+&= \hphantom{-} \lambda x_+ + N_{x_+}(z), \\ 
  \dot{y}_+&= -\lambda y_+ + N_{y_+}(z), \\
 \dot{c}_\ell &= \hphantom{-}  i\nu_l c_\ell +O_2(z), \quad \ell=1,\dots,N 
\end{split}
\end{equation*}
where for any saddle variable $v \in \{x_-,y_-,x_+,y_+\}$ we have
\[
  N_v(z)=\sum_{m \in M_{1,v}}  g_{v,m}(c_*)z^m  + \sum_{m \in M_{2,v}} g_{v,m}(z)z^m,  
\]
where   $g_{v,m}$ are continuous functions,
$M_{1,v},M_{2,v}$ are \emph{finite} sets of indices,
and any $z^m$ is a resonant monomial for the saddle variables,
satisfying
on the one hand $m_s:=m_{x_-}+m_{y_- }+m_{x_+ } +m_{y_+}\geq 3$  and $m_c =0$ if $m=\left(m_{x_-},m_{y_- },m_{x_+},m_{y_+},m_c\right) \in M_{ 1, v }$, and on the other hand $m_s=1$ and $m_c \geq 3$ if $m\in M_{2,v }$.



This change of coordinates preserves  the sign-symmetry and if the original system is Hamiltonian, then this change is symplectic.
\end{lemma}
\textbf{Proof:}
We take any $P=Q+2$, this gives us $P+1-Q=3$ which gives $m_c \geq 3$ in the second sum.   We have  $m_s\geq 3$ because there are no resonant terms of order 2.
\qed

Let us see what $q$ we need to guarantee $k=2$. For this we  use Remark~\ref{rem:Q0q0}.  We have (see Def.~\ref{def:Alambmu})
\begin{equation*}
  \mu_{\min}=\mu_{\max}=\lambda_{\min}=\lambda_{\max}=\lambda.
\end{equation*}
Hence from (\ref{eq:Q0-our}) for a general vector field we obtain
\begin{equation*}
  Q_0(k)=4k+2, \quad q_0(k)=4k+5, 
\end{equation*}
while for a Hamiltonian system (see (\ref{eq:Q0-ham-our})) we obtain
\begin{equation*}
  Q_0(k)=4k+3, \quad q_0(k)=4k+6. 
\end{equation*}
For a general vector field, we need $q \geq q_0$, $q-1\geq P \geq Q$, $q\geq Q+3$ (see the statement of Theorem~\ref{thm:polNForm-loccenter}). We take $P=Q_0(k)+2$ to obtain
\begin{equation*}
  P=4k+4, \quad q\geq 4k+5.
\end{equation*}
Hence to obtain $k=2$ we need $q \geq 13$. In the Hamiltonian case we need $q \geq 14$.

\subsection{Conditions of sign-symmetry for Hamiltonian systems}

For Hamiltonians with the symplectic form (\ref{eq:OmegaSym}) the sign-symmetry is
\begin{equation}
  R(z)=R(\mathcal{S}_s(z))  \label{eq:ham-signsym}
\end{equation}
for any choice of signs  $\{s\}$.
It is easy to see that the Hamiltonian vector field induced by $R$ using the symplectic form (\ref{eq:OmegaSym}) also has sign-symmetry.

However after straightening  center-stable and center-unstable manifolds the symplectic form $T_\ast \Omega$ is no longer in standard form (\ref{eq:OmegaSym}) and the relation between Hamiltonian and its induced vector field  given by (\ref{eq:def-ham-vf}) is   more  complicated than (\ref{eq:ham-eq}). However, we will show that the straightening
transformation $T_\varepsilon$ has also sign-symmetry, as well as $T_{\varepsilon \ast} \Omega$ and the condition on the sign-symmetry of $T_{\ast} R$ will be still
(\ref{eq:ham-signsym}).

\section{Preservation of sign symmetry for all steps leading to polynomial normal forms}

\subsection{The sign-symmetry for the Taylor expansion in saddle directions}


The goal of this section is to show that if a map or vector field or Hamiltonian function has  the sign-symmetry, then each term in the Taylor expansion has the sign-symmetry. This is a
content of Lemma~\ref{lem:eachTermSSym}. This result is needed to justify that removing the non-resonant terms can be done preserving the sign-symmetry.

Assume $G \in C^{Q+1}$ (vector field or map).
The Taylor formula with respect to saddle directions $(z_+,z_-)$ with an integral remainder might be written as follows
\begin{eqnarray}
  G_{\pm,l}(z_+,z_-,c) &=& \sum_{i=0}^{Q} \sum_{|\alpha|+|\beta|=i} a^{\alpha,\beta} \frac{\partial^{i} G_{\pm,l}}{\partial z_+^\alpha \partial z_-^\beta}(z_\pm=0,c) z_+^\alpha  z_-^\beta \notag \\
  &+&
  \sum_{|\alpha|+|\beta|=Q+1} d^{\alpha,\beta}_{\pm,l} (z_+,z_-,c)z_+^\alpha z_-^\beta,  \label{eq:G-Taylor}
\end{eqnarray}
where $a^{\alpha,\beta}$ are some constants  and $d^{\alpha,\beta}_{\pm,l}(z_+,z_-,c)$ are coefficients of the $(Q+1)$-linear map representing  the remainder
\begin{equation*}
 G_{\pm l,\mathrm{Rem}}(z_+,z_-,c)=  \left(\int_0^1 \frac{(1-t)^{Q}}{Q!} D_{(z_+,z_-)}^{(Q+1)}G_{\pm,l}(tz_+,tz_-,c) dt\right)(z_+,z_-)^{Q+1},
\end{equation*}
i.e., up to the constant depending on $\alpha$, $\beta$ and $Q$,
\begin{equation*}
   d^{\alpha,\beta}_{\pm,l}(z_+,z_-,c) = a^{\alpha,\beta} \int_0^1 (1-t)^{Q} \frac{\partial^{Q+1}G_{\pm,l}}{\partial z_+^\alpha \partial z_-^\beta}(tz_+,tz_-,c) dt.
\end{equation*}

For each term (\ref{eq:G-Taylor}) we can define a map as follows. For example the term $a^{\alpha,\beta} \frac{\partial^{i} G_+}{\partial z_+^\alpha \partial z_-^\beta}(0,0,c) z_+^\alpha  z_-^\beta$
induces the map $T(z_+,z_-,c)=(T_+(z_+,z_-,c),T_-(z_+,z_-,c),T_1(z_+,z_-,c),\dots,T_m(z_+,z_-,c))$, with
\begin{equation*}
  T_+(z_+,z_-,c)=a^{\alpha,\beta} \frac{\partial^{i} G_+}{\partial z_+^\alpha \partial z_-^\beta}(0,0,c) z_+^\alpha  z_-^\beta, \quad T_- \equiv 0, \quad T_l\equiv 0, 
\end{equation*}
and analogously for the other terms.

\begin{lemma}
\label{lem:eachTermSSym}
Assume that $G \in C^{Q+1}$ has  the sign-symmetry. Then each non-vanishing term (i.e. not identically equal to $0$) in the Taylor expansion (\ref{eq:G-Taylor}) also has the sign-symmetry.
\end{lemma}
\textbf{Proof:}
The sign-symmetry of all the terms with respect to the changes of sign on $c_k$ is immediate.   Hence it is enough to consider the changes of signs of $z_+$ (for  $z_-$ the proof is analogous).

Consider $G_l$ (for $G_-$ the argument is the same). The sign-symmetry with respect $z_+$ means that $G_l(z_+,z_-,c)=G_l(-z_+,z_-,c)$. This implies that
\begin{equation*}
  \frac{\partial^{|\alpha|+|\beta|} G_l}{\partial z_+^\alpha \partial z_-^\beta}(z_+=0,z_-,c)=0, \quad \mbox{if $|\alpha|$ is odd}.
\end{equation*}
Therefore all terms of order less than or equal to $Q$ have the symmetry $z_+ \to -z_+$.

Now we look at the remainder. Since for all $t$ it holds that
\begin{equation*}
  \frac{\partial^{|\alpha|+|\beta|} G_l}{\partial z_+^\alpha \partial z_-^\beta}(tz_+,tz_-,c)=(-1)^{|\alpha|}\frac{\partial^{|\alpha|+|\beta|} G_l}{\partial z_+^\alpha \partial z_-^\beta}(t(-z_+),tz_-,c),
\end{equation*}
we see that
\begin{equation*}
 \frac{\partial^{|\alpha|+|\beta|} G_l}{\partial z_+^\alpha \partial z_-^\beta}(tz_+,tz_-,c)z_+^{\alpha}=(-1)^{|\alpha|}\frac{\partial^{|\alpha|+|\beta|} G_l}{\partial z_+^\alpha \partial z_-^\beta}(t(-z_+),tz_-,c)(-z_+)^{\alpha},
\end{equation*}
and therefore
\begin{equation*}
d_l^{\alpha,\beta}(-z_+,z_-,c)=d_l^{\alpha,\beta}(z_+,z_-,c).
\end{equation*}
This establishes the sign-symmetry with respect to $z_+$ of all the remainder terms for $G_l$.
It remains to consider $G_+$.  The sign-symmetry with respect $z_+$ means that
$G_+(z_+,z_-,c)=-G_+(z_+,z_-,c)$. This implies that
\begin{equation*}
  \frac{\partial^{|\alpha|+|\beta|} G_+}{\partial z_+^\alpha \partial z_-^\beta}(z_+=0,z_-,c)=0, \quad \mbox{if $|\alpha|$ is even.}
\end{equation*}
Therefore all terms of order less than or equal than $Q$ have the symmetry $z_+ \to -z_+$.

Now we look at the remainder. Since
\begin{equation*}
  \frac{\partial^{|\alpha|+|\beta|} G_+}{\partial z_+^\alpha \partial z_-^\beta}(tz_+,tz_-,c)=(-1)^{|\alpha|+1}\frac{\partial^{|\alpha|+|\beta|} G_+}{\partial z_+^\alpha \partial z_-^\beta}(t(-z_+),tz_-,c),
\end{equation*}
we see that
\begin{equation*}
 \frac{\partial^{|\alpha|+|\beta|} G_+}{\partial z_+^\alpha \partial z_-^\beta}(tz_+,tz_-,c)z_+^{\alpha}=(-1)^{|\alpha|+1}\frac{\partial^{|\alpha|+|\beta|} G_+}{\partial z_+^\alpha \partial z_-^\beta}(t(-z_+),tz_-,c)(-z_+)^{\alpha},
\end{equation*}
and therefore
\begin{equation*}
d_+^{\alpha,\beta}(-z_+,z_-,c)=-d_+^{\alpha,\beta}(z_+,z_-,c).
\end{equation*}
This establishes the sign-symmetry with respect to $z_+$ of all the remainder terms for $G_+$.

\qed

\subsubsection{Sign-symmetry for Taylor expansions for Hamiltonians }

The Taylor formula with an integral remainder might be written as follows
\begin{eqnarray}
  G(z_+,z_-,c) &=& \sum_{i=0}^{Q} \sum_{|\alpha|+|\beta|=i} a^{\alpha,\beta} \frac{\partial^{i} G}{\partial z_+^\alpha \partial z_-^\beta}(0,0,c) z_+^\alpha  z_-^\beta +  \sum_{|\alpha|+|\beta|=Q+1} d^{\alpha,\beta} (z_+,z_-,c)z_+^\alpha z_-^\beta  \label{eq:G-ham-Taylor}
\end{eqnarray}
where $a^{\alpha,\beta}$ are some constants and $d_{\alpha,\beta}(z_+,z_-,c)$ are the coefficients of the $(Q+1)$-linear form of the remainder
\begin{equation*}
 G_{\mathrm{Rem}}(z_+,z_-,c)=  \left(\int_0^1 \frac{(1-t)^{Q}}{Q!} D_{(z_+,z_-)}^{(Q+1)}G(tz_+,tz_-,c) dt\right)(z_+,z_-)^{Q+1},
\end{equation*}
i.e., up to the constant depending on $\alpha$ and $\beta$,
\begin{equation*}
   d^{\alpha,\beta}(z_+,z_-,c) = a^{\alpha,\beta} \int_0^1 (1-t)^{Q} \frac{\partial^{Q+1}G}{\partial z_+^\alpha \partial z_-^\beta}(tz_+,tz_-,c) dt.
\end{equation*}

\begin{lemma}
\label{lem:HamAllTermSsym}
For $Q \geq 2$ assume that $G \in C^{Q+2}$  has the sign-symmetry. Then every non-vanishing term (i.e. not identically equal to $0$) in the Taylor expansion (\ref{eq:G-ham-Taylor}) also has the sign-symmetry
\end{lemma}
We omit an easy proof similar to the proof of Lemma~\ref{lem:eachTermSSym}.    

\subsection{The sign-symmetry  for the terms in the Taylor expansion in  $y$-direction}


The goal of this section is to show that all the terms defined in the decomposition obtained in Lemma~\ref{lem:prepR1R2-subtle} have  the sign-symmetry if the original vector field (or map) has it.

In this section we use coordinates $(x,y,c)$, $x=(x_+,x_-)$, $y=(y_+,y_-)$, so that $z_+=(x_+,y_+)$ and $z_-=(x_-,y_-)$.

Given a vector field (or map) $R$   we write it as (compare with (\ref{eq:RTexp-y}) )
\begin{eqnarray}
    R(x,y,c) &=& R(x,0,c) + D_y R(x,0,c)y + \dots + \frac{1}{\ell!}D_y^{\ell}R(x,0,c)(y)^{\ell} \notag \\
    & & + \mathrm{Rem}(x,y,c)(y)^{\ell+1} \label{eq:R-dcmp-Ty}
\end{eqnarray}

\begin{lemma}
\label{lem:R-subT-Vpm-SSym}
  Assume that $R \in C^{q(R)}$ has the sign-symmetry and $\ell+1 \leq q(R)$. Then every term in (\ref{eq:R-dcmp-Ty}) also has  the sign-symmetry.
\end{lemma}
We omit an easy proof.

\subsubsection{The sign-symmetry for the terms in the Taylor expansion
in the  $y$-direction for Hamiltonians}

Recall that the condition for the sign-symmetry for a Hamiltonian $G$ is given by (see (\ref{eq:ham-signsym}))
\begin{equation*}
  G(z)=G(\mathcal{S}_s(z))  
\end{equation*}
for any choice of signs  $\{s\}$.

Given a Hamiltonian $R$   we write it as (compare with (\ref{eq:RTexp-y}) )
\begin{eqnarray}
    R(x,y,c) &=& R(x,0,c) + D_y R(x,0,c)y + \dots + \frac{1}{\ell!}D_y^{\ell}R(x,0,c)(y)^{\ell} \notag \\
    & & + \mathrm{Rem}(x,y,c)(y)^{\ell+1}. \label{eq:H-dcmp-Ty}
\end{eqnarray}

\begin{lemma}
\label{lem:ham-R-subl-Ssym}
Assume that $R \in C^{q(R)}$  and $\ell+1 \leq q(R)$. If $R$ has the sign-symmetry, then every non-vanishing term (i.e. not identically equal to $0$) in the Taylor expansion (\ref{eq:H-dcmp-Ty}) also has the sign-symmetry.
\end{lemma}
We omit an easy  proof.

\subsection{The sign-symmetry is preserved in the preparation of the remainder}

Assume that $R$ has the sign-symmetry.

The decomposition $R=R_1 + R_2$ is performed using Lemma~\ref{lem:prepR1R2-subtle}, then from Lemma~\ref{lem:R-subT-Vpm-SSym} for general vector fields  and Lemma~\ref{lem:ham-R-subl-Ssym} for Hamiltonians we know that
 $R_1$ and $R_2$ have the sign-symmetry.

\subsection{The transformation removing non-resonant
terms preserves the the sign-symmetry}

Assume that vector field $\mathcal{Z}$ (or the Hamiltonian $Z$) has  the sign-symmetry.

For a general vector field, the transformation bringing $\mathcal{Z}$ to the normal form is a composition of maps
of the type $\mathrm{Id}$ plus some non-resonant terms (see \cite{A}). From Lemma~\ref{lem:eachTermSSym} it follows that each such term has the desired geometric properties, and the induced transformation also has.

For Hamiltonians the transformation removing non-resonant terms (see for example \cite{G,GDF}) is a composition of symplectic maps, which are the time shift $1$ of the Hamiltonian flow induced by the polynomial non-resonant term in the Hamiltonian to be  removed. From Lemma~\ref{lem:HamAllTermSsym} it follows that each such term  has the sign-symmetry. So the induced time-$1$ also has sign-symmetry.

For reference purposes, we formulate the above as a lemma.
\begin{lemma}
  The transformation of removing non-resonant terms can be performed in such a way that the sign-symmetry is preserved. Moreover, in the Hamiltonian case it is symplectic.
\end{lemma}

\subsection{The sign-symmetry  and preparation of  compact data}

When preparing ``compact data" in Section~\ref{subsec:prep-compt-data} we modified the vector fields and Hamiltonians.

The question is: if $\mathcal{Z}$ and $\mathcal{R}$ have the sign-symmetry, then $\tilde{\mathcal{Z}}$ and $\tilde{\mathcal{R}}$ introduced in(\ref{eq:dotx-lcd-app2}--\ref{eq:dotc-lcd-app2}) will also have the sign-symmetry? Analogously for the Hamiltonian $\tilde H$ introduced in~\eqref{eq:mod-ham}.
The answer is yes if the function $\eta$ used in both cases is constructed as follows.

Let $\alpha: \mathbb{R} \to [0,1]$ be $C^\infty$ with compact support and such that $\alpha(x)=\alpha(-x)$, so that $\alpha'(0)=0$. We define $\eta$, which will not depend on $c$, by
\begin{equation}
\eta(x_{1\pm},x_{2\pm},\dots,y_{1\pm},y_{2\pm},\dots)=\Pi_i \alpha(x_{i+})  \cdot \Pi_i \alpha(x_{i-})\Pi_i \alpha(y_{i+})  \cdot \Pi_i \alpha(y_{i-}). \label{eq:eta-wpm}
\end{equation}
Observe that
\begin{equation*}
\frac{\partial \eta}{\partial x_{i\pm}}(p)=0, \quad \mbox{if $x_{i\pm}(p)=0$},\qquad \frac{\partial \eta}{\partial y_{i\pm}}(p)=0, \quad \mbox{if $y_{i\pm}(p)=0$}. 
\end{equation*}

The following lemma is obvious
\begin{lemma}
  Let $\eta$ be as in (\ref{eq:eta-wpm}).
Assume that the vector field $\mathcal{Z}$ and $\mathcal{R}$ have the sign-symmetry. Then $\tilde{\mathcal{Z}}$ and $\tilde{\mathcal{R}}$ have it too. The same for the Hamiltonian $\tilde{Z}$ and $\tilde{R}$, as long as $Z$ and $R$ have the sign-symmetry.
\end{lemma}

\subsection{The transformation straightening the  invariant manifolds  preserves  the sign-symmetry}

Let $(x_-,x_+,y_-^u(x_-,x_+,c),y_+^u(x_-,x_+,c),c)$ be the  center-unstable manifold \newline and $(x_-^s(y_-,y_+,c),\allowbreak x_+^s(y_-,y_+,c),\allowbreak y_-,y_+,c)$ be the  center-stable  manifold as obtained in Lemma~\ref{lem:iso-block-mod-vekf}.

We have the following result about some  properties of the functions $y_\pm^u$ and $x_\pm^s$.
\begin{lemma}
If the vector field (\ref{eq:Z-with-rem}) has the sign-symmetry, then for any $s_+,s_- \in \{-1,1\}$ it holds that
\begin{eqnarray}
    y^u(x,c_k,c_*)&=&y^u(x,-c_k,c_*), \label{eq:yu-ssym-ck} \\
    x^s(y,c_k,c_*)&=&x^s(y,-c_k,c_*), \label{eq:xs-ssym-ck} \\
    y_+^u(s_+ x_+,s_- x_-,c)&=&s_+ y_+^u(x_+,x_-,c), \label{eq:yu+ssym} \\
    y_-^u(s_+x_+,s_-x_-,c)&=&s_-y_-^u(x_+,x_-,c), \label{eq:yu-ssym} \\
    x_+^s(s_+ y_+,s_- y_-,c)&=&s_+ x_+^s(y_+,y_-,c)\nonumber \\
     x_-^s(s_+y_+,s_-y_-,c)&=&s_-x_-^s(y_+,y_-,c). \nonumber 
\end{eqnarray}
\end{lemma}
\textbf{Proof:}
It is enough to consider the equalities involving the function $y^u$, because analogous symmetric arguments give us the statements for $x^s$.

The center-unstable manifold is the graph of the function  $(x_-,x_+,c) \mapsto  (y_-^u(x_-,x_+,c), y_+^u(x_-,x_+,c))$.
The sign-symmetry ($c_k \to -c_k$) implies that $(x,y^u(x,c_k,c_*),c_k,c_*) \in W^{cu}$ iff $(x,y^u(x,c_k,c_*),-c_k,c_*) \in W^{cu}$. Therefore
\begin{equation*}
  y^u(x,c_k,c_*)=y^u(x,-c_k,c_*).
\end{equation*}
This establishes (\ref{eq:yu-ssym-ck}) (and by symmetry (\ref{eq:xs-ssym-ck})).

The sign-symmetry ($z_- \mapsto -z_-$) implies that
if $(x_+,x_-,y_+^u(x_+,x_-,c),y_-^u(x_+,x_-,c),c)\in W^{cu}$,
 then  $(x_+,-x_-,y_+^u(x_+,x_-,c),-y_-^u(x_+,x_-,c),c)\in W^{cu}$. Therefore
\begin{eqnarray*}
  y_+^u(x_+,-x_-,c)=y_+^u(x_+,x_-,c),\quad y_-^u(x_+,-x_-,c)=-y_-^u(x_+,x_-,c).
\end{eqnarray*}
This argument and an analogous one for the sign-symmetry $z_+ \mapsto -z_+$ establishes (\ref{eq:yu+ssym}) and (\ref{eq:yu-ssym}).
\qed

From the above Lemma we immediately obtain the following result.
\begin{lemma}
If the  vector field (\ref{eq:Z-with-rem})  has the sign-symmetry then, for every $\varepsilon \in [0,1]$, the
transformation $T_\varepsilon$  defined by(\ref{eq:S-man-tran})   has the sign-symmetry.
\end{lemma}

\subsection{The solutions of the cohomological equation have the sign-symmetry }

Recalling that
\begin{itemize}
\item $\dfrac{\partial \varphi}{\partial z}(t,z)$ and $S(t,z)$ were defined in Section~\ref{sec:estm-gen-coheq},
\item $G_1$ and $G_2$ were defined by \eqref{eq:Gsolx} and \eqref{eq:Gsol}, respectively,
\end{itemize}
it is obvious that $G_1$ and $G_2$ have the sign-symmetry.

\subsection{Conclusion of the proof of Theorem~\ref{thm:polNForm-centerSubspaces}}
We have shown that all stages in the construction of the coordinate change of Theorem~\ref{thm:polNForm-loccenter}  have  the sign-symmetry. This finishes the proof of Theorem~\ref{thm:polNForm-centerSubspaces}.


\appendix

\section{Logarithmic norms and related topics}
\label{app:sec-lognorm}

In this section we state some facts about logarithmic norms and their applications to ODEs \cite{D,HNW,KZ,L}.

We define for a  linear map $A: V \to V$ of a normed space
\begin{equation}
  m(A)=\inf_{\|x\|=1} \|A(x)\|.  \label{eq:minA}
\end{equation}
For an interval matrix $\mathbf{A}\subset \mathbb{R}^{k\times n}$ we set
\begin{equation}
m(\mathbf{A})=\inf_{A\in \mathbf{A}}m(A). \label{eq:minA-intv}
\end{equation}

\begin{definition}
\label{def:lognorm}
For a square matrix $A\in \mathbb{R}^{n\times n}$ we define
the logarithmic norm of $A$ denoted by $\mu_{\rm{log}}(A)$ by \cite{D,HNW,KZ,L}
\begin{equation*}
\mu_{\rm{log}}(A)=\lim_{h\rightarrow 0^{+}}\frac{\Vert I+hA\Vert -\Vert I\Vert }{h}
\end{equation*}%
 and the logarithmic minimum of $A$ \cite{CZ17}
\begin{equation*}
m_{l}(A)=\lim_{h\rightarrow 0^{+}}\frac{m(I+hA)-\Vert I\Vert }{h}.
\end{equation*}
\end{definition}

It is known that $\mu_{\rm{log}}(A)$ is a continuous and  convex function. Moreover, see \cite[Lem. 3]{CZ17},
\begin{equation*}
  m_l(A)=-\mu_{\rm{log}}(-A).  
\end{equation*}

In the following theorem a bound on the distance between solutions of an ODE is established in
terms of the logarithmic norm. The proof of this result can be found in \cite{HNW} (the part involving $\mu_{\rm{log}}$) and in \cite[Thm. 5]{CZ17} for the lower bound
involving $m_l$.

\begin{theorem}
Consider an ODE
\begin{equation}
x^{\prime }=f(t,x),  \label{eq:non-auto-ode}
\end{equation}%
where $x\in \mathbb{R}^{n}$ and $f:\mathbb{R}\times \mathbb{R}%
^{n}\rightarrow \mathbb{R}^{n}$ is $C^{1}$.

Let $x(t)$ and $y(t)$ for $t\in [t_{0},t_{0}+T]$ be two solutions of (\ref{eq:non-auto-ode}). Let $W\subset \mathbb{R}^{n}$ such that for each $t\in [t_{0},t_{0}+T]$ the segment connecting $x(t)$ and $y(t)$ is
contained in $W$. Let
\begin{eqnarray*}
L&=&\sup_{x\in W,t\in [ t_{0},t_{0}+T]}\mu_{\rm{log}}\left( \frac{\partial f}{\partial x}(t,x)\right), \\
l&=&\inf_{x \in W, t \in [t_0,t_0+T]} m_l\left( \frac{\partial f}{\partial x}(t,x)\right).
\end{eqnarray*}
Then for $t\in [0,T]$ it holds that
\begin{equation*}
\exp (lt)\left\Vert x(t_{0})-y(t_{0})\right\Vert  \leq \left\Vert x(t_{0}+t)-y(t_{0}+t)\right\Vert \leq \exp (Lt)\left\Vert x(t_{0})-y(t_{0})\right\Vert .
\end{equation*}
\end{theorem}

From the above result one easy derives the following theorem.
\begin{theorem}
\label{thm:log-norm-derflow}Consider an ODE
\begin{equation}
x^{\prime }=f(x),  \label{eq:auto-ode}
\end{equation}%
where $x\in \mathbb{R}^{n}$ and $f: \mathbb{R}^{n}\rightarrow \mathbb{R}^{n}$ is $C^{1}$.

Let $x(t)$ be a solution of (\ref{eq:auto-ode}) and let $W\subset \mathbb{R}^{n}$ such that $x(t) \in W$ for each $
t\in [0,T]$. Let
\begin{eqnarray*}
L&=&\sup_{x\in W} \mu_{\rm{log}}\left(\frac{\partial f}{\partial x}(t,x)\right), \\
l&=&\inf_{x \in W} m_l\left( \frac{\partial f}{\partial x}(t,x)\right)
\end{eqnarray*}
Then for $t\in [0,T]$ it holds that
\begin{equation*}
  e^{lt} \leq m\left(\frac{\partial \varphi}{\partial x}(t,x) \right), \qquad  \left\|\frac{\partial \varphi}{\partial x}(t,x) \right\| \leq \exp (Lt), \quad t \in [0,T]
\end{equation*}
\end{theorem}

The following theorem  follows  from Lemma 4.1 in \cite{KZ}
\begin{theorem}
\label{thm:comptestm}
 Assume that either
\begin{itemize}
\item  $z:[0,T] \to {\mathbb
R}^k$ is a solution of the equation
\begin{equation}
  z'(t)=A(t) \cdot z(t) + \delta(t),  \label{eq:z'Atz}
\end{equation}
where $\delta:[0,T] \to {\mathbb R}^k$ and $A:[0,T]\to {\mathbb
R}^{k \times k}$ are continuous.
\end{itemize}
or
\begin{itemize}
 \item $z:[0,T] \to \mathbb{R}^{k\times k}$ is a solution of the equation
\begin{equation}
  z'(t)=A(t) \cdot z(t) + \delta(t),  \label{eq:z'Atz-mat}
\end{equation}
or the equation
\begin{equation}
  z'(t)= z(t) \cdot A(t) + \delta(t), \label{eq:z'Atz-mat-tran}
\end{equation}
where $\delta:[0,T] \to {\mathbb R}^{k\times k}$ and $A:[0,T]\to {\mathbb
R}^{k \times k}$ are continuous.
\end{itemize}

Assume that the continuous  functions  $J:[0,T] \to \mathbb{R}$ and $C:[0,T] \to \mathbb{R}_+$ satisfy the following inequalities for
all $t \in [0,T]$
\begin{equation*}
\mu_{\rm{log}}(A(t) \leq J(t), \quad |\delta(t)| \leq C(t).
\end{equation*}

Then
\begin{equation*}
  |z(t)| \leq y(t)
\end{equation*}
where  $y:[0,T] \to \mathbb{R}^n$ is a solution of the
problem
\begin{equation*}
   y'(t)=J(t)y(t) + C(t), \qquad y(0)=|z(0)|.
\end{equation*}
\end{theorem}

While the result for (\ref{eq:z'Atz}) is indeed a direct consequence of Lemma 4.1 in \cite{KZ}, the statements for (\ref{eq:z'Atz-mat}) and (\ref{eq:z'Atz-mat-tran}) are obtained by the same reasoning which led to Lemma 4.1 in \cite{KZ}.

From Theorem~\ref{thm:comptestm} we immediately obtain the following result.
\begin{theorem}
\label{thm:lognorm-Duhamel}
  Under the same assumptions about $z(t)$, $A(t)$, $\delta(t)$, $C(t)$ as in theorem  in Theorem~\ref{thm:comptestm},
  assume that there exists $\alpha \in \mathbb{R}$ such that
  \begin{equation*}
    \mu_{\rm{log}}(A(t)) \leq \alpha, \quad t \in [0,T].
\end{equation*}
Then
\begin{equation}
  |z(t)| \leq e^{\alpha t}|z(0)| + \int_0^t e^{\alpha(t-s)}C(s)ds.  \label{eq:ln-duhamel}
\end{equation}
\end{theorem}

In the previous results, the choice of norms has been arbitrary. We now apply these results to the case where the norm is Euclidean. In this case we have the following formula for the logarithmic norm of a matrix.

\begin{lemma}
If we choose the Euclidean norm for the computation of $\mu_{\rm{log}}(A)$, then
\begin{eqnarray}
\mu_{\rm{log}}(A) =\max \{\lambda \in \text{spectrum of }(A+A^{\top })/2\} =  \max_{\|x\|=1} (Ax,x)
\label{eq:eucl-log-norm} \\
m_l(A) =\min \{\lambda \in \text{spectrum of }(A+A^{\top })/2\} = \min_{\|x\|=1} (Ax,x)  \label{eq:eucl-m-l-norm}.
\end{eqnarray}
\end{lemma}

\section{Some results about normally hyperbolic invariant manifolds}
\label{sec:NHIM}

The content of this section is taken by adapting the results of \cite{CZ15} for the maps case and of \cite{CZ17} for the ODE case. In those papers, the normally hyperbolic invariant manifold was $\Lambda=\mathbb{T}^c$, while here we consider $\Lambda=\mathbb{R}^c$.

In \cite{CZ15,CZ17} the covering map for a $c$-dimensional torus
\begin{equation*}
\varphi :\mathbb{R}^{c}\rightarrow \Lambda = \mathbb{R}/( (2 R_\Lambda) \cdot \mathbb{Z})^{c},
\end{equation*}%
provided a global coordinate system for the whole torus.
Given a map $f$ or a vector field on $\Lambda \times \mathbb{R}^u \times \mathbb{R}^s$, it can be lifted to $\mathbb{R}^c \times \mathbb{R}^u \times \mathbb{R}^s$
and all the assumptions in \cite{CZ15,CZ17} are formulated in terms of the 'lifted' $f$.
The lifted $f$ are periodic in the first variable,
but nevertheless, all the techniques and conclusions also apply to the situation when we have uniform bounds in the center direction (i.e., $\mathbb{R}^c$), and are therefore applicable
to global maps or vector fields defined on $\mathbb{R}^c \times \mathbb{R}^u \times \mathbb{R}^s$. In addition, some topological conditions that appear in the case of the torus
can be dropped in the global case.

Throughout the this section we shall use the notation $z=(\lambda ,x,y)$ to denote
points in $\mathbb{R}^c \times \mathbb{R}^u \times \mathbb{R}^s$. This means that the notation $\lambda $ will stand for points on $%
\Lambda $, $x$ for points in $\mathbb{R}^{u}$, and $y$ for points
in $\mathbb{R}^{s}$. We will write $f$ as $(f_{\lambda },f_{x},f_{y})$,
where $f_{\lambda },f_{x},f_{y}$ stand for projections onto $\Lambda $, $%
\mathbb{R}^{u}$ and $\mathbb{R}^{s}$, respectively. On $\mathbb{R}^{c}\times
\mathbb{R}^{u}\times \mathbb{R}^{s}$ we will use the Euclidian norm.

\subsection{NHIMs for maps}

Given $U \subset  \mathbb{R}^c \times \mathbb{R}^{u}\times \mathbb{R}^{s}$
 we consider a $C^{k+1}$ map, for $k\geq 1$,
\begin{equation*}
f:U \rightarrow \mathbb{R}^c \times \mathbb{R}^{u}\times \mathbb{R}^{s}.
\end{equation*}
In \cite{CZ15} where we considered tori, a number $R_\Lambda$ being the radius (in fact period of lifted $f$ in $\lambda$-direction)
was introduced. Now when considering $\Lambda=\mathbb{R}^n$ we can take $R_\Lambda=\infty$.

For any $R>0$ (in \cite{CZ15} it was required that $R<\frac{1}{2}R_{\Lambda }$) denote by $D=D(R)$ the set
\begin{equation*}
D(R)=\Lambda \times \overline{B}_{u}(R)\times \overline{B}_{s}(R),
\end{equation*}%
where $\overline{B}_{n}(R)$ stands for a closed ball of radius $R$, centered
at zero, in $\mathbb{R}^{n}$.

We will now define several constants $\xi$ and $\mu$ related to the expansion and contraction rates of $f$ on $D$ along several directions.
These constants are defined on the set $D$ and depend on $L$ (this will measure slopes) and the notion of $P(z)$ for $z \in D$.

In \cite[page 6220]{CZ15} it was asked that $L \in \left( \frac{2R}{R_{\Lambda }},1\right)$, but in our global case we can have arbitrary $ L$ and we can set
\begin{equation*}
 P(q)=D.
\end{equation*}
The restriction on $L$ and the notion of $P(q)$ were introduced in order to have a projection of cones with slopes $L$ or $1/L$ to be on the same chart (defined by $\varphi$)
on $\Lambda$ and $P(q)$ was a preimage of all points that can be the same chart as point $q$. For the same reason we also drop the backward cone condition (Def. 14 in \cite{CZ15}).

Before presenting definition of various constants that describe expansion and contraction rates we need to introduce one  more notation.
\begin{definition}
Let $f:\mathbb{R}^{n}\rightarrow \mathbb{R}^{k}$ be a $C^{1}$ function. We
define the interval enclosure of the derivative $Df$ on $U\subset \mathbb{R}%
^{n}$ as the set $[Df(U)]\subset \mathbb{R}^{k\times n}$, defined as%
\begin{equation*}
\lbrack Df(U)]=\left\{ A=\left( a_{ij}\right) _{\substack{ i=1,\ldots ,k  \\ %
j=1,\ldots ,n}}:a_{ij}\in \left[ \inf_{x\in U}\frac{\partial f_{j}}{\partial
x_{j}}(x),\sup_{x\in U}\frac{\partial f_{j}}{\partial x_{j}}(x)\right]
\right\} .
\end{equation*}
\end{definition}

Recall that in (\ref{eq:minA}) and (\ref{eq:minA-intv}) we already defined $m(A)$ and $m([DF(U)])$. These notations are used below.

For $L>0$ and the set $D$ we  define (see \cite[page 6220]{CZ15})
\begin{align*}
\mu _{s,1}& =\sup_{z\in D}\left\{ \left\Vert \frac{\partial f_{y}}{\partial y%
}\left( z\right) \right\Vert +\frac{1}{L}\left\Vert \frac{\partial f_{y}}{%
\partial (\lambda ,x)}(z)\right\Vert \right\} , \\
\mu _{s,2}& =\sup_{z\in D}\left\{ \left\Vert \frac{\partial f_{y}}{\partial y%
}\left( z\right) \right\Vert +L\left\Vert \frac{\partial f_{\left( \lambda
,x\right) }}{\partial y}(z)\right\Vert \right\} ,
\end{align*}%
\begin{align*}
\xi _{u,1}& =\inf_{z\in D}\left\{ m\left( \frac{\partial f_{x}}{\partial x}%
(z)\right) -\frac{1}{L}\left\Vert \frac{\partial f_{x}}{\partial \left(
\lambda ,y\right) }(z)\right\Vert \right\} , \\
\xi _{u,1,P}& =\inf_{z\in D}m\left[ \frac{\partial f_{x}}{\partial x}(P(z))%
\right] -\frac{1}{L}\sup_{z\in D}\left\Vert \frac{\partial f_{x}}{\partial
\left( \lambda ,y\right) }(z)\right\Vert , \\
\xi _{u,2}& =\inf_{z\in D}\left\{ m\left( \frac{\partial f_{x}}{\partial x}%
\left( z\right) \right) -L\left\Vert \frac{\partial f_{(\lambda ,y)}}{%
\partial x}(z)\right\Vert \right\} ,
\end{align*}%
\begin{align}
\mu _{cs,1}& =\sup_{z\in D}\left\{ \left\Vert \frac{\partial f_{\left(
\lambda ,y\right) }}{\partial \left( \lambda ,y\right) }(z)\right\Vert
+L\left\Vert \frac{\partial f_{\left( \lambda ,y\right) }}{\partial x}%
(z)\right\Vert \right\} ,  \notag \\
\mu _{cs,2}& =\sup_{z\in D}\left\{ \left\Vert \frac{\partial f_{\left(
\lambda ,y\right) }}{\partial \left( \lambda ,y\right) }(z)\right\Vert +%
\frac{1}{L}\left\Vert \frac{\partial f_{x}}{\partial \left( \lambda
,y\right) }(z)\right\Vert \right\} ,  \notag
\end{align}%
\begin{align*}
\xi _{cu,1}& =\inf_{z\in D}\left\{ m\left( \frac{\partial f_{(\lambda ,x)}}{%
\partial (\lambda ,x)}(z)\right) -L\left\Vert \frac{\partial f_{(\lambda ,x)}%
}{\partial y}(z)\right\Vert \right\} , \\
\xi _{cu,1,P}& =\inf_{z\in D}m\left[ \frac{\partial f_{(\lambda ,x)}}{%
\partial (\lambda ,x)}(P(z))\right] -L\sup_{z\in D}\left\Vert \frac{\partial
f_{(\lambda ,x)}}{\partial y}(z)\right\Vert , \\
\xi _{cu,2}& =\inf_{z\in D}\left\{ m\left( \frac{\partial f_{(\lambda ,x)}}{%
\partial \left( \lambda ,x\right) }(z)\right) -\frac{1}{L}\left\Vert \frac{%
\partial f_{y}}{\partial (\lambda ,x)}(z)\right\Vert \right\} .
\end{align*}

\begin{definition}\cite[Def. 5]{CZ15}
\label{def:rate-conditions}We say that $f$ satisfies rate conditions of
order $k\geq 1$ if $\xi _{u,1},$ $\xi _{u,1,P},$ $\xi _{u,2},$ $\xi _{cu,1},$
$\xi _{cu,1,P},$ $\xi _{cu,2}$ are strictly positive, and for all $k\geq
j\geq 1$ it holds that
\begin{equation*}
\mu _{s,1}<1<\xi _{u,1,P},  
\end{equation*}%
\begin{align*}
\frac{\mu _{cs,1}}{\xi _{u,1,P}}& <1,\qquad \frac{\mu _{s,1}}{\xi _{cu,1,P}}<1, \\
\frac{\left( \mu _{cs,1}\right) ^{j+1}}{\xi _{u,2}}& <1,\qquad \frac{\mu
_{s,2}}{(\xi _{cu,1})^{j+1}}<1, \nonumber \\
\frac{\mu _{cs,2}}{\xi _{u,1}}& <1,\qquad \frac{\mu _{s,1}}{\xi _{cu,2}}<1.\nonumber
\end{align*}
\end{definition}

We introduce the following notation:
\begin{eqnarray*}
J_{s}(z,M) &=&\left\{ \left( \lambda ,x,y\right) :\left\Vert \left( \lambda
,x\right) -\pi _{\lambda ,x}z\right\Vert \leq M\left\Vert y-\pi
_{y}z\right\Vert \right\} , \\
J_{u}\left( z,M\right) &=&\left\{ \left( \lambda ,x,y\right) :\left\Vert
\left( \lambda ,y\right) -\pi _{\lambda ,y}z\right\Vert \leq M\left\Vert
x-\pi _{x}z\right\Vert \right\} .  
\end{eqnarray*}%
We shall refer to $J_{s}(z,M)$ as a stable cone of slope $M$ at $z$, and to
$J_{u}(z,M)$ as an unstable cone of slope $M$ at $z$.

\begin{definition}
We say that a sequence $\left\{ z_{i}\right\} _{i=-\infty }^{0}$ is a (full)
backward trajectory of a point $z$ if $z_{0}=z,$ and $f\left( z_{i-1}\right)
=z_{i}$ for all $i\leq 0.$
\end{definition}

\begin{definition}
We define the local center-stable set in $D$ as%
\begin{equation*}
W^{cs}_{\mathrm{loc}}=\left\{ z:f^{n}(z)\in D\text{ for all }n\in \mathbb{N}\right\} .
\end{equation*}
\end{definition}

\begin{definition}
We define the local center-unstable set in $D$ as%
\begin{equation*}
W^{cu}_{\mathrm{loc}}=\{z:\text{there is a full backward trajectory of }z\text{ in }D\}.
\end{equation*}
\end{definition}

\begin{definition}
We define the maximal invariant set in $D$ as
\begin{equation*}
\Lambda ^{\ast }_{\mathrm{loc}}=\{z:\text{there is a full trajectory of }z\text{ in }D\}.
\end{equation*}
\end{definition}

\begin{definition}
Assume that $z\in W^{cs}_{\mathrm{loc}}$. We define the local stable fiber
of $z$ as
\begin{equation*}
W_{z,\mathrm{loc}}^{s}=\left\{ p\in D:f^{n}\left( p\right) \in J_{s}\left(
f^{n}(z),1/L\right) \cap D\text{ for all }n\in \mathbb{N}\right\} .
\end{equation*}
\end{definition}

\begin{definition}
Assume that $z\in W^{cu}_{\mathrm{loc}}$. We define the local unstable fiber
of $z$ as%
\begin{eqnarray*}
W_{z,\mathrm{loc}}^{u} &=&\left\{p\in D:\exists \textnormal{ backward trajectory }\left\{
p_{i}\right\} _{i=-\infty }^{0}\text{ of }p\text{ in }D,
\textnormal{ and for any such }\right.\\
&&\left. \textnormal{backward trajectory }\left\{ z_{i}\right\} _{i=-\infty
}^{0}\textnormal{ of }z\textnormal{ in }D \textnormal{ it holds that}\left. p_{i}\in J_{u}\left( z_{i},1/L\right) \cap
D\right. \right\}.
\end{eqnarray*}
\end{definition}

The definitions of $W_{z,\mathrm{loc}}^{s}$ and $W_{z,\mathrm{loc}}^{u}$ are related to cones, which
is a nonstandard approach, as the standard one is through convergence rates. In \cite{CZ15}
it is shown that it implies the convergence rate as in the
standard theory.

Under above assumptions it will turn out that $f$ is injective on $W^{cu}_{\mathrm{loc}}$.
Therefore the backward orbit in the definition of $W_{z,\mathrm{loc}}^{u}$ is unique.

In \cite[page 6223]{CZ15} for given $\lambda \in \Lambda$ the set $D_\lambda^\pm$ was defined. Now when $\Lambda=\mathbb{R}^c$ we do not need to
`localize' $D$ in a good chart and instead we define the following sets:
\begin{eqnarray*}
D^{+} &=&\mathbb{R}^c
\times \overline{B}_{u}(R)\times \partial B_{s}(R), \\
D^{-} &=&\mathbb{R}^c
\times \partial \overline{B}_{u}(R)\times B_{s}(R).
\end{eqnarray*}

This is a `global' modification of Def. 15 in \cite{CZ15}
\begin{definition}
\label{def:covering}We say that $f$ satisfies covering conditions (on $D$) if  there exist a homotopy $h$
\begin{equation*}
h:[0,1] \times D\rightarrow \mathbb{R}^c \times \mathbb{R}^{u}\times \mathbb{R}^{s},
\end{equation*}%
a $\lambda^\ast \in \mathbb{R}^c$ and a linear map $A:\mathbb{R}^{u}\rightarrow \mathbb{R}^{u}$ which satisfy:

\begin{enumerate}
\item $h_{0}=f|_{D},$

\item for any $\alpha \in \left[ 0,1\right] $,
\begin{eqnarray*}
h_{\alpha }\left( D^-\right) \cap D  &=&\emptyset , \\
h_{\alpha }(D) \cap D^{+} &=&\emptyset ,
\end{eqnarray*}

\item $h_{1}\left( \lambda ,x,y\right) =\left( \lambda ^{\ast },Ax,0\right) $%
,

\item $A\left( \partial B_{u}(R)\right) \subset \mathbb{R}^{u}\setminus
\overline{B}_{u}(R).$
\end{enumerate}
\end{definition}

The following theorem is Theorem 16 in \cite{CZ15}. The only change is that now $\Lambda=\mathbb{R}^c$, we  assume the boundedness of $\|D^jf\|$ for $j=1,\dots,k+1$ on $D$
and we drop the backward cone condition.
\begin{theorem}
\label{thm:nhim-map} Let $\Lambda=\mathbb{R}^c$, $k\geq 1$ and $f:D\rightarrow \Lambda
\times \mathbb{R}^{u}\times \mathbb{R}^{s}$ be a $C^{k+1}$ map. If $f$
satisfies covering conditions and rate conditions of order $k$, and $\|D^j f\|$ are bounded on $D$ for $j=1,\dots,k$,
then $W^{cs}_{\mathrm{loc}},W^{cu}_{\mathrm{loc}}$ and $\Lambda ^{\ast }$ are $C^{k}$
manifolds, which are graphs of $C^{k}$ functions
\begin{align*}
w^{cs}_{\mathrm{loc}}& :\Lambda \times \overline{B}_{s}(R)\rightarrow \overline{B}_{u}(R),
\\
w^{cu}_{\mathrm{loc}}& :\Lambda \times \overline{B}_{u}(R)\rightarrow \overline{B}_{s}(R),
\\
\chi & :\Lambda \rightarrow \overline{B}_{u}(R)\times \overline{B}_{s}(R),
\end{align*}%
meaning that
\begin{align*}
W^{cs}_{\mathrm{loc}}& =\left\{ \left( \lambda ,w^{cs}_{\mathrm{loc}}(\lambda ,y),y\right) :\lambda \in
\Lambda ,y\in \overline{B}_{s}(R)\right\} , \\
W^{cu}_{\mathrm{loc}}& =\left\{ \left( \lambda ,x,w^{cu}_{\mathrm{loc}}(\lambda ,y)\right) :\lambda \in
\Lambda ,x\in \overline{B}_{u}(R)\right\} , \\
\Lambda ^{\ast }& =\left\{ \left( \lambda ,\chi (\lambda )\right) :\lambda
\in \Lambda \right\} .
\end{align*}%
Moreover, if $f$ restricted to $W^{cu}_{\mathrm{loc}}$ is an injection, $w^{cs}_{\mathrm{loc}}$ and $w^{cu}_{\mathrm{loc}}$ are Lipschitz
with constants $L$, and $\chi $ is Lipschitz with the constant $\sqrt{2}L/\sqrt{1-L^{2}}$. The manifolds $W^{cs}_{\mathrm{loc}}$ and $W^{cu}_{\mathrm{loc}}$ intersect
transversally on $\Lambda ^{\ast }$, and $W^{cs}_{\mathrm{loc}}\cap W^{cu}_{\mathrm{loc}}=\Lambda ^{\ast }$.

The manifolds $W^{cs}_{\mathrm{loc}}$ and $W^{cu}_{\mathrm{loc}}$ are foliated by invariant fibers $%
W_{z,\mathrm{loc}}^{s}$ and $W_{z,_\mathrm{loc}}^{u}$, which are graphs of $%
C^{k}$ functions%
\begin{eqnarray*}
w_{z,\mathrm{loc}}^{s} &:&\overline{B}_{s}(R)\rightarrow \Lambda \times \overline{B}%
_{u}(R), \\
w_{z,\mathrm{loc}}^{u} &:&\overline{B}_{u}(R)\rightarrow \Lambda \times \overline{B}%
_{s}(R),
\end{eqnarray*}%
meaning that%
\begin{eqnarray*}
W_{z,\mathrm{loc}}^{s} &=&\left\{ \left( w_{z,\mathrm{loc}}^{s}\left( y\right) ,y\right) :y\in
\overline{B}_{s}(R)\right\} , \\
W_{z,\mathrm{loc}}^{u} &=&\left\{ \left( \pi _{\lambda }w_{z,\mathrm{loc}}^{u}\left( x\right) ,x,\pi
_{y}w_{z,\mathrm{loc}}^{u}\left( x\right) \right) :x\in \overline{B}_{u}(R)\right\} .
\end{eqnarray*}%
The functions $w_{z,\mathrm{loc}}^{s}$ and $w_{z,\mathrm{loc}}^{u}$ are Lipschitz with constants $1/L$. Moreover,
\begin{eqnarray*}
W_{z,\mathrm{loc}}^{s} &=&\left\{p\in D:f^{n}(p)\in D\textnormal{ for all }n\geq 0,\textnormal{ and }
 \exists n_{0},\exists C>0\textnormal{ (which can depend on p) }\right.\\
&&\left.\textnormal{ s.t. for }n\geq n_{0}, \left\Vert f^{n}(p)-f^{n}(z)\right\Vert \leq C\mu
_{s,1}^{n}\right\} ,
\end{eqnarray*}%
and if $\{z_{i}\}_{i=-\infty }^{0}$ is the unique backward trajectory of $z$
in $D$, then
\begin{eqnarray*}
W_{z,\mathrm{loc}}^{u} &=&\left\{p\in W^{cu}_{\mathrm{loc}}:\textnormal{such that the unique backward trajectory }
\{p_{i}\}_{i=-\infty }^{0}\textnormal{ of }p\textnormal{ in }D\textnormal{ satisfies that }\right. \\
&&\left. \exists n_{0}\geq 0,\exists C>0\textnormal{ (which can depend on p) s.t. for }
n\geq n_{0}, \left\Vert p_{-n}-z_{-n}\right\Vert \leq C\xi _{u,1,P}^{-n}\right\}.
\end{eqnarray*}
\end{theorem}
\begin{rem}
The above theorem says nothing about the dependence of $W_{z,\mathrm{loc}}^{s,u}$ on $z$.
\end{rem}

From the proof of Theorem~\ref{thm:nhim-map} we immediately obtain the following result.
\begin{rem}
In the context of Theorem~\ref{thm:nhim-map} the following holds. For $j=1,\dots,k$ the derivatives of the functions
$w^{cu}_{\mathrm{loc}}$, $w^{cs}_{\mathrm{loc}}$, $w^{u,s}_{z,\mathrm{loc}}$ of order $j$
are bounded by some constants $K_j$ depending on $\|D^s f\|$ for $s=1,\dots,j$ and constants $L$, $\xi$ and $\mu$.
\end{rem}
While in \cite{CZ15} there is no explicit general formula for $K_j$, the  explicit bound for second derivatives has been worked out in Theorem 23 in \cite{CZ17}.

Observe that we obtain $L$ or $1/L$ as  bounds for the Lipschitz constants for the functions $w^{cu}_{\mathrm{loc}}$, $w^{cs}_{\mathrm{loc}}$,
$w^{u}_{\mathrm{loc}}$, $w^{s}_{\mathrm{loc}}$. Hence when $L$ is small we get quite big bounds for some slopes. This is clearly an overestimate for the case when $
\mathbb{T}\times \{0\}\times \{0\}$ is our NHIM.  This a
consequence of the choices that have been made when formulating Theorem \ref{thm:nhim-map},
as the authors in \cite{CZ15} did not introduce different constants for each type of cones,
plus several inequalities between them. However, below following \cite{CZ15} we give conditions
which allow us to obtain better Lipschitz constants.

\begin{theorem} \cite[Thm. 17]{CZ15}
\begin{eqnarray*}
\mu &=&\sup_{z\in D}\left\{ \left\Vert \frac{\partial f_{y}}{\partial y}%
\left( z\right) \right\Vert +M\left\Vert \frac{\partial f_{y}}{\partial
(\lambda ,x)}(z)\right\Vert \right\} , \\
\xi &=&\inf_{z\in D}m\left( \left[ \frac{\partial f_{(\lambda ,x)}}{\partial
\left( \lambda ,x\right) }(P(z))\right] \right) -\frac{1}{M}\sup_{z\in
D}\left\Vert \frac{\partial f_{(\lambda ,x)}}{\partial y}(z)\right\Vert .
\end{eqnarray*}%
If assumptions of Theorem \ref{thm:nhim-map} hold true and also $\xi/\mu>1$,
 then the function $w_{z,\mathrm{loc}}^{s}$ from Theorem \ref{thm:nhim-map} is Lipschitz
with constant $M.$
\end{theorem}

\begin{theorem} \cite[Thm. 18]{CZ15}
\begin{eqnarray*}
\xi &=&\inf_{z\in D}m\left( \left[ \frac{\partial f_{x}}{\partial x}(P(z))%
\right] \right) -M\sup_{z\in D}\left\Vert \frac{\partial f_{x}}{\partial
\left( \lambda ,y\right) }(z)\right\Vert , \\
\mu &=&\sup_{z\in D}\left\{ \left\Vert \frac{\partial f_{\left( \lambda
,y\right) }}{\partial \left( \lambda ,y\right) }(z)\right\Vert +\frac{1}{M}%
\left\Vert \frac{\partial f_{\left( \lambda ,y\right) }}{\partial x}%
(z)\right\Vert \right\} .
\end{eqnarray*}%
If assumptions of Theorem \ref{thm:nhim-map} hold true and also $\xi/\mu>1$,
 then the function $w_{z,\mathrm{loc}}^{u}$ from Theorem \ref{thm:nhim-map} is Lipschitz
with constant $M.$
\end{theorem}

\begin{theorem} \cite[Thm. 19]{CZ15}
\begin{eqnarray*}
\xi &=&\inf_{z\in D}m\left[ \frac{\partial f_{(\lambda ,x)}}{\partial
(\lambda ,x)}(P(z))\right] -M\sup_{z\in D}\left\Vert \frac{\partial
f_{(\lambda ,x)}}{\partial y}(z)\right\Vert , \\
\mu &=&\sup_{z\in D}\left\{ \left\Vert \frac{\partial f_{y}}{\partial y}%
\left( z\right) \right\Vert +\frac{1}{M}\left\Vert \frac{\partial f_{y}}{%
\partial (\lambda ,x)}(z)\right\Vert \right\} .
\end{eqnarray*}%
If assumptions of Theorem \ref{thm:nhim-map} hold true and also $\xi/\mu>1$, then the function $w^{cu}_{\mathrm{loc}}$ from Theorem \ref{thm:nhim-map} is Lipschitz with
constant $M$.
\end{theorem}

\begin{theorem} \cite[Thm. 20]{CZ15}
\begin{eqnarray*}
\xi &=&\inf_{z\in D}m\left[ \frac{\partial f_{x}}{\partial x}(P(z))\right] -%
\frac{1}{M}\sup_{z\in D}\left\Vert \frac{\partial f_{x}}{\partial \left(
\lambda ,y\right) }(z)\right\Vert , \\
\mu &=&\sup_{z\in D}\left\{ \left\Vert \frac{\partial f_{\left( \lambda
,y\right) }}{\partial \left( \lambda ,y\right) }(z)\right\Vert +M\left\Vert
\frac{\partial f_{\left( \lambda ,y\right) }}{\partial x}(z)\right\Vert
\right\}
\end{eqnarray*}%
If assumptions of Theorem \ref{thm:nhim-map} hold true and also $\xi/\mu>1$, then the function $w^{cs}_{\mathrm{loc}}$ from Theorem \ref{thm:nhim-map} is Lipschitz with
constant $M$.
\end{theorem}

\subsubsection{Comments on the inequalities}

Let $J_{s}^{c}(z,M)$ and $J_{u}^{c}(z,M)$ stand for the complements of $J_{s}(z,M)$ and $J_{u}(z,M),$ respectively.
We now discuss the meaning of several inequalities on Definition \ref{def:rate-conditions} of rate
conditions as well as where they are needed (see section 3.3 in \cite{CZ15} for more details):

\begin{itemize}
\item $\mu _{cs,1}<\xi _{u,1,P}$: the forward invariance of $J_{u}(z,1/L)$. $\xi _{u,1,P}>1$: the expansion in $J_{u}(z,1/L)$ for $x$ - coordinate. This is
needed for the proof of the existence of $W^{cs}$.

\item $\xi _{cu,1,P}>\mu _{s,1}$: the forward invariance of $%
J_{s}^{c}(z,1/L) $. $\mu _{s,1}<1$: the
contraction in $y$-direction in $J_{s}(z,1/L)$. This is needed for the proof of the existence of $W^{cu}$.

\item $\frac{\mu _{s,2}}{(\xi _{cu,1})^{j+1}}<1$, $j=1,\dots ,k$: the $C^{k}$-smoothness of $W^{cu}$.

\item $\frac{\left( \mu _{cs,1}\right) ^{j+1}}{\xi _{u,2}}<1$, $j=1,\dots ,k$
: the $C^{k}$-smoothness of $W^{cs}$.

\item $\frac{\mu _{cs,1}}{\xi _{u,1,P}}<1$: the existence of fibers $W_{q}^{u}$. $\frac{\mu _{cs,2}}{\xi _{u,1}}<1$:
the $C^{k}$ smoothness of $W_{q}^{u}$.

\item $\frac{\mu _{s,1}}{\xi _{cu,1,P}}<1$: the existence of fibers $W_{q}^{s}$. $\frac{\mu _{s,1}}{\xi _{cu,2}}<1$:
the $C^{k}$ smoothness of $W_{q}^{s}$.
\end{itemize}

\subsection{NHIMs for ODEs}
\label{subsec:NHIM-ode}
The content of this section  is partially based on Section 5 in \cite{CZ17}, where the results from \cite{CZ15} about maps were adapted to the context
of ODEs. The goal of this section is to describe how these results extend to the case of $\Lambda=\mathbb{R}^c$.

We consider an ODE
\begin{equation}
q^{\prime }=f(q) , \label{eq:ode-Wcu}
\end{equation}%
where%
\begin{equation*}
f:\Lambda \times \mathbb{R}^{u}\times \mathbb{R}^{s}\rightarrow \mathbb{R}%
^{c}\times \mathbb{R}^{u}\times \mathbb{R}^{s}.
\end{equation*}%
We denote by $\Phi \left( t,q\right) $ the flow induced by (\ref{eq:ode-Wcu}).

We  define the notion of an isolating block, which is an ODE version of Def.~\ref{def:covering}.

\begin{definition} \cite[Def. 19]{CZ17}
\label{def:isolating-segment}
We say that $D(r)=\Lambda \times \overline{B}_{u}(r) \times \overline{B}_{s}(r) $ is an isolating block for $f$ if

\begin{enumerate}
\item Exit.  For any $q\in \Lambda \times \partial \overline{B}_{u}(r) \times \overline{B}_{s}(r) $,
\begin{equation}
\left( \pi _{x}F_\varepsilon(q)|\pi _{x}q\right) >0. \label{eq:isoblock-exit}
\end{equation}

\item Entry. For any $q\in \Lambda \times \overline{B}_{u}(r) \times \partial \overline{B}_{s}(r) $,%
\begin{equation}
\left( \pi _{y}F_\varepsilon(q)|\pi _{y}q\right) <0. \label{eq:isoblock-entry}
\end{equation}
\end{enumerate}
\end{definition}

Isolating blocks are important constructs in the Conley index theory \cite{MM}. Intuitively,
in Definition~\ref{def:isolating-segment} the set $\Lambda \times \partial \overline{B}_{u}(r) \times \overline{B}_{s}(r) $ plays the role of the exit set, and $\Lambda \times
\overline{B}_{u}(r) \times \partial \overline{B}_{s}(r) $ of the entry set.

\begin{definition}
We define the center-unstable set of ode (\ref{eq:ode-Wcu})
in $D$ as
\begin{equation*}
W_{\mathrm{loc},D}^{cu}=\{q:\Phi( t,q) \in D\text{ for all }%
t<0\},
\end{equation*}
and the center-stable set of ode (\ref{eq:ode-Wcu})
in $D$ as
\begin{equation*}
W_{\mathrm{loc},D}^{cs}=\{q:\Phi( t,q) \in D\text{ for all }%
t>0\}.
\end{equation*}
\end{definition}

We shall consider the set
\begin{equation*}
D=\Lambda \times \overline{B}_{u}\left( R\right) \times \overline{B}%
_{s}\left( R\right) .
\end{equation*}

Since the set $D$ will remain fixed throughout the discussion, from now on
we will simplify notation by writing $W^{cu}$ instead of $W_{\mathrm{loc}%
,D}^{cu}.$

We will now define the constants corresponding to $\xi$ and $\mu$ for the map case.

The rule is that the corresponding
constant for $\xi$ will be $\overrightarrow{\xi}$, and analogously for $\mu$.

The definitions are such that if we consider a map $\Phi(h,\cdot)$ i.e., the flow induced by (\ref{eq:ode-Wcu}) by $h>0$,
then  for all constants $\xi$ and $\mu$ it holds that (see Thm. 31 in \cite{CZ17})
\begin{eqnarray*}
  \mu = 1 + h \overrightarrow{\mu} + O(h^2), \\
  \xi = 1 + h \overrightarrow{\xi} + O(h^2).
\end{eqnarray*}

Therefore following \cite{CZ17} we introduce the following constants ($\mu_{\rm{log}}()$ and $m_l()$ are defined in Appendix~\ref{app:sec-lognorm})
\begin{align*}
\overrightarrow{\mu _{s,1}}& =\sup_{z\in D}\left\{ \mu_{\rm{log}}\left( \frac{\partial
f_{y}}{\partial y}(z)\right) +\frac{1}{L}\left\Vert \frac{\partial f_{y}}{%
\partial (\lambda ,x)}(z)\right\Vert \right\} , \\
\overrightarrow{\mu _{s,2}}& =\sup_{z\in D}\left\{ \mu_{\rm{log}}\left( \frac{\partial
f_{y}}{\partial y}\left( z\right) \right) +L\left\Vert \frac{\partial
f_{\left( \lambda ,x\right) }}{\partial y}(z)\right\Vert \right\} ,
\end{align*}%
\begin{align*}
\overrightarrow{\xi _{u,1}}& =\inf_{z\in D}\left\{ m_{l}\left( \frac{%
\partial f_{x}}{\partial x}(z)\right) -\frac{1}{L}\left\Vert \frac{\partial
f_{x}}{\partial \left( \lambda ,y\right) }(z)\right\Vert \right\} , \\
\overrightarrow{\xi _{u,1,P}}& =\inf_{z\in D}m_{l}\left( \frac{\partial f_{x}%
}{\partial x}(P(z))\right) -\frac{1}{L}\sup_{z\in D}\left\Vert \frac{%
\partial f_{x}}{\partial \left( \lambda ,y\right) }(z)\right\Vert ,
\end{align*}%
\begin{align}
\overrightarrow{\mu _{cs,1}}& =\sup_{z\in D}\left\{ \mu_{\rm{log}}\left( \frac{\partial
f_{\left( \lambda ,y\right) }}{\partial \left( \lambda ,y\right) }(z)\right)
+L\left\Vert \frac{\partial f_{\left( \lambda ,y\right) }}{\partial x}%
(z)\right\Vert \right\} ,  \notag \\
\overrightarrow{\mu _{cs,2}}& =\sup_{z\in D}\left\{ \mu_{\rm{log}}\left( \frac{\partial
f_{\left( \lambda ,y\right) }}{\partial \left( \lambda ,y\right) }(z)\right)
+\frac{1}{L}\left\Vert \frac{\partial f_{x}}{\partial \left( \lambda
,y\right) }(z)\right\Vert \right\} ,  \notag
\end{align}%
\begin{align*}
\overrightarrow{\xi _{cu,1}}& =\inf_{z\in D}\left\{ m_{l}\left( \frac{%
\partial f_{(\lambda ,x)}}{\partial (\lambda ,x)}(z)\right) -L\left\Vert
\frac{\partial f_{(\lambda ,x)}}{\partial y}(z)\right\Vert \right\} , \\
\overrightarrow{\xi _{cu,1,P}}& =\inf_{z\in D}m_{l}\left( \frac{\partial
f_{(\lambda ,x)}}{\partial (\lambda ,x)}(P(z))\right) -L\sup_{z\in
D}\left\Vert \frac{\partial f_{(\lambda ,x)}}{\partial y}(z)\right\Vert, \\
\overrightarrow{\xi _{u,2}}& =\inf_{z\in D}\left\{ m_{l}\left( \frac{\partial f_x}{\partial x}(z)\right) -L\left\Vert\frac{\partial f_{(\lambda ,y)}}{\partial x}(z)\right\Vert \right\} , \\
\overrightarrow{\xi _{cu,2}}& =\inf_{z\in D}m_{l}\left( \frac{\partial f_{(\lambda ,x)}}{\partial (\lambda ,x)}(z)\right) -\frac{1}{L}\sup_{z\in
D}\left\Vert \frac{\partial f_y}{\partial (\lambda ,x)}(z)\right\Vert .
\end{align*}

We
define the rate conditions for ODEs as follows.
\begin{definition} \cite[Def. 28]{CZ17}
\label{def:rate-cond-ode}We say that the vector field $f$ satisfies rate
conditions of order $k\geq 1$ if for all $k\geq j\geq 1$ it holds that
\begin{equation*}
\overrightarrow{\mu _{s,1}}<0<\overrightarrow{\xi _{u,1,P}},
\end{equation*}%
\begin{gather}
\overrightarrow{\mu _{cs,1}}<\overrightarrow{\xi _{u,1,P}},\qquad \qquad
\overrightarrow{\mu _{s,1}}<\overrightarrow{\xi _{cu,1,P}},\nonumber\\
 \overrightarrow{\mu _{s,2}}<(j+1) \overrightarrow{\xi _{cu,1}},\qquad \overrightarrow{\mu _{cs,2}}< \overrightarrow{\xi _{u,1}}, \nonumber \\
 (j+1) \overrightarrow{\mu_{cs,1}} < \overrightarrow{\xi_{u,2}} \qquad \overrightarrow{\mu_{s,1}} < \overrightarrow{\xi_{cu,2}}. \label{eq:rate-cond-4-ode}
\end{gather}
\end{definition}
In fact in \cite[Def. 28]{CZ17} conditions (\ref{eq:rate-cond-4-ode})are missing. Because there the focus was just on $W^{cu}$, and
dealing with $W^{cs}$ requires (\ref{eq:rate-cond-4-ode}).

There is no ODE-analogue of Theorem~\ref{thm:nhim-map} in \cite{CZ17}, only a part of it related to $W^{cu}$ is stated there as Theorem 30. However, the full result stated
below can in the same way be derived from Theorem~\ref{thm:nhim-map}.
\begin{theorem}
\label{thm:nhim-ode}Let $k\geq 1$ and assume that $f$ is $C^{k+1}$.  Assume  that $D=\Lambda \times \overline{B}%
_{u}\left( R\right) \times \overline{B}_{s}\left( R\right) $ is an isolating block for $f$ and $f$ has bounded derivatives up to order $k+1$ on $D$ and
rate conditions of order $k$ are satisfied on $D$ for some constant $L$.

Then $W^{cs},W^{cu}$ and $\Lambda ^{\ast }$ are $C^{k}$
manifolds, which are graphs of $C^{k}$ functions
\begin{align*}
w^{cs}& :\Lambda \times \overline{B}_{s}(R)\rightarrow \overline{B}_{u}(R),
\\
w^{cu}& :\Lambda \times \overline{B}_{u}(R)\rightarrow \overline{B}_{s}(R),
\\
\chi & :\Lambda \rightarrow \overline{B}_{u}(R)\times \overline{B}_{s}(R),
\end{align*}%
meaning that
\begin{align*}
W^{cs}& =\left\{ \left( \lambda ,w^{cs}(\lambda ,y),y\right) :\lambda \in
\Lambda ,y\in \overline{B}_{s}(R)\right\} , \\
W^{cu}& =\left\{ \left( \lambda ,x,w^{cu}(\lambda ,y)\right) :\lambda \in
\Lambda ,x\in \overline{B}_{u}(R)\right\} , \\
\Lambda ^{\ast }& =\left\{ \left( \lambda ,\chi (\lambda )\right) :\lambda
\in \Lambda \right\} .
\end{align*}%
Moreover,  $w^{cs}$ and $w^{cu}$ are Lipschitz
with constants $L$, and $\chi $ is Lipschitz with the constant $\sqrt{2}L/\sqrt{1-L^{2}}$. The manifolds $W^{cs}$ and $W^{cu}$ intersect
transversally on $\Lambda ^{\ast }$, and $W^{cs}\cap W^{cu}=\Lambda ^{\ast }$.

The manifolds $W^{cs}$ and $W^{cu}$ are foliated by invariant fibers $W_{z}^{s}$ and $W_{z}^{u}$, which are graphs of $C^{k}$ functions
\begin{eqnarray*}
w_{z}^{s} &:&\overline{B}_{s}(R)\rightarrow \Lambda \times \overline{B}_{u}(R), \\
w_{z}^{u} &:&\overline{B}_{u}(R)\rightarrow \Lambda \times \overline{B}_{s}(R),
\end{eqnarray*}%
meaning that%
\begin{eqnarray*}
W_{z}^{s} &=&\left\{ \left( w_{z}^{s}\left( y\right) ,y\right) :y\in
\overline{B}_{s}(R)\right\} , \\
W_{z}^{u} &=&\left\{ \left( \pi _{\lambda }w_{z}^{u}\left( x\right) ,x,\pi
_{y}w_{z}^{u}\left( x\right) \right) :x\in \overline{B}_{u}(R)\right\} .
\end{eqnarray*}%
The functions $w_{z}^{s}$ and $w_{z}^{u}$ are Lipschitz with constants $1/L$. Moreover,
\begin{eqnarray*}
W_{z}^{s} &=&\left\{p\in D: \Phi(t,p)\in D\textnormal{ for all }t\geq 0,\textnormal{ and }
\exists T_{0}\geq 0,\exists C>0\textnormal{ (which can depend on p)}\right. \\
&&\left. \textnormal{s.t. for }t\geq T_{0}, \left\Vert \Phi(t,p)-\Phi(t,z)\right\Vert \leq C\exp(t \overrightarrow{\mu_{s,1}})\right\} ,
\end{eqnarray*}%
and
\begin{eqnarray*}
W_{z}^{u} &=&\left\{p\in W^{cu}: \exists T_{0}\leq0,\exists C>0\textnormal{ (which can depend on p)}\right. \\
 && \left. \textnormal{s.t. for }t\leq T_{0},  \left\Vert \Phi(t,p)-\Phi(t,z)\right\Vert \leq C \exp (t \overrightarrow{\xi _{u,1,P}})\right\}.
\end{eqnarray*}

For $j=1,\dots,k$ the derivatives of the functions $w^{cu}$, $w^{cs}$, $w^{u,s}_z$ of order $j$
are bounded by some constants $C_j=O(\sup_{s=1,\dots,j} \sup_{z \in D} \|D^s f\|)$ and  constants $L$, $\overrightarrow{\xi}$ and $\overrightarrow{\mu}$.
\end{theorem}

\section{Comparison with other works}
\label{sec:comparison}
\subsection{Comparison with \cite{BK95,BK96}}

In \cite{BK95} we have the following formulas for $Q_0$ and $q_0$ (denoted there by $K$) are given (see \cite[page~63, Theorem 3.20 in Chapter II Sec. 2.3]{BK95}
and (\cite[II.3.2 page 48]{BK95})
\begin{eqnarray}
  Q_0(k)&=&\left[\frac{\lambda_{\max}}{\lambda_{min}} + k \left( \frac{\mu_{max}}{\lambda_{min}}+1\right)\right]+\left[\frac{\mu_{\max}}{\mu_{min}} + k\left( \frac{\lambda_{\max}}{\mu_{min}}+1\right)\right] +2,  \label{eq:Q0-BK95}\\
  q_0 (k)&=&Q_0(k)+k.\nonumber
\end{eqnarray}
It is clear that the value of $Q_0$ obtained by us (i.e., (\ref{eq:Q0-our})) is smaller  than that of (\ref{eq:Q0-BK95}). The difference seems to be equal to two. Our value of $q_0$ is also better.

The paper \cite{BK96}, inspired by \cite{BLW}, is an improvement of results from \cite{BK95}, as the authors consider the preservation of volume, symplectic form or contact structure.  In terms of $Q_0$ and $q_0$ the results appear to be  better than the ones from \cite{BK95}, however the case of a general vector field is not considered.
For the Hamiltonian vector field they obtain (see \cite[Theorem 8.1]{BK96} ) the following inequalities (we use our notation)
\begin{equation}
Q_0(k)+1= q_0(k) > \frac{2 \lambda_{max}(k+1)}{\lambda_{min}} + 2.  \label{eq:Q0Ham-B96}
\end{equation}
Since for Hamiltonian systems $\lambda_{max}=\mu_{max}$ and  $\lambda_{min}=\mu_{min}$ we see that (\ref{eq:Q0-BK95}) gives
\begin{equation*}
  Q_0(k) \geq \frac{2 \lambda_{max}(k+1)}{\lambda_{min}} + 2k,
\end{equation*}
hence (\ref{eq:Q0Ham-B96}) is an improvement for Hamiltonian systems.   As a result their numbers for Hamiltonian systems are better than ours.

Let us also stress that in \cite{BK95,BK96} the authors show also how some resonant terms can be removed by finitely smooth  coordinate change, which we are not discussing
in our work. In our paper~\cite{DZ} (and also in \cite{BLW}) we deal with the question of removing the remainder, only.

\subsection*{Comparison with \cite{BLW}}
Theorem 1.2 in \cite{BLW} gives (we use our notation) the following.
Let
\begin{eqnarray*}
  A= \frac{\lambda_{min}}{\mu_{max}} \frac{\mu_{min}}{\lambda_{max}+\mu_{min}}
\end{eqnarray*}
and let $B$ for the general vector field be given by
\begin{equation*}
  B_{no-struct}=\frac{\mu_{max}^2 + \mu_{min}(\lambda_{max}-\lambda_{min})}{\mu_{max}(\lambda_{max}+\mu_{min})}
\end{equation*}
and for Hamiltonian vector fields
\begin{equation*}
  B_{ham}=1 - 2 A.
\end{equation*}

 Then the relation between $k$ and $Q$ claimed by the authors is
 \begin{equation*}
 1 \leq k \leq Q A - B.
 \end{equation*}
 Therefore
 \begin{equation*}
    Q_0(k) \geq \frac{k + B}{A}.
 \end{equation*}
 Moreover,
 \begin{eqnarray*}
  q_0 \geq Q_0+2.
\end{eqnarray*}
It appears that the above bounds are considerably better than the ones obtained by us or in \cite{BK95,BK96}. For example when $\mu_{min}=\mu_{max}=\lambda_{min}=\lambda_{max}$
using the above formulas in the general vector field case one obtains $A=B=\frac{1}{2}$, hence
\begin{eqnarray*}
  Q_0(k) = 2k+1, \qquad q_0 = 2k+3.
\end{eqnarray*}

While our formulas give (see Section~\ref{subsec:Q0-spec-system})
\begin{eqnarray*}
  Q_0(k)=4k+2.
\end{eqnarray*}

Therefore the technique from \cite{BLW} appears to produce better constants for the pure saddle case. We have found it difficult at adopt it to the case when center 
directions are present. 




\end{document}